\documentclass[a4paper,10pt]{article}

\usepackage{amsmath,amsfonts,amssymb,amsthm}
\usepackage{graphics,epsfig,color,enumerate}
\usepackage{setspace,mathtools,enumitem,dsfont,verbatim}
\usepackage{textcomp}

\usepackage{geometry}
\usepackage[utf8]{inputenc}
\usepackage{lmodern}
\usepackage[utf8]{inputenc}
\usepackage[T1]{fontenc}
\usepackage{bbm}
\usepackage{xcolor}
\usepackage{mathtools}
\usepackage{bm}

\usepackage{todonotes}
\usepackage[english]{babel}
\usepackage[utf8]{inputenc}
\usepackage[T1]{fontenc}
\usepackage{csquotes}

\usepackage[nocompress]{cite}
\usepackage[hidelinks]{hyperref}

\usepackage{subcaption}
\usepackage{cases}
\usepackage{placeins}
\usepackage{mathrsfs}

\usepackage{letltxmacro}
\newcommand{\reffont}{\rm} 
\AtBeginDocument{%
\LetLtxMacro\oldref\ref
\renewcommand{\ref}[1]{%
{\reffont\oldref{#1}}}
}

\makeatletter
\newcommand{\customitem}[1]{%
\item[\rm#1]\protected@edef\@currentlabel{#1}%
}
\makeatother

\makeatletter
\def\@captype{figure}
\@addtoreset{equation}{section}
\makeatother

\newcounter{extralabel}[section]

\newtheoremstyle{fancyremark}{\topsep}{\topsep}{\rm}{}{\bfseries}{.}{ }{}
\newtheorem{ittheorem}{Theorem}
\newtheorem{itlemma}{Lemma}
\newtheorem{itproposition}{Proposition}
\newtheorem{itcorollary}{Corollary}
\newtheorem{italgorithm}{Algorithm}
\newtheorem{itassumption}{Assumption}
\newtheorem{itconstruction}{Construction}
\newtheorem{itobservation}{Observation}
\newtheorem{itconjecture}{Conjecture}
\newtheorem{itexample}{Example}

\theoremstyle{plain}
\newtheorem{itcondition}{Condition}
 
\theoremstyle{fancyremark}
\newtheorem{itremark}{Remark}
 
\theoremstyle{definition}
\newtheorem{itdefinition}{Definition}

 \newenvironment{theorem}{\addtocounter{extralabel}{1}
 \begin{ittheorem}}{\end{ittheorem}}

 \newenvironment{lemma}{\addtocounter{extralabel}{1}
 \begin{itlemma}}{\end{itlemma}}

 \newenvironment{proposition}{\addtocounter{extralabel}{1}
 \begin{itproposition}}{\end{itproposition}}

 \newenvironment{definition}{\addtocounter{extralabel}{1}
 \begin{itdefinition}}{\end{itdefinition}}

 \newenvironment{corollary}{\addtocounter{extralabel}{1}
 \begin{itcorollary}}{\end{itcorollary}}

 \newenvironment{example}{\addtocounter{extralabel}{1}
 \begin{itexample}}{\end{itexample}}

\newcommand\xqed[1]{
\leavevmode\unskip\penalty9999 \hbox{}\nobreak\hfill
\quad\hbox{#1}}
\newcommand\RemarkEnd{\xqed{$\blacksquare$}}
\newenvironment{remark}{\addtocounter{extralabel}{1}

\begin{itremark}}{\RemarkEnd\end{itremark}}

\usepackage{xpatch}
\makeatletter
\AtBeginDocument{\xpatchcmd{\@thm}{\thm@headpunct{.}}{\thm@headpunct{}}{}{}}
\makeatother
\makeatletter
\patchcmd{\thmhead@plain}%
{\thmnote{ {\the\thm@notefont(#3)}}}%  original form
{\thmnote{ {\the\thm@notefont\bf[#3]}}}%  new form
{}{}
\let\thmhead\thmhead@plain
\makeatother
\newcommand{\intv}[2][1]{\{#1, \ldots, #2\}}
\newcommand{\ER}{Erdős-Rényi }
\newcommand{\distr}[2]{\mu^{#1}_{#2}}

\newcommand{\cycle}[2]{\gamma_{#2}(#1)}
\newcommand{\MC}[1]{\gamma_{\rm max}(#1)}

\newcommand{\graphcycle}[2]{\mathscr{C}_{#1}(#2)}
\newcommand{\maxgraph}[2]{\mathscr{C}^{\rm #2}_{\rm max}(#1)}
\newcommand{\secgraph}[2]{\mathscr{C}^{\rm #2}_{\rm sec}(#1)}
\newcommand{\conncomp}[1]{\bm{\mathcal{CC}}(#1)}

\newcommand{\xcycle}[3]{\mathcal{C}^{#3}_{#1}(#2)}
\newcommand{\Unif}[1]{\mathsf{Unif}(#1)}
\newcommand\restrict[1]{\raisebox{-.5ex}{$|$}_{#1}}
\newcommand{\norm}[1]{\left\lVert#1\right\rVert}
\newcommand{\vertex}{\mathcal{V}}
\newcommand{\edgeset}{\mathcal{E}}
\newcommand{\startvertex}{v_0}

\newcommand{\typicalevent}[2]{\Omega^{\rm #1}_n(#2)}

\newcommand{\TVD}{d_{\mathrm{TV}}}

\newcommand{\AG}[1]{A_{#1}}
\DeclareMathOperator{\supp}{supp}
\newcommand{\fsupp}[1]{\mathrm{supp}\left(#1\right)}
\newcommand{\deq}{\stackrel{d}{=}}
\newcommand{\convdist}{\stackrel{d}{\rightarrow}}
\renewcommand{\P}{\mathbb{P}}
\newcommand{\convprob}{\stackrel{\P}{\rightarrow}}

\newcommand{\E}{\mathbb{E}}
\newcommand{\N}{\mathbb{N}}
\newcommand{\1}{\mathbbm{1}}
\newcommand{\dd}{\mathrm{d}}
\newcommand{\ee}{\mathrm{e}}
\newcommand{\tddalt}{T^\Downarrow_{n,\startvertex}}
\newcommand{\Mtddalt}{T^\Downarrow_{n,\startvertex}}
\newcommand{\tlocalmix}[1]{T^{{\rm\sss(LM)},{#1}}_{n,\startvertex}}
\newcommand{\cond}{\,\mid}
\newcommand{\dnTVD}[1]{\mathcal{D}_n^{\startvertex}(#1)}
 
\newcommand{\sss}{\scriptscriptstyle}

\newcommand{\ISRW}{n,v_0}
\newcommand{\ISRWdistr}[1]{\distr{\ISRW}{#1}}

\newcommand{\vep}{\varepsilon}
\newcommand{\prob}{{\mathbb P}}
\newcommand{\eqn}[1]{\begin{equation}#1\end{equation}}

\newcommand{\op}{o_{\sss \prob}}
\newcommand{\convp}{\stackrel{\sss\P}{\longrightarrow}}
\newcommand{\convd}{\stackrel{\sss d}{\longrightarrow}}

\newcommand{\cmpl}{{\rm{c}}}

\newcommand{\Xievent}{\Xi}

\title{Mixing of fast random walks on\\ 
dynamic random permutations}
\author{}

\author{

Luca Avena
\footnotemark[1]
\\

Remco van der Hofstad
\footnotemark[2]
\\

Frank den Hollander
\footnotemark[3]
\\

Oliver Nagy
\footnotemark[3]
}

\footnotetext[1]{
Dipartimento di Matematica e Informatica `Ulisse Dini',
Università degli Studi di Firenze,
Viale Morgagni 67/a, 50134, Firenze, Italy
}

\footnotetext[2]{
Department of Mathematics and Computer Science, Eindhoven University of Technology, 
P.O.\ Box 513, 5600 MB Eindhoven, The Netherlands
}

\footnotetext[3]{
Mathematical Institute, Leiden University, 
P.O.\ Box 9512, 2300 RA Leiden, The Netherlands
}

\date{February~29, 2024}

\allowdisplaybreaks

\begin{document}

\maketitle

\begin{abstract}
We analyse the mixing profile of a random walk on a dynamic random permutation, focusing on the regime where the walk evolves much faster than the permutation. Two types of dynamics generated by random transpositions are considered: one allows for coagulation of permutation cycles only, the other allows for both coagulation and fragmentation. We show that for both types, after scaling time by the length of the permutation and letting this length tend to infinity, the total variation distance between the current distribution and the uniform distribution converges to a limit process that drops down in a single jump. This jump is similar to a one-sided cut-off, occurs after a random time whose law we identify, and goes from the value $1$ to a value that is a strictly decreasing and deterministic function of the time of the jump, related to the size of the largest component in \ER random graphs. After the jump, the total variation distance follows this function down to $0$.

\vspace{0.2cm}
\\
\noindent
\emph{Acknowledgment.}
The work in this paper was supported by the Netherlands Organisation for Scientific Research (NWO) through Gravitation-grant NETWORKS-024.002.003. The authors thank Nathana\"el Berestycki and Sam Olesker-Taylor for useful discussions. The paper was finalized while the authors were participating in the ICTS-NETWORKS workshop “Challenges in Networks” (ICTS/NETWORKS2024/01) at the International Centre for Theoretical Sciences (ICTS) in Bengaluru, India. The authors are thankful to ICTS for their hospitality.
\normalsize

\end{abstract}

\newpage 

\small
\tableofcontents
\normalsize

\newpage

%%%%%%%%%%% SECTION 1 %%%%%%%%%%%%%%%%%%%

\section{Introduction and main results}

%%%

\subsection{Target}

The goal of this paper is to identify the mixing profile of a fast random walk on a \emph{dynamic random permutation}, where fast means that the random walk instantly achieves local equilibrium, i.e., fully mixes on the cycle of the permutation it sits on before the next change in the permutation occurs. The focus is on two types of dynamics for the permutation, both starting from the identity permutation and consisting of successive applications of \emph{random transpositions}. The first type -- called \emph{coagulative dynamics} -- imposes the constraint that transpositions leading to fragmentation of a permutation cycle are ignored. The second type -- called \emph{coagulative-fragmentative dynamics} -- does not impose this constraint. A major feature of dynamic random permutations is that they represent a \emph{disconnected} geometry, which marks a departure from the setting that was considered in earlier work.

\begin{figure}[htbp]
\centering
\begin{minipage}{.48\textwidth}
\centering
\begin{tikzpicture}
\node (img) {\includegraphics[width=\columnwidth-2.5cm]{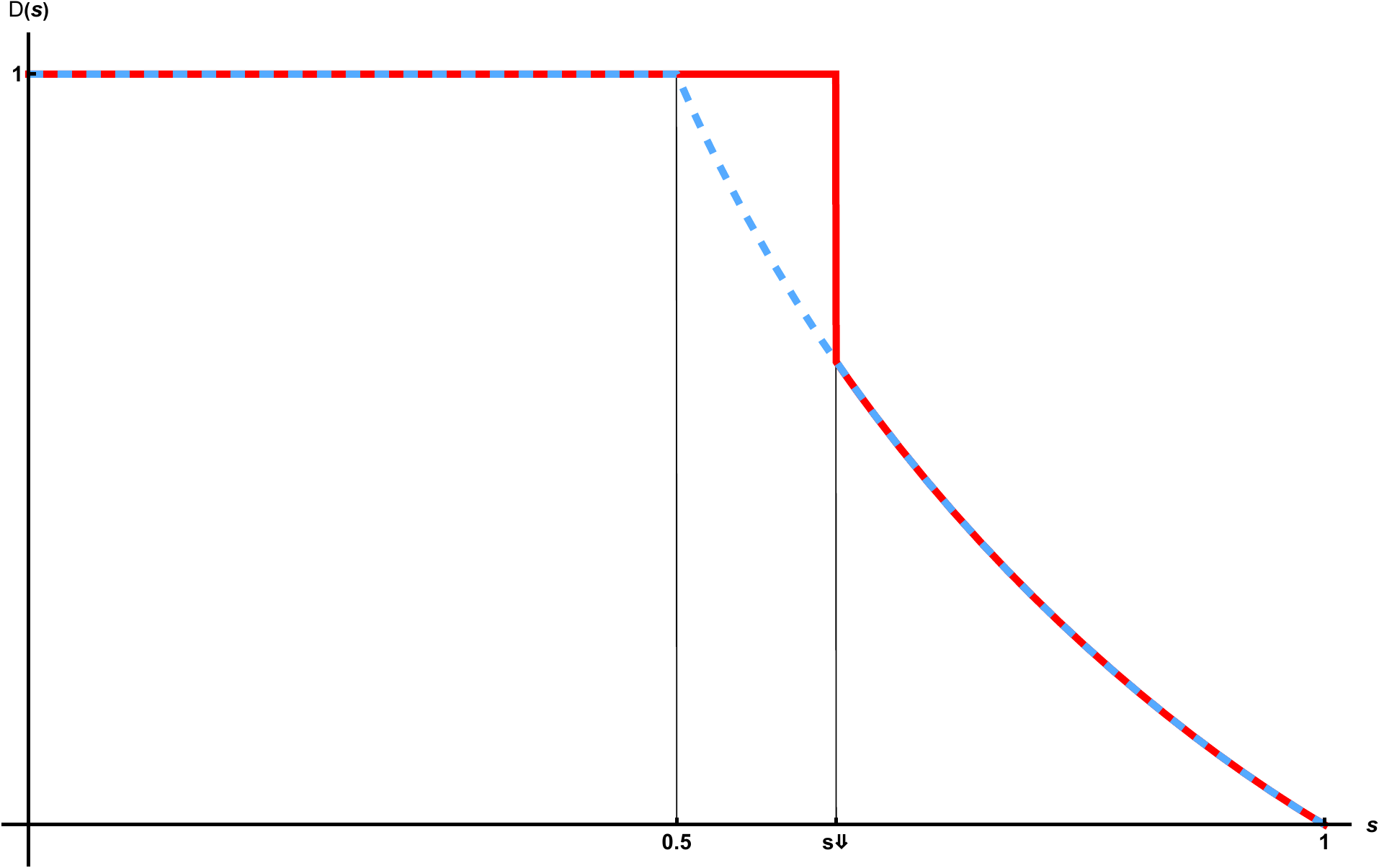}};
\node[below=of img, node distance=0cm, yshift=1cm,font={\tiny\color{black}}] {scaled time};
\node[left=of img, node distance=0cm, rotate=90, anchor=center,yshift=-0.7cm,font={\tiny\color{black}}] {TVD};
\end{tikzpicture}
\captionof{figure}{\small The red curve is a typical evolution of the total variation distance for an infinitely fast random walk on a \emph{coagulative dynamic permutation}. The blue curve is a plot of the deterministic function of the scaled time to which the total variation distance drops at a random time and subsequently sticks to.}
\label{fig:intro}
\end{minipage}%
%%%
\hfill
%%%
\begin{minipage}{.49\textwidth}
\centering
\includegraphics[width=\columnwidth-2cm]{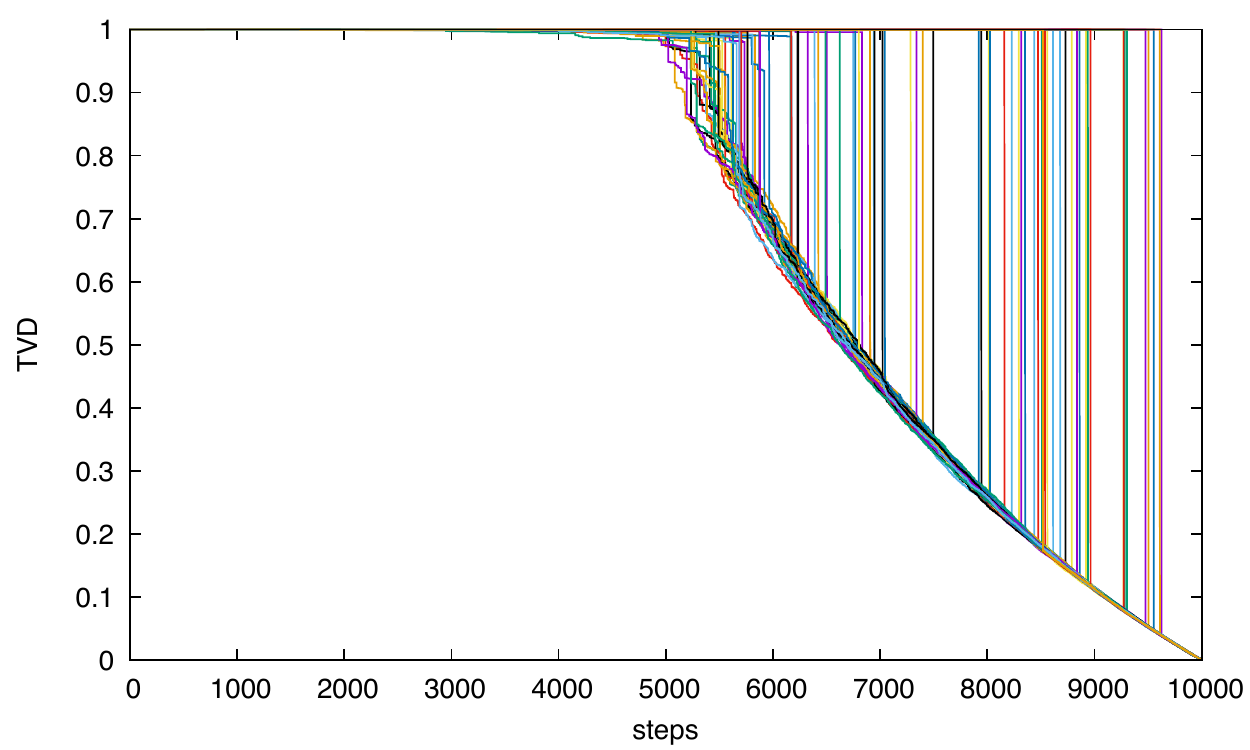}
\captionof{figure}{\small Simulations of the evolution of the total variation distance for $10^2$ different realisations of a \emph{coagulative dynamic permutation} of $10^4$ elements and an infinitely fast random walk on top. Each simulation run corresponds to a single coloured curve.}
\label{fig:sim-coag}
\end{minipage}
\begin{minipage}{.49\textwidth}
\centering
\begin{tikzpicture}
\node (img)  {\includegraphics[width=\columnwidth-2.5cm]{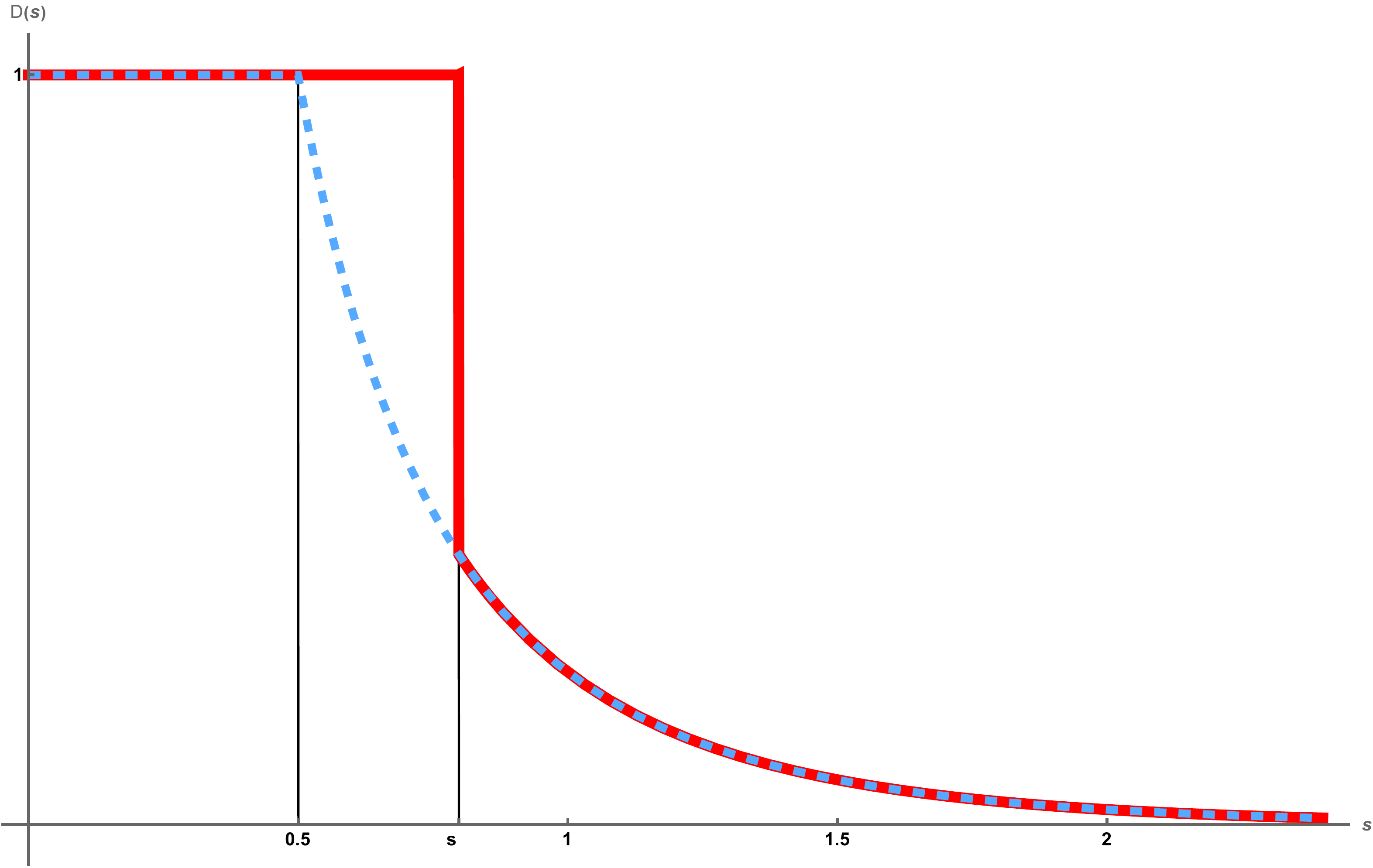}};
\node[below=of img, node distance=0cm, yshift=1cm,font={\tiny\color{black}}] {scaled time};
\node[left=of img, node distance=0cm, rotate=90, anchor=center,yshift=-0.7cm,font={\tiny\color{black}}] {TVD};
\end{tikzpicture}
\captionof{figure}{\small The same as Fig.~\ref{fig:intro} for a \emph{coagulative-fragmentative dynamic permutation}.}
\label{fig:intro-coagfrag}
\end{minipage}
%%%
\hfill
%%%
\begin{minipage}{.49\textwidth}
\centering
\includegraphics[width=\columnwidth-2cm]{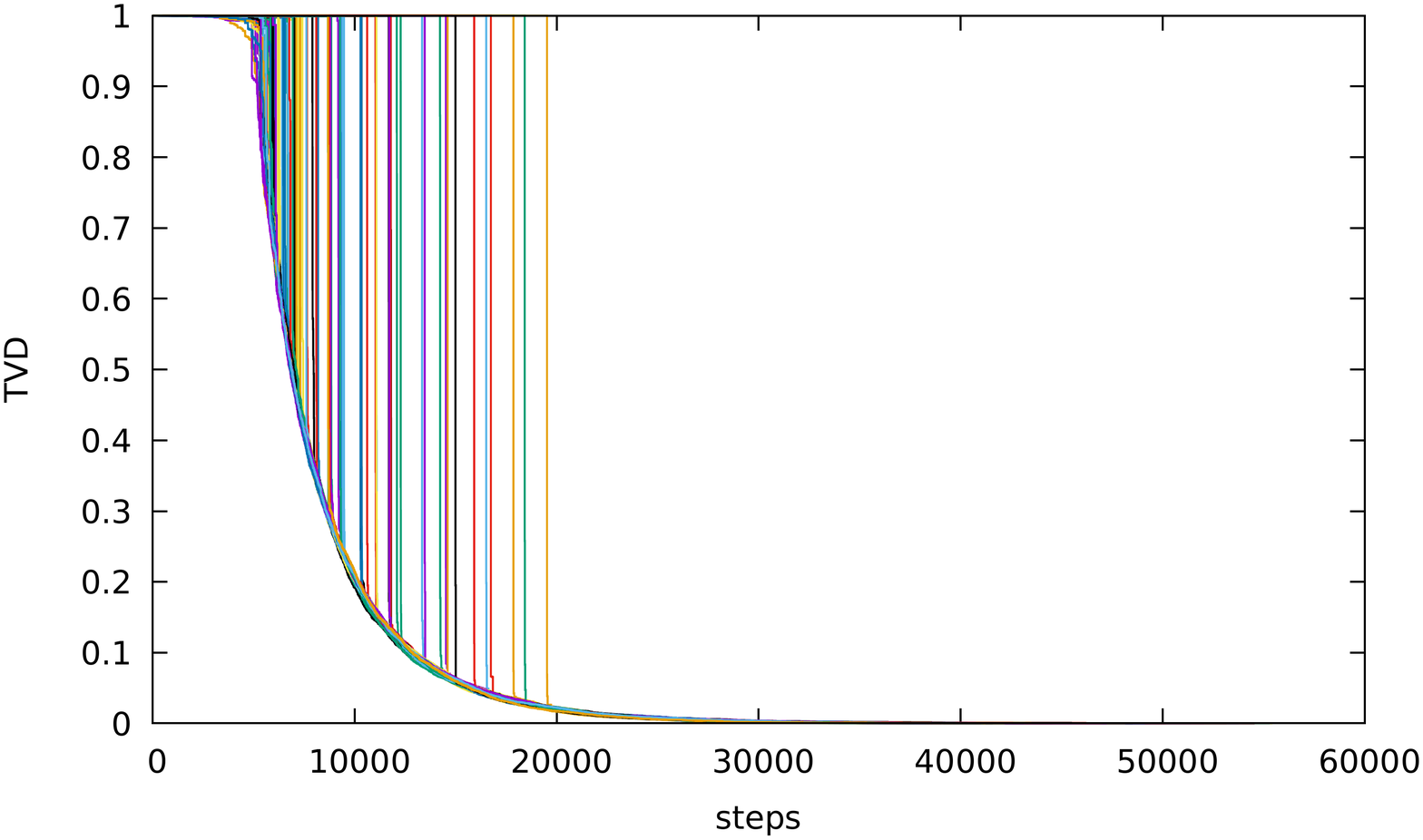}
\captionof{figure}{\small The same as Fig.~\ref{fig:sim-coag} for a \emph{coagulative-fragmentative permutation}.}
\label{fig:sim-coagfrag}
\end{minipage}
\end{figure}

We show that for both dynamics, after scaling time by the length of the permutation and letting this length tend to infinity, the total variation distance between the current distribution and the uniform distribution converges to a limit process that makes a \emph{single jump} down from the value $1$ to a value on a deterministic curve and subsequently follows this curve on its way down to $0$. The aforementioned curve is strictly decreasing in time and is related to the size of the largest component in the \ER random graph. The jump down to this curve, which is similar to a \emph{one-sided cut-off}, occurs after a \emph{random time} whose law we identify. The law of the drop-down time and the function describing the deterministic curve are different for the two types of dynamics. Visual representations of the mixing profiles are given in Figs.~\ref{fig:intro} and \ref{fig:intro-coagfrag}, while simulations are shown in Figs.~\ref{fig:sim-coag} and \ref{fig:sim-coagfrag}.

%%%%%%%%%%%%%%%%%%%%%%%%%%%%%%%%%%%%%%%%%%%
\FloatBarrier
%%%%%%%%%%%%%%%%%%%%%%%%%%%%%%%%%%%%%%%%%%%%

The remainder of this section is organised as follows. Section~\ref{sec:intro:background} provides background and recalls earlier work. Section~\ref{sec:introdef} fixes the setting and introduces relevant definitions and notations. Section~\ref{sec:intro:ER} lists some preliminaries for \ER random graphs that are needed along the way. Section~\ref{ss.agp} introduces a graph process associated with the dynamics that serves as a tool for analysing the dynamics. Section~\ref{sec:intro:results} contains two main theorems, one for each type of dynamics, describing the evolution of the total variation distance between the current distribution of the random walk and its equilibrium distribution, which is the uniform distribution on $[n]=\intv{n}$. Section~\ref{sec:intro:discussion} discusses the importance of the main theorems, places them in their proper context, and provides an outline of the remainder of the paper.  

%%%%

\subsection{Background and earlier work}
\label{sec:intro:background}

While over the past years random walks on \emph{static} random graphs have received a lot of attention, and the scaling properties of quantities like mixing times, cover times and metastable crossover times have been identified, much less is known about random walks on \emph{dynamic} random graphs. In the static setting, a two-sided cut-off on scale $\log n$  has been established for a general class of \emph{undirected} sparse graphs with good expansion properties \cite{BHS2017,BLPS2018,LS2010}. Similar results have been obtained for \emph{directed} sparse graphs \cite{BCS2018,BCS2019} and for graphs with a \emph{community structure} \cite{BH2018}.

 In the dynamic setting, predominantly the focus has been on dynamic percolation, \ER random graphs with edges switching on and off randomly, and configuration models with random rewiring of edges. Both directed and undirected graphs have been considered, as well as backtracking and non-backtracking random walks. In \cite{PSS2015,PSS2017,HS2019} random walks on \emph{dynamic percolation} clusters on a $d$-dimensional \emph{discrete torus} were considered. Mixing times were identified for several parameter regimes controlling the rate of the random walk and the rate of the random graph dynamics. Similar results were obtained for dynamic percolation on the \emph{complete graph} \cite{ST2018,SZ2019}. Some further advances were achieved in \cite{AKL2018}, where general bounds on mixing times, hitting times, cover times and return times were derived for certain classes of dynamic random graphs under appropriate \emph{expansion assumptions}. Non-backtracking random walks on configuration models that with high probability are connected were studied in a series of papers \cite{AGHH2018,AGHH20182,AGHHN2020}, which culminated in a general framework for studying mixing times of non-backtracking random walks on dynamic random graphs subject to mild regularity conditions. Mixing of random walks on \emph{directed} configuration models was treated in \cite{CQ2019}.

Random permutations generated by random transpositions have attracted plenty of interest as well. An important starting point is \cite{DS1981}, where a \emph{cut-off} in the total variation distance was established after the application of $\frac12 n\log n + O(n)$ random transpositions. The \emph{sharp} constant in front of the leading-order term was achieved with the help of \emph{representation theory} for the symmetric group. This paper led to a flurry of follow-up work, of which we mention \cite{S2005}, where the structure of large cycles of an evolving random permutation was studied. Similar results were obtained in \cite{BD2006}, including sharp control on the number of observed fragmentations. An important aspect of both \cite{S2005} and \cite{BD2006} is the representation of an evolving random permutation, starting from the identity permutation, in terms of a random graph process that can be studied by using the theory of random graphs. For the \emph{coagulative-fragmentative} dynamics considered in the present paper, also called \emph{transposition dynamics}, this graph process representation yields a graph-growth model that at every step adds an edge drawn uniformly at random. This graph-growth model is closely related to the standard \enquote{combinatorial} \ER model, whose study is by now a classical topic in the theory of random graphs (see, for example, \cite{FK2016} or \cite{vdH2016}). Yet another important feature of \cite{S2005} is the introduction of \emph{Schramm's coupling} as a tool to study the cycle structure of evolving random permutations. In a follow-up article \cite{BSZ2011}, a modified version of this coupling is used to study the mixing of dynamic permutations endowed with a more general dynamics, of which the transposition dynamics is a special case. We also mention \cite{BKLM2019}, which contains a detailed account of Schramm's coupling. The works cited above  each highlight one particular facet of the random transposition model, but close relatives have been studied extensively under different names: \emph{mean-field T\'oth model} \cite{B1993}, the \emph{interchange process} on the complete graph (see \cite{AK2013,BKLM2019,HS2021} and references therein), or \emph{multi-urn Bernoulli–Laplace diffusion models} \cite{S2021}, where our setting corresponds to a particular choice of the model parameters.

The \emph{coagulative} dynamics considered in the present paper can be recast, in the spirit of \cite{S2005}, as a graph-valued random process that starts with an empty graph on $n$ vertices and describes a forest that progressively merges into a spanning tree on $n$ vertices through the addition of edges that do not create a cycle. The study of this process and its close relatives has a somewhat twisted history. It is similar to the \emph{standard additive coalescent} (see \cite{A1999} for an overview), but it is also interesting in its own right (see \cite{A1990,LP1992}). Finally, there is a wealth of results on \emph{minimal spanning trees} and \emph{Kruskal's algorithm}, which is another closely related process. In particular, we mention \cite{JS2019}, since this work implies some facts that we list in Section~\ref{sec:CDP}. We derive these facts independently, using different techniques in a different setting. 

%%%
\subsection{Setting, definitions and notation}
\label{sec:introdef}

For $n \in \N$, let $S_n$ denote the set of permutations of $[n]$, i.e., bijections from $[n]$ to itself. Recall that $S_n$ endowed with the operation of permutation composition $\circ$ forms a group. Write $\cycle{\pi}{v}$ to denote the cycle of the permutation $\pi$ that contains the element $v$. 

\begin{definition}[Dynamic permutation]
A sequence of permutations of $[n]$, denoted by $\Pi_n = (\Pi_n(t))_{t=0}^{t^{\rm max}}$ with $t^{\rm max} \in \N_0 \cup \{\infty\}$, is called a dynamic permutation. \hfill$\spadesuit$
\end{definition}

\begin{example}[Transpositions may fragment cycles or coagulate cycles]
\label{ex:coalfrag}
{\rm Pick $n=7$ and consider the permutation
$$
\left(\begin{array}{lllllll} 1 &2 &3 &4 &5 &6 &7\\ 2 &3 &4 &5 &6 &7 &1 \end{array}\right) 
\quad \textrm{with cycle structure} \quad (1,2,3,4,5,6,7).
$$
The transposition $(1,5)$ turns this into the permutation
$$
\left(\begin{array}{lllllll} 1 &2 &3 &4 &5 &6 &7\\  6 &3 &4 &5 &2 &7 &1 \end{array}\right)  
\quad \textrm{with cycle structure} \quad (1,6,7)(2,3,4,5).
$$
Another application of the same transposition acts in reverse. Note that $S_n$ is a \emph{non-commutative group} for any $n \geq 3$.}\hfill$\blacksquare$
\end{example}

We consider two types of dynamic permutations:

\begin{definition}[Coagulative dynamic permutation]
\label{def:cdp}
$\Pi_n = (\Pi_n(t))_{t=0}^{n-1}$ is called a \emph{coagulative dynamic permutation} (CDP) when $\Pi_n(0) = \mathrm{Id}$ (i.e., the identity permutation) and
\begin{equation}
\Pi_n(t) = \Pi_n(t-1) \circ (a,b), \qquad t \in [n-1],
\end{equation}
where, for each $t \in [n-1]$, $(a,b)$ is a random transposition sampled uniformly at random from the set of all transpositions of $[n]$ that satisfy the \emph{constraint}
\begin{equation}
\label{eq:coag-restriction}
\cycle{\Pi_n(t-1)}{a} \neq \cycle{\Pi_n(t-1)}{b}. 
\end{equation}
The latter guarantees that no cycle of $\Pi_n(t-1)$ is fragmented by the transposition $(a,b)$.  
\hfill$\spadesuit$
\end{definition}

\begin{definition}[Coagulative-fragmentative dynamic permutation]
\label{def:cfdp}
$\Pi_n = (\Pi_n(t))_{t=0}^{\infty}$ is called a \emph{coagulative-fragmentative dynamic permutation} (CFDP) when the same holds as in Definition~\ref{def:cdp}, but without the constraint in \eqref{eq:coag-restriction}. \hfill$\spadesuit$
\end{definition}

\begin{remark}[Time horizon for dynamic permutations and cycle structure]
Since CDP starts from the identity permutation, it becomes a permutation with a single cycle after exactly $n-1$ steps. Once this happens, there is no permutation that satisfies \eqref{eq:coag-restriction} and the dynamics is trapped. CFDP has no traps and can evolve forever. The structure of cycles is random and their sizes, scaled by $1/n$, converge in distribution to the Poisson-Dirichlet distribution with parameter 1 (see e.g.\ \cite[Theorem 1.1]{S2005} for a precise statement).
\end{remark}

Our aim is to study mixing of fast random walks on both CDP and CFDP. To simplify our analysis, we work with \emph{infinite-speed} random walks, as defined next:

\begin{definition}[Infinite-speed random walk on $\Pi_n$]
\label{def:ISRW}
Fix $\Pi_n$ and an element $\startvertex \in [n]$. Recall that $\cycle{\Pi_n(t)}{v}$ is the cycle of $\Pi_n(t)$ that contains $v$. Formally, the infinite-speed random walk (ISRW) starting from $\startvertex$ is a sequence of probability distributions $(\ISRWdistr{}(t))_{t\in\N_0}$ supported on $[n]$, with initial distribution at time $t=0$ given by
\begin{equation}
\ISRWdistr{}(0) = \left(\ISRWdistr{w}(0)\right)_{w\in[n]},
\end{equation}
where $\ISRWdistr{w}(0)$, the mass at $w \in [n]$ at time $t=0$, is given by
\begin{equation}
\ISRWdistr{w}(0) = 
\begin{cases}
\frac{1}{|\cycle{\Pi_n(0)}{w} |},
&w \in \cycle{\Pi_n(0)}{\startvertex},\\ 
0, 
&w \notin \cycle{\Pi_n(0)}{\startvertex},
\end{cases}
\end{equation}
and with distribution at a later time $t \in \N$ given by
\begin{equation}
\ISRWdistr{}(t) = \left(\ISRWdistr{w}(t)\right)_{w\in[n]},
\end{equation}
where
\begin{equation} 
\ISRWdistr{w}(t) = \frac{1}{|\cycle{\Pi_n(t)}{w}|} \sum_{u \in \cycle{\Pi_n(t)}{w}}\ISRWdistr{u}(t-1).
\end{equation}
Informally, ISRW \emph{spreads infinitely fast over the cycle in the permutation it resides on}. \hfill$\spadesuit$ 
\end{definition}

\begin{figure}[htbp]
\vspace{0.3cm}
\centering
\includegraphics[width=\textwidth]{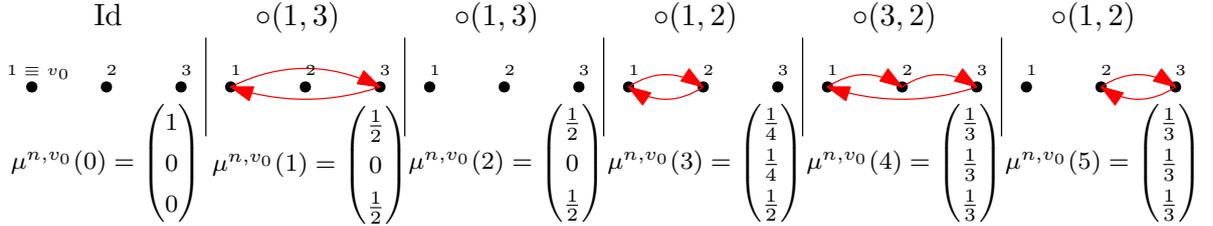}
\vspace{0.1cm}
\caption{\small Example of an evolution of an ISRW on top of a CFDP with three elements starting from the identity permutation. The first row shows the transpositions that generate the next permutation. The second row is a graphical representation of the cycles of this permutation. The third row shows the evolution of the ISRW distribution, given that it started from the element  $1$.}
\label{fig:exISRW}
\end{figure}

In Appendix~\ref{apx:limit} we show that the infinite-speed random walk arises as the limit of a standard random walk whose stepping rate relative to the rate of the permutation dynamics tends to infinity. Note that the evolution of the ISRW distribution is fully determined by the initial position of the random walk and the realisation of the dynamic permutation. See Figure~\ref{fig:exISRW} for an illustration.

%%%%%%%%%%%%%%%%%%%%%%%%%%%%%%%%%%%%%%%%%%%%%%%%%%%%%%%%%%%%
%%%%%%%%%%%%%%%%%%%%%%%%%%%%%%%%%%%%%%%%%%%%%%%%%%%%%%%%%%%%

%%%

\subsection{Preliminaries for \ER random graphs}
\label{sec:intro:ER}

The arguments in this paper frequently make use of results on the structure of \ER random graphs. This section provides what is needed to state the main theorems in Section~\ref{sec:intro:results}.

\begin{definition}[Standard \ER multi-graph process]
\label{def:ER1}
The \emph{standard \ER multi-graph process} on $n$ vertices is the discrete-time process $\left(G(n,t) \right)_{t=0}^{t_{\rm max}}$ constructed as follows:
\begin{enumerate}
\item 
$G(n,0)$ is the graph with $n$ vertices and no edges. 
\item 
At each time $t\in\N_0$, pick an edge $e_t$ uniformly at random from the $\binom{n}{2}$ possible edges, and let $G(n,t)$ be the graph obtained by adding $e_t$ to $G(n, t-1)$.
\end{enumerate}
Note that we do not allow for self-loops, but do allow for multiple edges. \hfill$\spadesuit$
\end{definition}

\begin{remark}[Versions and asymptotic equivalence]
There are versions of the \ER multi-graph process that differ in how edges are deployed and whether or not multiple edges and self-loops are allowed. With respect to monotone properties, notably the expected size of connected components, the \enquote{random growth} $G(n,t)$ model described in Definition~\ref{def:ER1} is asymptotically equivalent to the \enquote{combinatorial} model $G(n,M)$ with $M = t$ edges at times $t=O(n)$, which in turn is asymptotically equivalent to the \enquote{bond percolation} model $G(n,p)$ with $\smash{p = M {\binom{n}{2}}^{-1}}$. For details, see~\cite[Sections~1.1, 1.3]{FK2016}. Since we work on time scales of order $n$, we will use this asymptotic equivalence without further notice.
\end{remark}

Definition~\ref{def:ER1} allows for some natural modifications, of which one is important for the study of CDP:

\begin{definition}[Cycle-free \ER graph process]
\label{def:cfER}
The \emph{cycle-free \ER graph process} on $n$ vertices is the graph process $\smash{\left(G(n,t)\right)_{t=0}^{t_{\rm max}}}$ starting from the empty graph with $n$ vertices, such that at each time $t$ an edge is added that is chosen uniformly at random from the set of edges that do not create a cycle, a multi-edge or a self-loop. Thus, $G(n,t)$ is a forest for all $0 \leq t\leq t_{\rm max}$. \hfill$\spadesuit$
\end{definition}

To understand the typical evolution of CDP, we make use of two \emph{couplings}: one between CDP and cycle-free \ER graph processes, the other between cycle-free \ER graph processes and their standard counterparts. To explain how, we need to introduce three functions that describe key structural properties of these processes: 

\begin{definition}[Functions related to the structure of \ER random graphs]
\label{def:ERstruct}
$\mbox{}$\\
(1) Define $\zeta\colon\, [0, \infty) \to [0,1)$ as $\zeta(u) = 0$ for $u\in[0, \tfrac12]$ and as the unique positive solution of the equation $1 - \zeta(u) = \ee^{-2u\zeta(u)}$ for $u \in (\tfrac12,\infty)$. Note that $\zeta$ is non-decreasing and continuous on $[0,\infty)$, and analytic on $(\tfrac12,\infty)$.\\
(2) Define $\phi\colon\, [0, \infty) \to [0,1)$ as
\begin{equation}
\label{def:phi}
\phi(v) = \int_0^v \dd u\,[1-\zeta^2(u)], \qquad v \in [0,\infty).
\end{equation}
Note that $\phi$ is strictly increasing and continuous on $[0,\infty)$, and hence has a well-defined inverse $\phi^{-1}$. Furthermore, the function $\phi$ is properly normalised in the sense that $\phi(\infty) = 1$ (see Appendix~\ref{apx:norm}).\\ 
(3) Define $\eta\colon\, [0,1) \to [0,1)$ as  
\begin{equation}\label{eq:eta}
\eta(w) = \zeta(\phi^{-1}(w)), \qquad w\in[0,1).
\end{equation}
\hfill$\spadesuit$
\end{definition}

\begin{figure}[htbp]
\centering
\includegraphics[width=.28\textwidth]{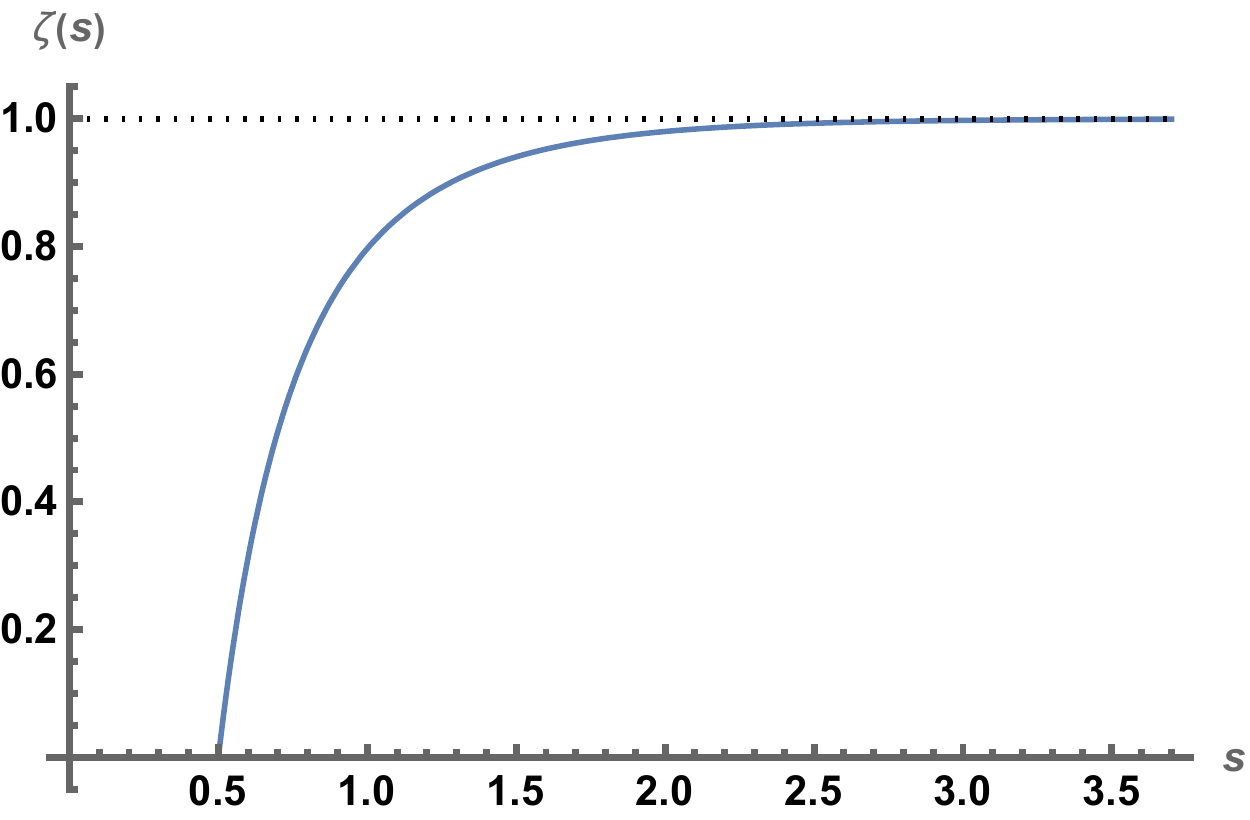}
\qquad
\includegraphics[width=.28\textwidth]{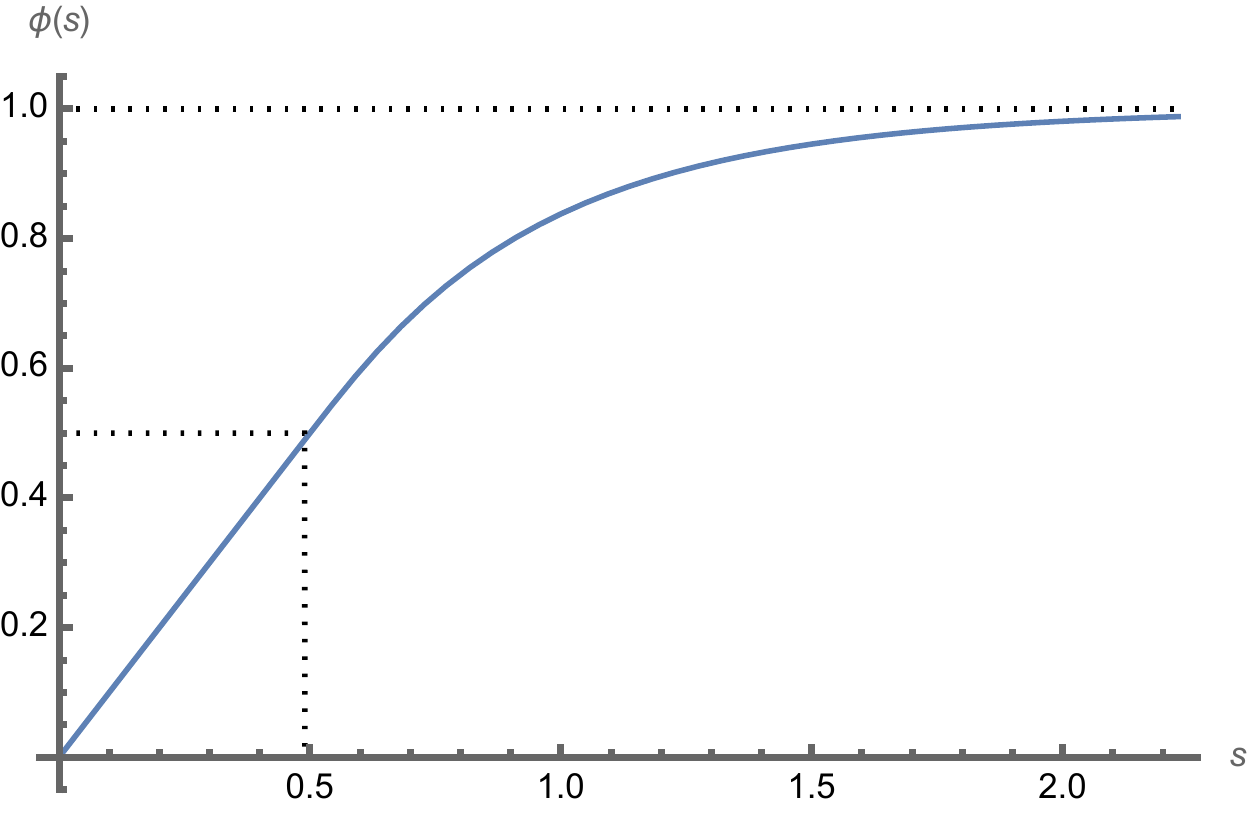}
\qquad
\includegraphics[width=.28\textwidth]{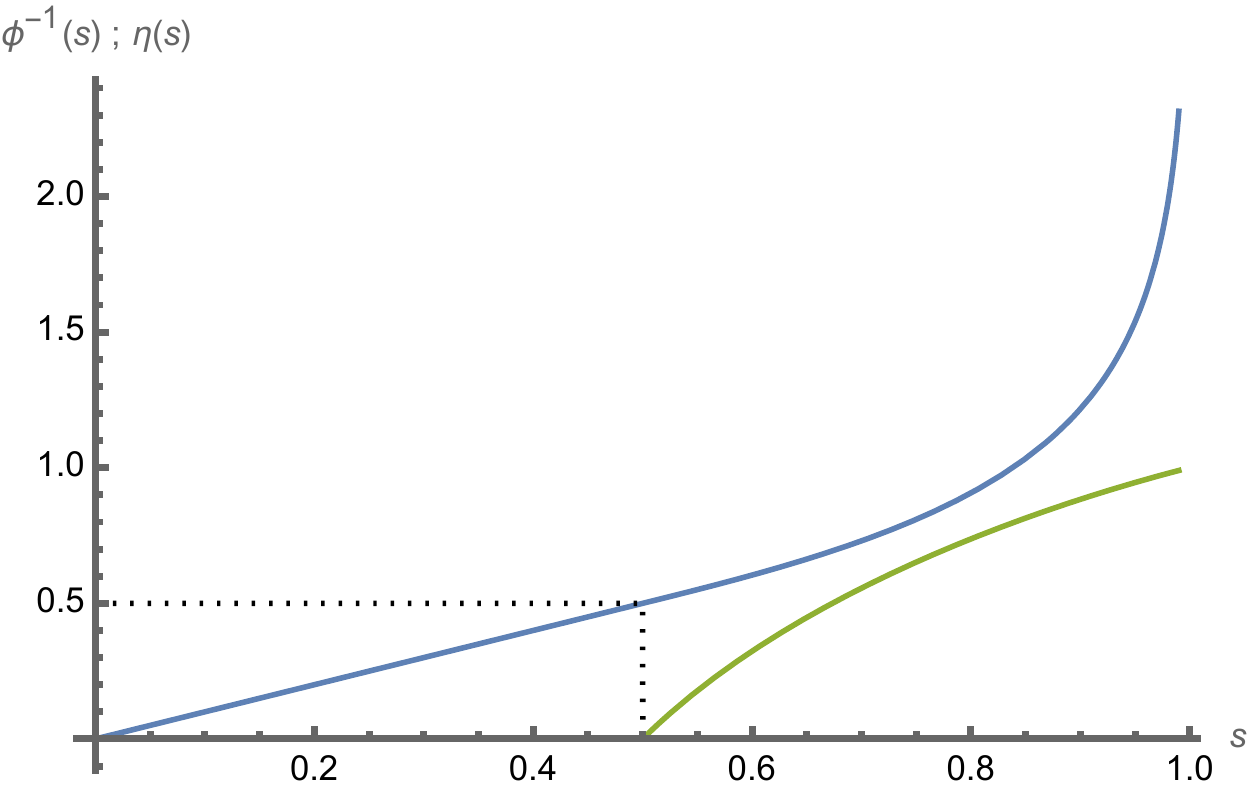}
\caption{\small Graphs of the functions introduced in Definition~\ref{def:ERstruct}: $\zeta$, $\phi$, respectively, $\phi^{-1}$ (upper curve), $\eta$ (lower curve).}
\label{fig:graphs}
\end{figure}

The functions defined in Definition~\ref{def:ERstruct} are illustrated in Fig.~\ref{fig:graphs} and have the following interpretation: 
\begin{itemize}
\item[(1)]
$\zeta(u)$ describes the expected size of the largest component of the \ER random graph at time $un$. For $u \in [0,\infty)$, denote by $|\mathscr{C}_{\rm max}^{\rm ER}(n,un)|$ the size of the largest connected component in the \ER random graph with $n$ vertices and $un$ edges, and $|\mathscr{C}_{\rm sec}^{\rm ER}(n,un)|$ the size of the second-largest connected component. Then, as $n\to\infty$, 
\begin{equation}
\label{eq:ERsizes}
\frac{|\mathscr{C}_{\rm max}^{\rm ER}(n,un)|}{n} \convprob \zeta(u), \qquad \frac{|\mathscr{C}_{\rm sec}^{\rm ER}(n,un)|}{n} \convprob 0.
\end{equation}
\item[(2)]
The function $\phi$ provides the link between the standard and the cycle-free \ER graph process (see Lemmas~\ref{lem:cc}--\ref{lem:et} below).
\item[(3)]
$\eta(u)$ is the analogue of $\zeta(u)$ for the cycle-free \ER graph process at time $un$, $u\in [0,1]$ (see Lemma~\ref{lem:CFmax} below).
\end{itemize}

\noindent
Note the change in behaviour of $\zeta,\phi,\phi^{-1},\eta$ at $\tfrac12$. Note that $\phi^{-1}$ blows up at $1$.

%%%%%%%%%%%%%% FIGURE %%%%%%%%%%%%%%%%%
%%%%%%%%%%%%%%%%%%%%%%%%%%%%%%%%%%%%%%%%%

%%%

\subsection{Associated graph process}
\label{ss.agp}

For any dynamic permutation starting from the identity permutation, define the associated graph process as follows:

\begin{definition}[Graph process associated with $\Pi_n$]
\label{def:AG}
Let $\Pi_n=(\Pi_n(t))_{t=0}^{t_{\rm max}}$ with $t_{\rm max} \in \N \cup \{ \infty \}$ be a dynamic permutation starting from the identity permutation. Construct the \emph{associated graph process}, denoted by $\AG{\Pi_n}$, as follows:
\begin{enumerate}
\item 
At time $t=0$, start with the empty graph on the vertex set $\vertex = [n]$.
\item 
At times $t\in\N$, add the edge $\{a,b\}$, where $a,b$ are such that $\Pi_n(t) = \Pi_n(t-1)\circ (a,b)$.
\end{enumerate}
\hfill$\spadesuit$
\end{definition}

\noindent
Associated graph processes were used in \cite{S2005,BD2006} and follow-up articles to represent the evolution of a general dynamic permutation in terms of a dynamics generated by applying a single transposition at every time step.

A crucial role will be played by the first time when the support of the random walk distribution \emph{intersects the largest connected component of the associated graph process}:

\begin{definition}[Largest component of the associated graph process]
\label{def:maxcomp}
Denote by $\maxgraph{\AG{\Pi_n}(t)}{} $ the set of vertices in the largest connected component in the associated graph process at time $t$. If such a connected component is not unique, then take all the vertices in all the largest connected components. \hfill$\spadesuit$
\end{definition}

\begin{remark}[Possible non-uniqueness of the largest connected component]
In situations where we employ Definition~\ref{def:maxcomp}, the largest connected component is unique with high probability. Situations where it is not unique will be of no importance.
\end{remark}

\begin{definition}[Drop-down time]
\label{def:Mtdd} 
Fix any $\varepsilon_n>0$ such that $\varepsilon_n = \omega(n^{-1/3})$ and $\varepsilon_n=o(1)$ as $n\to\infty$. The \emph{drop-down time} is defined as
\begin{equation}\label{eq:Mtdddef}
\Mtddalt = \inf \Big\{t > \tfrac{n}{2}[1+\varepsilon_n]\colon\, \fsupp{\ISRWdistr{}(t)} \cap \mathscr{C}_{\rm max}(\AG{\Pi_n}(t))\neq \emptyset\Big\}.
\end{equation}
\hfill$\spadesuit$
\end{definition}

\begin{remark}[Drop-down time and hitting time of the largest permutation cycle]
At first sight it might seem unintuitive that the time $\Mtddalt$ from Definition~\ref{def:Mtdd} plays an important role. Given the diffusive nature of ISRW, an arguably more natural candidate would be the first time when the ISRW is supported on the largest \emph{permutation cycle}. However, the above definition in terms of the associated graph process allows for a unified presentation of our results in different settings, even when the associated graph process at a single time \emph{does not} provide all the information about the structure of permutation cycles.
\end{remark}

For CDP, the drop-down time is the first time when the cycle that contains $v_0$ merges with the largest cycle. For CFDP, however, this is not necessarily true because cycles fragment. We therefore define the drop-down time to be the first time when the random walk `sees' the maximal component, see \eqref{eq:Mtdddef}. Later, we will see that, in fact, afterwards the mass spreads over $\mathscr{C}_{\rm max}(\AG{\Pi_n}(t))$ quickly.

\begin{remark}[Properties of drop-down time]
Clearly, $\tddalt$ is random. However, if we condition on a particular realisation of $\Pi_n$, then $\tddalt$ is a deterministic function of the starting point of the random walk. The role of $\varepsilon_n$ is to ensure that $\tddalt$ represents the first time in the \emph{supercritical regime} when the \emph{largest component} in the associated \emph{\ER graph process} coincides with the support of the ISRW, see Section~\ref{sec:agp}. The choice of $\varepsilon_n$ ensures that the drop-down time \emph{avoids the critical window}, which corresponds to $\smash{\tfrac{n}{2} + O(n^{2/3})}$, yet covers the entire supercritical regime.
\end{remark}

%%%

\subsection{Main results}
\label{sec:intro:results}

For convenience, we introduce the following shorthand notation:

\begin{definition}[Total variation distance away from equilibrium]
For $v\in[n]$, define
\begin{equation}
\dnTVD{t} = d_{\rm TV} \left( \ISRWdistr{}(t) , \Unif{[n]} \right), \qquad t \in \N_0.
\end{equation}
\hfill$\spadesuit$
\end{definition}

\noindent
Our main results are the following two theorems ($\convdist$ denotes convergence in distribution):

\begin{theorem}[Mixing profile for ISRW on CDP]
\label{thm:main} 
$\mbox{}$
\begin{enumerate}
\item[\rm{(1)}] 
Uniformly in $\startvertex\in[n]$,
\begin{align}
\frac{\tddalt}{n} \convdist s^\Downarrow,
\end{align}
where $s^\Downarrow$ is the $[0,1]$-valued random variable with distribution  (recall \eqref{eq:eta})
\begin{equation}
\P(s^\Downarrow \leq s) =  \eta(s),\qquad s\in[0,1].
\end{equation}
\item[{\rm (2)}] 
Uniformly in $\startvertex\in[n]$, 
\begin{align}
(\dnTVD{sn})_{s \in [0,1]} \convdist \big(1-\eta(s)\1_{\{s>s^\Downarrow\}}\big)_{s \in [0,1]} 
\quad \mbox{in the Skorokhod $M_1$-topology}.
\end{align}
\end{enumerate}
\end{theorem}

\begin{theorem}[Mixing profile for ISRW on CFDP]
\label{thm:frag:main} 
$\mbox{}$ 
\begin{enumerate}
\item[{\rm (1)}] 
Uniformly in $\startvertex\in[n]$,
\begin{align}
\frac{\tddalt}{n} \convdist u^\Downarrow,
\end{align}
where $u^\Downarrow$ is the non-negative random variable with distribution (recall Definition~\ref{def:ERstruct}(1))
\begin{equation}
\P(u^\Downarrow \leq u) =  \zeta(u), \qquad u\in[0,\infty).
\end{equation}
\item[{\rm (2)}] 
Uniformly in $\startvertex\in[n]$, 
\begin{align}
(\dnTVD{un})_{u \in [0,\infty)} \convdist \big(1-\zeta(u)\1_{\{u>u^\Downarrow\}}\big)_{u \in [0,\infty)}
\quad \mbox{in the Skorokhod $M_1$-topology}.
\end{align}
\end{enumerate}
\end{theorem}

\noindent
The proofs of these theorems are given in Sections~\ref{sec:CDP} and \ref{sec:CFDP}, respectively.

%%%

\subsection{Discussion}
\label{sec:intro:discussion}

\medskip\noindent
{\textbf 1.}
Despite the similarity of Theorems~\ref{thm:main}--\ref{thm:frag:main}, the latter is \emph{far more delicate}. For CDP, mixing is simply induced by the ISRW entering the ever-growing largest cycle. For CFDP, the presence of fragmentations breaks the direct link between the dynamic permutation and its associated graph process: \emph{a single connected component may carry more than one permutation cycle}. Specifically, the largest component of the associated graph process carries a large number of permutation cycles and, at the drop-down time, the distribution of the ISRW is supported on only one of them. It is not a priori clear how many steps the dynamics needs to spread out the ISRW distribution over all the elements that lie on the largest component of the associated graph process. Therefore, a major hurdle in the proof of Theorem~\ref{thm:frag:main} is to show that such \emph{local mixing} happens on time scale $o(n)$. We actually show a stronger statement, namely, that local mixing occurs on an arbitrarily small but diverging time scale (see Section~\ref{sec:CFDP:localmix} for details). The core of the proof is to show that on the largest component over time there is a diverging count of appearances of permutation cycles that span almost the entire largest component of the associated graph process.

\begin{figure}[htbp]
\centering
\includegraphics[width=0.6\linewidth]{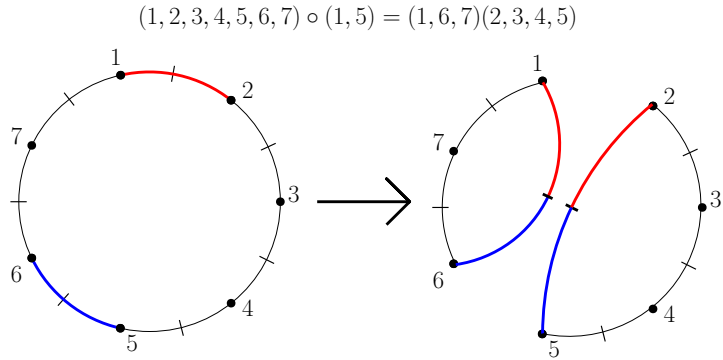}
\caption{\small Dynamic permutations are similar to rewirings in the configuration model, where all degrees are two. Recall Example~\ref{ex:coalfrag}. Consider the permutation $\Pi(0) = (1,2,3,4,5,6,7)$, which consists of a single cycle and corresponds to a degree-two graph that has a single connected component. Apply the transposition $(1,5)$ to get a new permutation $\Pi(1) = \Pi(0) \circ (1,5) = (1,6,7) (2,3,4,5)$, which consists of two cycles and corresponds to a degree-two graph that has two connected components, obtained by sampling the edges $(1,\Pi(0)(1)) = (1,2)$ and $(5,\Pi(0)(5)) = (5,6)$ and rewiring them. }
\label{FigPerRew}
\vspace{-.3cm}
\end{figure}

\medskip\noindent
{\textbf 2.} 
Theorems \ref{thm:main}--\ref{thm:frag:main} extend our earlier results for the total variation distance of a (non-backtracking) random walk on a configuration model subject to random rewirings \cite{AGHHN2020}. There we assumed that \emph{all the degrees are at least three}, which corresponds to a \emph{supercritical} configuration model that with high probability is connected (see \cite[Chapter~4]{vdH2018}). Our model with evolving permutation cycles is closely related to the setting where \emph{all the degrees are two}, which in turn corresponds to a special kind of configuration model that with high probability is disconnected (see Fig.~\ref{FigPerRew}). In this setting, even small perturbations of the degree sequence can lead to significantly different behaviour (see \cite{F2020} for details). In Appendix \ref{apx:deg2} we comment further on the connection between permutations and degree-two graphs. More concretely, we show that in the setting of dynamic degree-two graphs with rewiring, we obtain an ISRW-mixing profile analogous to the one described in Theorem~\ref{thm:frag:main} (see Theorem~\ref{thm:cm2:main}). We stress that in the present work the starting configuration is fixed to be the identity permutation, which would correspond to a graph with only self-loops, whereas in our previous work the starting configuration was sampled from the configuration model.

%%%%%%%%%%%%%%%%%%%%%%%%%%%%%%%%%%%%%%%%%%%%%%%%%%%%%%%%%%%%
%%%%%%%%%%%%%%%%%%%%%%%%%%%%%%%%%%%%%%%%%%%%%%%%%%%%%%%%%%%%

\medskip\noindent
{\textbf 3.}
The mixing profile in Theorems \ref{thm:main}--\ref{thm:frag:main} is unusual: the total variation distance makes a \emph{single jump} down from the value $1$ to a value on a deterministic curve and subsequently follows this curve on its way down to $0$. This jump, which is similar to a \emph{one-sided cut-off}, occurs after a \emph{random time}. The law of the drop-down time and the function describing the deterministic curve depend on the choice of dynamics.

\medskip\noindent
{\textbf 4.}
The pathwise statements in part (2) of Theorems \ref{thm:main}--\ref{thm:frag:main} imply the following pointwise statements ($\sim$ denotes equality in distribution): 
\begin{equation}
\begin{array}{lll}
&\dnTVD{sn} \convdist 1 - \eta(s)Y(s), \quad s \in [0,1], &\text{ with } Y(s)\sim \mathsf{Bernoulli}(\eta(s)),\\
&\dnTVD{sn} \convdist 1 - \zeta(s)\bar Y(s), \quad s \in [0,\infty), &\text{ with } \bar Y(s)\sim \mathsf{Bernoulli}(\zeta(s)).
\end{array}
\end{equation}
Through the function $\phi$ plotted in Fig.~\ref{fig:graphs}, we can view the two mixing profiles as a continuous deformation of one another. Slower mixing for CFDP is intuitive: fragmentation slows down the mixing, while coagulation induces it.

\medskip\noindent
{\textbf 5.}
Note the similarities between the mixing profiles described by Theorems~\ref{thm:main}--\ref{thm:frag:main}. Both feature a single macroscopic jump at a random time to a deterministic curve that depends on the choice of the dynamics. We expect this type of behaviour to occur for any permutation dynamics whose associated graph process exhibits scaling behaviour similar to that of the \ER graph process. A class of graph processes that fits this criterion is the class of \emph{Achlioptas processes with bounded-size rules} (see \cite{RW2017} or \cite{SW2007}).

\medskip\noindent
{\textbf 6.}
We can formulate  conjectures about finite-speed random walks as well. Settings where the random walk rate dominates are easy to handle. If the random walk is fast enough to ensure local mixing (e.g.\  $\gg n^2$ steps of the random walk occur for every step of the random permutation), then our theorems should remain the same with negligible error terms. In this regime, the mixing is fully driven by the underlying geometry. However, once these rates are commensurate, we would have to deal with random walk distributions that are \emph{partially mixed over cycles}, meaning that the distribution of the random walk would not be uniform over its supporting cycle before this cycle is affected by the permutation dynamics.

\medskip\noindent
{\textbf 7.}
Dynamic permutations are a natural model for discrete dynamic random environments, which typically are \emph{disconnected} but nonetheless allow for interaction between their constitutive elements. We believe this setting to be interesting for other stochastic processes on random graphs as well, such as the \emph{voter model} or \emph{the contact process}.

%%%

\paragraph{\bf Organisation of the paper.}
\label{sec:intro:outline}

Section~\ref{sec:CDP} starts by establishing a link between CDP and cycle-free \ER random graphs. A \emph{coupling} construction is employed to describe the \emph{cyclic structure} of a typical CDP. These results are used to prove Theorem~\ref{thm:main}. Section~\ref{sec:CFDP} deals with CFDP, where the main problem is that the associated graph process provides weaker control over permutation cycles than for CDP. After this discrepancy is settled, we employ arguments analogous to those in Section~\ref{sec:CDP} to prove Theorem~\ref{thm:frag:main}. 

Appendices~\ref{apx:limit}--\ref{apx:norm} contain supplementary material that is not needed in Sections~\ref{sec:CDP}--\ref{sec:CFDP}. Namely, Appendix~\ref{apx:limit} shows that the ISRW arises as a fast-speed limit of the standard random walk.  Appendix~\ref{apx:norm} proves that the laws of the jump-down times in Theorems~\ref{thm:main}--\ref{thm:frag:main} are properly normalised. Appendix~\ref{apx:schramm} contains the key coupling that is used to study the cycle structure of CFDP, which is technical and of interest in itself. This coupling is needed in Section~\ref{sec:CFDP}. Appendix~\ref{app:sc} contains a technical computation that is needed in Section~\ref{sec:CFDP} as well. Finally,  Appendix~\ref{apx:deg2} elucidates the connection between random permutations and graphs with all the degrees equal to two and extends Theorem~\ref{thm:frag:main} to the setting of dynamic degree-two graphs.

%%%%%%%%% SECTION 2 %%%%%%%%%%%%%%%%%%%%%%%%%

\section{Coagulative dynamic permutations}
\label{sec:CDP}

In this section, we establish a link between dynamic permutations and evolving graphs. To do so, we couple a CDP with a cycle-free \ER graph process (Section~\ref{sec:agp}), and couple the latter with the standard \ER graph process (Section~\ref{sec:cfcc}) by making use of well-known results on the structure of connected components of \ER random graphs (recall Section~\ref{sec:intro:ER}). We use the couplings to prove Theorem~\ref{thm:main} (Section~\ref{sec:CDPmixing}).

%%%

\subsection{Representation via associated graph process}
\label{sec:agp}

Note that for the dynamics generated by transpositions sampled uniformly at random from the set of all transpositions of $n$ elements, the associated graph process is equal in distribution to the \ER process defined in Definition~\ref{def:ER1}. In the setting of \emph{coagulative} dynamic transpositions, this leads us to the following observation:

\begin{lemma}[Representation of CDP as cycle-free graph process]
\label{prop:ER}
If $\Pi_n$ is a CDP, then its associated graph process $\AG{\Pi_n}$ is the cycle-free \ER graph process defined in Definition~\ref{def:cfER}.
\end{lemma}

\begin{proof}
Recall that the change between two successive permutations in a CDP is generated by applying a \emph{single} transposition. Furthermore, note that the only transpositions causing a split of a permutation cycle are the ones that transpose two elements from the same cycle. Recall Definition~\ref{def:AG}, and note that if $\AG{\Pi_n}(t)$ is a forest, then its connected components correspond to cycles of $\Pi_n(t)$. Furthermore, observe that  cycle-splitting transpositions correspond to edges that join two vertices from the same connected component. Thus, if $\AG{\Pi_n}(t)$ is a forest, then any transposition causing a fragmentation of a permutation cycle corresponds to an edge that creates a cycle in the associated graph process.

Observe that the associated graph process always starts as a forest. Since fragmentations of permutation cycles are not allowed, there can be no edges that lead to graph cycles in the associated graph process. Since the associated graph process for a dynamic permutation with no constraints is the \ER graph process, the associated graph process for a CDP is the \ER graph process constrained to be a forest (see Definition~\ref{def:cfER}).
\end{proof}

%%%

\subsection{Connected components of the cycle-free \ER graph process}
\label{sec:cfcc}

%%%

\paragraph{$\bullet$ Coupling of \ER graph processes.}
We construct a coupling of the standard and the cycle-free \ER graph process that allows us to study the structure of the connected components of the cycle-free process.

\begin{definition}[Coupling between cycle-free and standard \ER graph process]
\label{def:graph-coupling}
Let $G_n = (G_n(t))_{t\in\N_0}$ be the \ER graph process on $[n]$ defined in Definition~\ref{def:ER1}, and denote the edge set of $G_n(t)$ by $\edgeset_{G_n(t)}$. Based on $G_n$, construct a graph-valued process $F_n = (F_n(t))_{t\in\N_0}$ as follows:
\begin{enumerate}
\item 
$F_n(0)$ is the empty graph with vertex set $[n]$.
\item 
At times $t\in\N$, define $e^\star(t) = \edgeset_{G_n(t)} \setminus \edgeset_{G_n(t-1)}$, which is the edge added at time $t$ to $G_n(t)$.
\begin{enumerate}
\item 
Construct the candidate graph at time $t$, defined as $F_{n}^{\star}(t) = (\vertex, \edgeset_{F_n(t-1)} \cup \{e^\star(t)\})$.
\item 
If $F_{n}^{\star}(t)$ is a forest, then set $F_n(t) = F_{n}^{\star}(t)$.
\item 
Otherwise, set $F_{n}(t) = F_{n}(t-1)$.
\end{enumerate}
\end{enumerate}
Define the \emph{effective time} $\tau_n(t)$ of the coupled process $(F_n(t))_{t\in\N_0}$ by setting $\tau_n(0) = 0$ and, recursively for $t \in \N$,
\begin{equation}
\tau_n(t) = 
\begin{cases}
\tau_n(t-1) + 1,  &\qquad\text{if $F_n(t) \neq F_n({t-1})$, i.e., the proposed edge has been accepted,}\\
\tau_n(t-1),        &\qquad\text{if $F_n(t) = F_n({t-1})$, i.e., the proposed edge has been rejected.}
\end{cases}
\end{equation}
Note that $\tau_n(t)$ is a random variable because it is a function of a random graph process. We suppress the dependence of $\tau_n(t)$  on $F_n$, since we will never work with more than one set of coupled processes at a time. \hfill$\spadesuit$
\end{definition}

\begin{remark}[Relation between $F_n$ and $\AG{\Pi_n}$]
\label{rem:CC}
By the definition of the coupling, if there are edge-rejections at times $\{t, t+1, t+2, \ldots, t+k\}$, then a string of $k+1$ copies of the same graph is observed in $F_n$, i.e., $ F_n(t-1) = F_n(t) = F_n(t+1) = \cdots = F_n(t+k)$. On the other hand, the associated graph process $\AG{\Pi_n}$ is by construction a sequence of graphs such that no two graphs are the same. To recover $\AG{\Pi_n}$ from $F_n$, from every string of copies of the same graph choose only one copy of that graph.
\end{remark}

The reason why this construction is useful to control the connected components of the cycle-free \ER graph process is stated in the following lemma:

\begin{lemma}[Connected components of $F_n$]
\label{lem:cc}
Let $H$ be a graph with vertex set $\vertex$, and define $\conncomp{H}$ to be the partition of $\vertex$ induced by the connected components of $H$. Let $G_n, F_n$ be as in Definition~\ref{def:graph-coupling}. Then, at every time $t\in\N_0$, $\conncomp{G_n(t)} = \conncomp{F_n(t)}$.
\end{lemma}

\begin{proof}
Note that any edge creating a cycle does not influence the size of the connected components.
\end{proof}

%%%

\paragraph{$\bullet$ Effective time.}
To use the above observation, we need to control the effective time $\tau_n(t)$. The following lemma shows that with high probability and after scaling by $1/n$, there is a simple relation between the standard time $t$ and the effective time $\tau_n(t)$:

\begin{lemma}[Effective time of a cycle-free \ER graph process]
\label{lem:et}
Let $G_n$, $F_n$ and $\tau_n$ be as in Definition~\ref{def:graph-coupling}, and $\phi$ as in Definition~\ref{def:phi}. Then, for any for $u \in [0,\infty)$,
\begin{equation}\label{eq:et}
\frac{\tau_n(un)}{n} \convprob \phi(u).
\end{equation}
\end{lemma}

\begin{proof}
The proof of Lemma~\ref{lem:et} consists of two separate lines of argument. First, we show that the left-hand side in \eqref{eq:et} concentrates around a deterministic quantity. Afterwards, the value of this quantity is computed. 

\medskip\noindent
\textsc{Part 1: Concentration of the associated martingale.}
Observe that 
\begin{equation}
\label{eq:tau2}
\tau_n(t) = t - \sum_{s=0}^t \1_{R(s)}, \qquad t \in \N_0,
\end{equation}
where $\1_{R(s)}$ is the edge-rejection indicator at time $s$. Let $\mathcal{F} = (\mathcal{F}_n(t))_{t\in\N_0}$ with $\mathcal{F}_n(t) = \sigma((G_n(q))_{q=0}^{t})$ be the natural filtration with respect to the \ER graph process. By the construction of the coupling, a rejection occurs whenever there is an edge that creates a cycle within a connected component. Therefore
\begin{equation}
\label{eq:rejcondprob}
\1_{R(0)} = 0, \qquad \left[\1_{R(s)}\cond \mathcal{F}_n(s-1)\right] \sim \mathsf{Bernoulli}(p_s), \qquad s \in \N,
\end{equation} 
where the success probabilities are given by
\begin{equation}
\label{eq:rejprob}
p_s = \sum \limits_{\substack{ \mathscr{C} \in \conncomp{G_n(s-1)}}} \frac{|\mathscr{C}|(|\mathscr{C}|-1)}{n(n-1)}, \qquad s \in \N,
\end{equation}
i.e., the number of edges that can join two vertices from the same connected component at time $s-1$ divided by the total number of edges. Introduce the shorthand notation
\begin{align}
\E_{t}[\cdot] = \E\left[\,\cdot \cond \mathcal{F}_n(t) \right],
\end{align}
and define two sequences of random variables $(D_t)_{t\in\N_0}$ and $(S_t)_{t\in\N_0}$ such that
\begin{equation}
\label{eq:martdef}
\begin{aligned}
S_0 &= D_0 = 0,\\
D_t &= \1_{R(t)} - \E_{t-1}\left[ \1_{R(t)}\right], \qquad t\in\N,\\ 
S_t &= \sum_{s=0}^t D_s = \sum_{s=0}^t \1_{R(s)} - \sum_{s=0}^t \E_{s-1}\left[ \1_{R(s)} \right], \qquad t\in\N.
\end{aligned}
\end{equation}
Note that, for any $t\in\N_0$,
\begin{align}
\E[S_t] &\leq t < \infty, \nonumber \\
\E_t\left[ S_{t+1}\right] &= \E_t\left[ \1_{R(t+1)} - \E_{t}\left[ \1_{R(t+1)} \right]\right] + \E_t\left[ S_t \right] = S_t, \\[0.1cm] \nonumber
|S_t - S_{t-1}| &= |D_t| \leq 1.
\end{align}
Hence, $(S_t)_{t\in\N_0}$ is a martingale with bounded differences with respect to the natural filtration of the \ER graph process. Using the Azuma-Hoeffding inequality, we can estimate 
\begin{equation}
\label{eq:azuma}
\P \left( |S_t| \geq \varepsilon \right) \leq 2 \exp\left(-\frac{\varepsilon^2}{2t}\right).
\end{equation}
Pick $t = un$, $u\in [0,\infty)$, and $\varepsilon_n = n^{\frac{1+\delta}{2}}$, $c>0$, $\delta \in (0,1)$. Introduce the event
\begin{equation}
\Xievent(un) = \{ |S_{un}| < n^{\frac{1+\delta}{2}} \}.
\end{equation}
By \eqref{eq:azuma},
\begin{equation}
\label{eq:bound}
\P\left( \Xievent^\cmpl(un) \right) \leq 2 \exp\left(-\frac{c^2 n^\delta}{2u}\right) = o(1)
\end{equation}
and hence $\P\left( \Xievent(un) \right) = 1-o(1)$. By the definition of $S_{t}$ in \eqref{eq:martdef}, we see that, on the event $\Xievent(un)$,
\begin{equation}
\1_{\Xievent(un)} \left[\frac{1}{n} \sum_{s=0}^{{un}} \1_{R(s)} -\frac{1}{n} \sum_{s=0}^{{un}} \E_{s-1}\left[\1_{R(s)} \right]\right]
= \1_{\Xievent(un)}\, o(1), \qquad \text{with  } |o(1)| \leq n^{\frac{1+\delta}{2}},
\end{equation}
which establishes the concentration of ${\tau_n(un)}/{n}$ as $n\to\infty$.

\medskip\noindent
\textsc{Part 2: Computation of limit.}
We compute
\begin{align}
\label{eq:Rlimit}
\frac{1}{n} \sum \limits_{s=0}^{{un}} \1_{R(s)} 
&= \frac{(\1_{\Xievent(un)} + \1_{\Xievent^\cmpl(un)})}{n} \sum_{s=0}^{{un}} \1_{R(s)}\\
&= \frac{1}{n} \left[\1_{\Xievent(un)} \left( \sum_{s=0}^{{un}} \E_{s-1}\left[ \1_{R(i)} \right] + o(n)\right) 
+  \1_{\Xievent^\cmpl(un)} \sum_{s=0}^{{un}}\1_{R(s)} \right].
\label{eq:betasplit}
\end{align}

\medskip\noindent
{\bf 1.} 
Observe that
\begin{equation}
\frac{1}{n} \1_{\Xievent^\cmpl(un)} \sum \limits_{s=0}^{{un}}\1_{R(s)} \convprob 0,
\end{equation}
because $\smash{\sum_{s=0}^{{un}}\1_{R(s)} \leq un}$ and $\smash{\1_{\Xievent^\cmpl(un)} \convprob 0}$. To understand the first summand in \eqref{eq:betasplit}, we need to introduce another event. Fix a sequence $\varepsilon_n$ such that $\varepsilon_n = o(1)$ and $\varepsilon_n = \omega(n^{-1/3})$. For any $n\in\N$, $t\in\N_0$, define the \ER \emph{typicality} event 
\begin{align}\label{eq:ERtypical}
\typicalevent{}{un} &= \big\{ |\maxgraph{n,vn}{ER}| \in (n[\zeta(v) - \varepsilon_n], n[\zeta(v) + \varepsilon_n])~\forall\, v\in [\tfrac{1}{2},u] \big\} 
\cap \big\{ 0 \leq \max_{s\leq un} |\secgraph{n,s}{ER}| \leq n\varepsilon_n \big\}.
\end{align}
Then we can write
\begin{align}
\1_{\Xievent(un)} \sum_{s=0}^{{un}} \E_{s-1}\left[ \1_{R(s)} \right] 
= \1_{\Xievent(un)} \left( \1_{\typicalevent{}{un}}  
+ \1_{\typicalevent{c}{un}}\right)\sum_{s=0}^{{un}} \E_{s-1}\left[ \1_{R(s)} \right].
\end{align}
Again, we see that 
\begin{equation}
\frac{1}{n}\1_{\typicalevent{c}{un}}\sum_{s=0}^{{un}} \E_{s-1}\left[ \1_{R(s)} \right] \convprob 0,
\end{equation}
because $\smash{\sum_{s=0}^{{un}} \E_{s-1}[\1_{R(s)}] \leq un}$ and $\smash{\1_{\typicalevent{c}{un}} \convprob 0}$.

\medskip\noindent
{\bf 2.}
It remains to compute $\smash{\frac{1}{n}}\1_{\Xievent(un)} \1_{\typicalevent{}{un}} \sum_{s=0}^{{un}} \E_{s-1}[\1_{R(s)}]$, which is the only term that will be non-zero after we take the limit $n\to\infty$ in \eqref{eq:Rlimit}. Recall that $\E_{s-1}[\1_{R(s)}] \sim \mathsf{Bernoulli}(p_s)$, where the success probabilities $p_s$ were introduced in \eqref{eq:rejprob}.

% We will use the following well-known results on the sizes of \ER connected components (see, for example, \cite{FK2016} or \cite{vdH2016}). Choose $\varepsilon_n>0$ such that $\varepsilon_n =\omega(n^{-1/3})$ and $\varepsilon_n = o(1)$. . Then the following hold:
%\begin{equation}
%\begin{aligned}
%\label{eq:ERtyp}
%|\maxgraph{n,t}{ER}| &=
%\zeta( t/n)n + o_{\P}(n), \qquad t \in [0, \infty),\\
%|\secgraph{n,t}{ER}| &\leq 
%\begin{cases}
%O_{\P}(n^{2/3}), &t \leq t_c,\\
%O_{\P}\left(\frac{1}{\varepsilon_n^2}\log (n\varepsilon_n^2) \right), &t\geq t_c.
%\end{cases}
%\end{aligned}
%\end{equation}
%Using \eqref{eq:ERtyp} and the fact that $|\vertex| = n$, we see that
Since $|\vertex| = n$,  we have the following bounds:
\begin{equation}
\label{eq:rejectionbounds}
p_s = \sum \limits_{\substack{\mathscr{C} \in \conncomp{G_n(s-1)}}} \frac{|\mathscr{C}|(|\mathscr{C}|-1)}{n(n-1)} 
\leq
\begin{cases}
\frac{|\maxgraph{n, n/2}{ER}|}{n}, &0\leq s  \leq \tfrac12 n,\\
%\frac{[n\zeta({s/n})]^2 + o_{\P}(n^2)}{n(n-1)} + o_{\mathbb{P}}(1),  &s > \tfrac12 n.
 \frac{|\maxgraph{n, s}{ER}|^2}{n^2} + \max\limits_{0\leq s\leq un}\frac{|\secgraph{n,s}{ER}|}{n} &  \tfrac12 n < s <un ,
\end{cases}
\end{equation}
where the last bound is uniform in $s\leq un$. The first line, for times below $n/2$, holds since
\eqn{
\sum \limits_{\substack{\mathscr{C} \in \conncomp{G_n(s-1)}}} \frac{|\mathscr{C}|(|\mathscr{C}|-1)}{n(n-1)}\leq \frac{|\maxgraph{n, s}{ER}|}{n}
\leq \frac{|\maxgraph{n, n/2}{ER}|}{n}.
}
%via this argument: component sizes are growing in time, and therefore their sizes can be bounded uniformly over any time interval by the size of the largest one at the end of this interval. Furthermore, the bound is essentially  the square root of expression implied by \eqref{eq:rejcondprob}. 
The second line separates the contribution of the maximal component and all the other components, and bound the non-maximal component similarly as in the first line.

In the \emph{supercritical} regime, we separately describe the contribution of the unique largest component and give an upper bound only on the probability of rejection due to the other components. On the \ER typicality event $\typicalevent{}{un}$ (recall \eqref{eq:ERtypical}),  the size of all the components in the subcritical and critical regime and all the components but the largest one in the supercritical regime can be uniformly bounded by $|\mathscr{C}| \leq Zn^{2/3}$, where $Z$ is a positive random variable. From \eqref{eq:rejcondprob}, \eqref{eq:ERtypical} and  \eqref{eq:rejectionbounds}, it follows that, on the event $\typicalevent{}{un}$,
\begin{align}
\1_{\typicalevent{}{un}} \sum_{s=0}^{{un}} \E_{s-1}\left[ \1_{R(s)} \right] = \1_{\typicalevent{}{un}} \sum_{s=0}^{{un}} p_s  
\begin{cases}
 \leq u |\maxgraph{n,n/2}{ER}|\leq \vep_n u, &0 \leq u \leq \tfrac{1}{2},\\
= \sum_{s= \frac{n}{2}}^{un} (\zeta(i/n)+\vep_n)^2 + \mathcal{R}(un),  &u > \tfrac{1}{2},
\end{cases}
\end{align}
where we use that $\zeta( \tfrac{1}{2})=0$. Since $|\maxgraph{n,n/2}{ER}|\leq \vep_n n$ on $\typicalevent{}{un}$, 
 the remainder term $\mathcal{R}(un)$ can be bounded as (recall \eqref{eq:rejectionbounds})
\begin{align}
\mathcal{R}(un) \leq un \max\limits_{0\leq s\leq un}\frac{|\secgraph{n,s}{ER}|}{n} \leq \vep_n u n = o(n).
\end{align}
%Specifically, the error term $o(n)$ in the last expression follows from the bounds in \eqref{eq:ERtypical} and \eqref{eq:rejectionbounds}.
%\begin{equation}
%\1_{\typicalevent{}{un}}  \sum_{s = \frac{n}{2}+ O(n^{2/3})}^{ un } \frac{o(n^2)}{n(n-1)} 
%= \sum_{s = \frac{n}{2}+ O(n^{2/3})}^{ un } o(1) = o(n).
%\end{equation}
 
\medskip\noindent
{\bf 3.}
Before wrapping up, let us note that
\begin{align}
 \frac{ o(n) + \sum_{s=0}^{{un}} \zeta(s/n)^2}{n} 
\to \int_0^u \mathrm{d}v\,\zeta^2(v),
\end{align}
because $\zeta^2$ is bounded and hence Riemann integrable over compact intervals $[0,u]$, and $\smash{\frac{1}{n}\sum_{s=0}^{{un}} \zeta(s/n)^2}$ is a Riemann sum of $\zeta^2$ over a regular partition of $[0,u]$ into subintervals of length $1/n$. This allows us to finish our previous computation, namely,
\begin{align}
\frac{1}{n} & \sum_{s=0}^{{un}} \1_{R(s)} \nonumber \\
&= \frac{1}{n} \left[\1_{\Xievent(un)} \left( \1_{\typicalevent{}{un}} 
+ \1_{\typicalevent{c}{un}}\right) \left( \sum_{s=0}^{{un}} \E_{s-1}\left[ \1_{R(s)} \right] + o(n)\right) 
+  \1_{\Xievent^\cmpl(un)} \sum_{s=0}^{{un}}\1_{R(s)} \right]\\ \nonumber
&=\frac{1}{n}\left(o(n) + \1_{\Xievent(un)}\1_{\typicalevent{}{un}} \sum_{s=0}^{{un}} \E_{s-1}\left[ \1_{R(s)} \right] \right) 
\convdist \int_0^u \mathrm{d}v\,\zeta^2(v) = u - \phi(u) 
\end{align}
(recall \eqref{def:phi}), from which the desired result follows.
\end{proof}

%%%

\paragraph{$\bullet$ Mapping between times.}
The main purpose of Lemma \ref{lem:et} is to show that on time scales of order $n$ there is a function $\phi$ (recall Definition~\ref{def:ERstruct}) capturing the correspondence, in the limit as $n\to\infty$, between the times at which the standard \ER graph process and its cycle-free counterpart have certain quantities distributed equally, notably, the sizes of their connected components. Since $\phi$ is strictly monotone, it admits a proper inverse, which allows us to relate the cycle-free graph process to the standard \ER graph growth process.

%%%

\paragraph{$\bullet$ Largest component.}
To conclude the analysis of the cycle-free graph process, we combine the above results to obtain a characterisation of the largest component of the cycle-free \ER graph process:

\begin{lemma}[Size of the largest component]
\label{lem:CFmax}
For $s\in[0,1]$, let $|\maxgraph{n, sn}{cfER}|$ be the size of the largest component of the cycle-free \ER graph process on $n$ vertices at time ${sn}$. Then
\begin{equation}
\frac{|\maxgraph{n,{sn}}{\rm cfER}|}{n} \convprob \eta(s), \qquad \eta(s) = \zeta(\phi^{-1}(s)).
\end{equation}
\end{lemma}

\begin{proof}
In Lemma \ref{lem:cc} we have established that for every $n\in\N$ the connected components of the cycle-free graph $G_n^{\rm{cf}}(t)$ correspond exactly to the connected components of the standard \ER graph process at some random time, denoted by $\tau_n^{-1}(t)$. In Lemma \ref{lem:et} we have established that $\tau_n^{-1}(t) = n\phi^{-1}(t/n) + o_{\mathbb{P}}(n)$ and $|\mathscr{C}_{\rm max}^{\rm ER}(n, sn)| = n\zeta(s) + o_{\mathbb{P}}(n)$. Hence 
\begin{equation}
|\maxgraph{n,sn}{\rm cfER}| \deq |\maxgraph{n,n(\phi^{-1}(s) + o_{\mathbb{P}}(1))}{\rm ER}| = n\zeta(\phi^{-1}(s) + o_{\mathbb{P}}(1)) + o_{\mathbb{P}}(n),
\end{equation}
from which the claim follows.
\end{proof}

%%%

\subsection{Drop-down time and mixing profile}
\label{sec:CDPmixing}

As stated in Theorem~\ref{thm:main}, for CDP the mixing profile exhibits a cut-off at a \emph{random} time. From that moment onwards, the total variation distance follows a \emph{deterministic} curve that is related to the typical structure of CDP. The following lemma gives the distribution of the \emph{drop-down} time and settles Theorem~\ref{thm:main}(1):

\begin{lemma}[Limit distribution of drop-down time for ISRW on CDP]
\label{lem:cdp:tdddistr}
Recall $\tddalt$ from Definition~\ref{def:Mtdd}. There exists a $[0,1]$-valued random variable $s^{\Downarrow}$ with a distribution function $\P(s^{\Downarrow} \leq s) = \eta(s)$, $s\in [0,1]$, such that
\begin{equation}
\frac{\tddalt}{n} \convdist s^{\Downarrow}.
\end{equation}
\end{lemma}

\begin{proof}
By the arguments in the proof of Lemma~\ref{prop:ER}, the sizes of the connected components of the cycle-free associated graph process exactly correspond to the sizes of the permutation cycles of the CDP at a given time $t$. Therefore we must study the probability that a uniform vertex lies on the largest component of a cycle-free \ER graph process.

Let $\P_n$ denote the law of CDP on $[n]$ and $\P_n^{\rm cfER}$ the law of the associated graph process, which is a cycle-free \ER graph process (recall Definition~\ref{def:AG}  and Lemma~\ref{prop:ER}). Fix a sequence $\varepsilon_n $ such that $\varepsilon_n = o(1)$ and $\varepsilon_n =\omega(n^{-1/3})$, and for $n\in\N$, $t\geq0$ define the \emph{typicality} event 
\begin{align}
\label{cfER-event}
\typicalevent{cfER}{t} &= \left\{ |\maxgraph{n,t}{cfER}| \in n(\eta(t/n) - \varepsilon_n, \eta(t/n) + \varepsilon_n, ) \right\} 
\cap \left\{ |\secgraph{n,t}{cfER}| \leq n\varepsilon_n \right\}.
\end{align}
For $s\in [0,\tfrac{1}{2})$, we have $\P_n(\tddalt \leq sn) =0$ by the definition of $\tddalt$. For $s\in[\tfrac{1}{2},1]$, by Lemma~\ref{prop:ER},
\begin{equation}
\begin{aligned}
&\P_n(\tddalt \leq sn) = \P_n^{\rm cfER} \big( \startvertex \in \maxgraph{sn}{cfER} \big) \\
&\qquad= \P_n^{\rm cfER} \big( \{\startvertex \in \maxgraph{sn}{cfER}\} \cap \typicalevent{cfER}{sn} \big)
+ \P_n^{\rm cfER} \big( \{\startvertex \in \maxgraph{sn}{cfER}\} \cap [\typicalevent{cfER}{sn}]^\cmpl \big). 
\end{aligned}
\end{equation}
Since the event $\typicalevent{cfER}{sn}$ occurs with probability $1-o(1)$ and $\P_n^{\rm cfER}(\{ \startvertex \in \maxgraph{sn}{cfER} \} \cap \typicalevent{cfER}{sn})$ $= \eta(s) + o(1)$, we see that, for any $s\in[0,1]$,
\begin{align}
\label{eq:Pndef}
\P_n(\tddalt \leq sn) \stackrel{n\to\infty}{\to} \eta(s).
\end{align}
Since $\eta$ is continuous, non-negative and non-decreasing on $[0,1]$ such that $\int_0^1 \mathrm{d}u\,\eta(u) = 1$, \eqref{eq:Pndef} defines a proper distribution function (recall Definition~\ref{def:ERstruct}). See Appendix~\ref{apx:norm} for a detailed computation.
\end{proof}

With the above results in hand, we are ready to prove the pointwise version of Theorem~\ref{thm:main}(2), characterising the mixing profile of ISRW on CDP:

\begin{lemma}[Pointwise limit of mixing profile for ISRWs on CDP]
\label{lem:CDP:profile}
Uniformly in $\startvertex \in [n]$,
\begin{equation}
\dnTVD{sn} \stackrel{d}{\to} 1-\eta(s)Y(s), \qquad s\in[0,1),
\end{equation}   
where $Y(s)\sim\mathsf{Bernoulli}(\eta(s))$.
\end{lemma}

\begin{proof}
Given a permutation $\pi$, let
\begin{equation}
|\MC{\pi}| = \max  \{ |\cycle{\pi}{v}|\colon\, v \in [n] \}
\end{equation}
denote the size of the largest cycle of $\pi$. For every $n\in\N$,
\begin{align}
\dnTVD{sn} = \TVD\Big(\Unif{[n]}, \Unif{[\MC{sn}]}\Big) = 1 - \frac{|\cycle{\startvertex}{sn}|}{n}
\end{align}
by Definition~\ref{def:ISRW} and the definition of total variation distance. Using Lemma~\ref{prop:ER}, we see that
\begin{align}
&\dnTVD{sn} \stackrel{d}{=} 1 - \frac{1}{n} \left(|\graphcycle{\startvertex}{sn}|\1_{ \typicalevent{cfER}{sn}} 
+  |\graphcycle{\startvertex}{sn}|\1_{[\typicalevent{cfER}{sn}]^\cmpl} \right) \nonumber\\
\label{eq:Dsplit} 
&= 1 - \frac{1}{n} \Bigg(|\graphcycle{\startvertex}{sn}|\1_{ \typicalevent{cfER}{sn}}  \left( \1_{\{\tddalt \leq sn\}} + \1_{\{\tddalt > sn\}}\right) 
+  |\graphcycle{\startvertex}{sn}|\1_{[\typicalevent{cfER}{sn}]^\cmpl} \Bigg). 
\end{align}
Standard results for the size of \ER connected components (recall \eqref{eq:ERsizes}) imply that 
\begin{align}
\1_{[\typicalevent{cfER}{sn}]^\cmpl} \convp 0, \qquad
\1_{\typicalevent{cfER}{sn}} \convp 1, \qquad  \frac{|\graphcycle{\startvertex}{sn}| 
\1_{\typicalevent{cfER}{sn}}\1_{\{\tddalt > sn\}}}{n} \convp 0, 
\label{eq:Pconv}
\end{align}
where the last limit follows from the fact that at times $0 \leq sn < \tddalt$ the initial vertex $\startvertex$ lies on a non-largest component, and hence the numerator scales as $o(n)$. Similarly, the size of the largest \ER component has a well-known limit (recall \eqref{eq:ERsizes}), on the event $\{\tddalt<sn\}$, namely,
\begin{align}
\frac{|\graphcycle{\startvertex}{sn}| \1_{\typicalevent{cfER}{sn}} }{n}  \convp \eta(s).
\label{eq:ERgiantconv}
\end{align}
Finally, the only random variable that converges to a non-degenerate random variable is
\begin{align}
\1_{\{\tddalt \leq sn\}} \stackrel{d}{\to} Y(s), \qquad Y(s)\sim\mathsf{Bernoulli}(\eta(s)),
\label{eq:Berconv}
\end{align}
which follows  from Lemma~\ref{lem:cdp:tdddistr}. Note that $\dnTVD{sn}$ is a sum of several random variables, and we established the convergence of each in \eqref{eq:Dsplit}--\eqref{eq:Berconv}. Hence, the claim follows via Slutsky's theorem.
\end{proof}

To conclude this section, we use Lemma~\ref{lem:CDP:profile} to prove the pathwise convergence part of Theorem~\ref{thm:main}:

\begin{proof}[Proof of Theorem~\ref{thm:main}(2)]
Observe that, for every $n\in\N$, any realisation of $\dnTVD{\cdot}$ is a monotone càdlàg path on the compact set $[0,1]$. In this special situation, the pointwise convergence proven in Lemma~\ref{lem:CDP:profile} implies pathwise convergence in the Skorokhod $M_1$-topology. For details, see \cite[Corollary 12.5.1]{W2002}.
\end{proof}

%%%%%%%% SECTION 3 %%%%%%%%%%%%%%%%%%%%%%%%%%%%
 
\section{Coagulative-fragmentative dynamic permutations}
\label{sec:CFDP}

In Section~\ref{sec:CDP}, for CDP it took effort to control the structure of the associated graph process, while the mixing profile was obtained via an easy argument. For CFDP the opposite is true: the associated graph process, introduced in Section~\ref{sec:agp}, is the \ER graph process defined in Definition~\ref{def:ER1} (which is one of the key facts used in \cite{S2005}), while the link between the cycles of the underlying permutation and the connected components of the associated graph process is far less clear. Indeed, each non-tree connected component of the associated graph process may represent \emph{multiple} permutation cycles, which brings a \emph{substructure} into the problem that needs to be controlled. Moreover, it is not a priori clear whether or not this substructure influences the mixing profile, since \emph{immediately} after the drop-down time the distribution of the ISRW is uniform over a component that spans only a random fraction of the largest component of the associated graph process. 

The key result in this section is that ISRW on CFDP exhibits \emph{fast local mixing on the largest component of the associated graph process upon drop-down}. After scaling, this leads to results that are qualitatively similar to those obtained for CDP, namely, the occurrence of a \emph{single} jump in the total variation distance, from $1$ to a deterministic value on a curve related to the largest component of the associated graph process, at a \emph{random time} whose distribution is again connected to the largest component of the associated graph process. The scaled time now takes values in $[0,\infty)$ instead of $[0,1]$.
% Because of this similarity, in some cases we do not write out proofs for CFDP in full, but only point out the differences with respect to their counterparts for CDP.

In Section~\ref{sec:CFDP:ddt} we identify the drop-down time and prove Theorem~\ref{thm:frag:main}(1). In Section~\ref{sec:CFDP:1cycle} we show that the support of ISRW lies on a single permutation cycle before the drop-down time. In Section~\ref{sec:CFDP:localmix} we prove fast local mixing after the drop-down time. In Section~\ref{sec:CFDP:mixprofile} we identify the mixing profile and prove Theorem~\ref{thm:frag:main}(2).

\begin{remark}[Permutation elements and graph vertices representing them]
Throughout this section we will (with a slight abuse of notation) identify the vertices in the associated graph process with the permutation elements they represent.
\end{remark}

%%%

\subsection{Drop-down time}
\label{sec:CFDP:ddt}

Recall that the central object for the identification of the limit distribution of the drop-down time for CDP in Section~\ref{sec:CDPmixing} was the function $\eta$ (recall Definition~\ref{def:ERstruct}), which describes the size of the largest component in the cycle-free \ER graph process. In the setting of CFDP, we formulate a result for $\Mtddalt$ analogous to Lemma~\ref{lem:cdp:tdddistr}, with the role of $\eta$ taken over by $\zeta$, which describes the size of the largest component in the standard \ER graph process:

\begin{lemma}[Limiting distribution of drop-down time for ISRW on CFDP]
\label{lem:cfdp:tdddistr}
Recall $\Mtddalt$ from Definition~\ref{def:Mtdd}. There exists a $[0,\infty)$-valued random variable $u^{\Downarrow}$ with distribution function $\P(u^{\Downarrow} \leq u) = \zeta(u)$, $u\in [0,\infty)$, such that
\begin{equation}
\frac{\Mtddalt}{n} \convdist u^{\Downarrow}.
\end{equation}
\end{lemma}

\begin{proof}
The proof is the same as that of Lemma~\ref{lem:cdp:tdddistr}, but uses the laws of CFDP and its associated graph processes, and uses $\zeta$ in place of $\eta$.
\end{proof}

%%%

\subsection{Drop-down in a single permutation cycle}
\label{sec:CFDP:1cycle}

In principle, it could happen that the ISRW support has experienced fragmentation before the drop-down time, which would significantly complicate our analysis. The main point of this section is to show that, with high probability, this \emph{does not} occur.

\begin{lemma}[ISRW support lies on a single permutation cycle before $\Mtddalt$]
\label{lem:1cycle}
Fix $\varepsilon_n>0$ such that $\varepsilon_n = \omega(n^{-1/3})$ and $\varepsilon_n=o(1)$ as $n\to\infty$. Let $\Omega^{\sss\rm (SC)}(t)$ denote the event that the support of the ISRW at time $t$ lies on a single permutation cycle. Then, uniformly in $v_0$ and $t = cn$ with $c \in (\tfrac12,\infty)$,
\begin{align}
\P\left( \Omega^{\sss\rm (SC)}(t) \mid \Mtddalt \geq t \right) = 1 - o(1).
\end{align}
\end{lemma}

\begin{proof}
Recall the associated graph process introduced in Definition~\ref{def:AG}, and the fact that the associated graph process of CFDP is equal in distribution to the standard \ER graph process. As explained in the proof of Lemma~\ref{prop:ER}, tree components in the associated graph process correspond to permutation cycles that have never experienced fragmentation. The idea of the proof is to show that, conditionally on the event $\{\Mtddalt \geq t\}$, the event $\Omega^{\rm tree}(t)$ that ISRW at time $t$ is supported on a single tree-component in the associated graph process occurs with high probability. Observe that $\Omega^{\rm tree}(t) \subseteq \Omega^ {\rm \sss(SC)}(t)$. First we condition on the event $\{\Mtddalt > t\}$. Afterwards we extend to the event $\{\Mtddalt \geq t\}$. 

Recall that in Lemma~\ref{lem:cfdp:tdddistr} we identified the limiting distribution of $\Mtddalt/n$. Since
\begin{align}
\label{eq:condprob}
\P\left( \Omega^ {\rm tree}(t) \mid \Mtddalt > t \right) = 1 - \P\left( \left[\Omega^{\rm tree}(t)\right]^\cmpl \mid \Mtddalt > t \right) 
= 1- \frac{\P\left(\left[\Omega^{\rm tree}(t)\right]^\cmpl \cap \{\Mtddalt > t\} \right)}{\P(\Mtddalt > t)},
\end{align} 
and the denominator is bounded away from 0 (recall Lemma~\ref{lem:cfdp:tdddistr}), it suffices to show that $\P([\Omega^{\rm tree}(t)]^\cmpl \cap \{\Mtddalt > t\}) = o(1)$. By the law of total probability, we can take the sum over all possible realisations of the underlying dynamics to obtain
\begin{align}
\P \left( \left[\Omega^{\rm tree}(t)\right]^\cmpl \cap \{\Mtddalt > t\} \right) 
&= \E \left[\P\left(\left[\Omega^{\rm tree}(t)\right]^\cmpl\cap \{\Mtddalt > t\} \mid \left(\Pi_n(t)\right)_{s=0}^t\right)\right].
\end{align}
By \cite[Theorem 5.10]{JLR2011}, with high probability the connected components of the associated graph process at time $t$ consist of the \emph{unique largest component, unicyclic connected components and trees}. By Definition~\ref{def:Mtdd}, conditionally on $\{ \Mtddalt > t\}$, the support of the ISRW in the associated graph process \emph{does not} lie on the largest component. It therefore lies, with high probability, on either a unicyclic component or a tree. Denote by $N^{\rm uc}(t)$ the number of vertices in an \ER graph process that are in unicyclic connected components at time $t$, and recall from \eqref{eq:ERsizes} that $\mathscr{C}_{\rm max}^{\rm ER}(n,t)$ denotes the size of the largest component of an  \ER graph on $n$ vertices with $t$ edges. It follows that
\begin{equation}
\begin{aligned}
&\E\left[\P\left(\left[\Omega^{\rm tree}(t)\right]^\cmpl \cap \{\Mtddalt > t\} ~\Big|~ \left(\Pi_n(t)\right)_{s=0}^t\right)\right] \\
&=  \E \left[\P\left(\left[\Omega^{\rm tree}(t)\right]^\cmpl \cap \left\{\startvertex \not\in \maxgraph{\AG{\Pi_n}(t)}{}\right\} 
~\Big|~ \left(\Pi_n(t)\right)_{s=0}^t\right)\right] + o(1)\\
&\leq \E \left[ \min \left( \frac{N^{\rm uc}(t)}{n - |\mathscr{C}_{\rm max}^{\rm ER}(n,t)|}, 1 \right) \right] + o(1).
\end{aligned}
\end{equation}
From \cite[Theorem 5.11]{JLR2011} it follows that  $N^{\rm uc}(t) = O_{\P}(n^{2/3})$ uniformly in $t=cn$ with $c \in (\tfrac12,\infty)$. Since the size of the largest component is $\zeta(\tfrac{t}{n})n + o_{\P}(n)$ (recall Definition~\ref{def:ERstruct}) and the number of vertices is $n$, it follows that the number of vertices outside the largest component at time $t$ is $(1-\zeta(\tfrac{t}{n}))n + o_{\P}(n) = \Theta_{\P}(n)$, again uniformly in $t=cn$ with $c \in (\tfrac12,\infty)$. This gives
\begin{align}
\E \left[ \min \left( \frac{N^{\rm uc}(t)}{n - |\mathscr{C}_{\rm max}^{\rm ER}(n,t)|}, 1 \right) \right] = o(1).
\end{align}
Putting the above estimates together, we get
\begin{align}
\P\left(\Omega^ {\rm\sss (SC)}(t) \mid \Mtddalt > t \right) = 1 - o(1).
\end{align}

It remains to show that this estimate holds not only conditionally on $\{\Mtddalt > t\}$, but also conditionally on $\{\Mtddalt \geq t\}$. Given that the ISRW support lies on a single cycle before time $\Mtddalt$, at time $\Mtddalt$ this cycle merges with exactly one other cycle whose elements are represented by the vertices in the largest component of the associated graph process, and so the ISRW remains supported on a single permutation cycle. It therefore follows that
\begin{align}
\P\left(\Omega^{\sss\rm (SC)}(t) \mid \Mtddalt \geq t \right) = 1 - o(1).
\end{align}
\end{proof}

%%%

\subsection{Local mixing upon drop-down}
\label{sec:CFDP:localmix}

The main difference with the setting in Section~\ref{sec:CDP} is that each non-tree connected component of the associated graph process of CFDP may represent \emph{multiple} permutation cycles. We show that, after scaling of time, this fine structure is not felt because the distribution of the ISRW rapidly becomes uniform over the elements of the permutation represented by the vertices of the relevant connected component of the associated graph process. A consequence of this \emph{fast mixing} is the occurrence of the same phenomenon as observed for CDP, namely, at time $\Mtddalt$ there is a single drop in the total variation distance.

\paragraph{$\bullet$ Local mixing.} To formalise the arguments, we first introduce the notion of local mixing on the largest component of the associated graph process:

\begin{definition}[Local mixing time]
\label{def:tlocalmix}
Consider an ISRW with distribution $\ISRWdistr{}$ started from the element $\startvertex$ and running on top of CFDP $\Pi_n$ (recall Definition~\ref{def:cfdp}), and let $\AG{\Pi_n}$ be the associated graph process. For $\varepsilon \in (0,1)$, define the stopping time
\begin{equation}\label{eq:LMatDD}
\tlocalmix{\varepsilon} = \min \left\{t>\Mtddalt\colon\, \TVD\left(\ISRWdistr{}(t), \Unif{\maxgraph{\AG{\Pi_n}(t)}{}}\right)< \varepsilon\right\}.
\end{equation}
\hfill$\spadesuit$
\end{definition}

At time $\tlocalmix{\varepsilon}$, the ISRW is well mixed on the giant $\maxgraph{\AG{\Pi_n}(t)}{}$. In the following statements we illustrate and quantify the influence of large-enough permutation cycles on ISRW-mixing. Below we play with three parameters $n,\vep,\delta$ and take limits in the order $n\to\infty$, $\vep \downarrow 0$ and $\delta \downarrow 0$. We also play with a time scale $a_n$ satisfying $\lim_{n\to\infty} a_n = \infty$ and $a_n=o(n)$. Along the way we need some facts established in Appendices~\ref{apx:schramm} and \ref{app:sc} that require more stringent conditions on $a_n$, namely, $a_n = o(n^{1/26})$, respectively, $a_n = o(n^{1/3})$. We summarise this by saying that $a_n$ grows slowly enough.

We will often use the following \ER typicality event, which occurs with high probability:
\begin{definition}[\ER typicality event]\label{def:ER-event}
Take $\varepsilon_n$ such that $\varepsilon_n = o(1)$ and $\varepsilon_n = \omega(n^{-1/3})$. Define the following event
\begin{align}\label{ER-event}
\typicalevent{\sss(ER)}{t} &= \left\{ |\maxgraph{n,t}{ER}| = n(\zeta(\tfrac{t}{n}) - \varepsilon_n, \zeta(\tfrac{t}{n}) + \varepsilon_n, ) \right\} 
\cap \left\{ |\secgraph{n,t}{ER}| \leq n\varepsilon_n \right\}.
\end{align}
Note that this event is different from the event $\typicalevent{}{un}$ defined in \eqref{eq:ERtypical}.
\hfill$\spadesuit$ 
\end{definition}

\begin{definition}[Events $\bm{\mathcal{M}_1(\varepsilon, \delta)}$, $\bm{\mathcal{M}_2(\varepsilon)}$]\label{def:obsmix}
Denote by $\mathfrak{X}_1^{(n)}(t) $ the normalised size of the largest cycle at time~$t$ (see \eqref{eq:frakX}). Recall the event $\Omega^{\sss\rm (SC)}(t)$ from Lemma~\ref{lem:1cycle}, the \ER typicality event $\typicalevent{\rm\sss (ER)}{t}$ from Definition~\ref{def:ER-event} (which both occur with high probability for any $t =cn$ with $c \in (\tfrac12,\infty)$), and introduce the abbreviation $M = |\maxgraph{\AG{\Pi_n}(\Mtddalt)}{}|$. Define the events
\begin{equation}
\label{eq:Mevents}
\begin{aligned}
\mathcal{M}_1(\varepsilon, \delta) &= \big\{|\supp(\ISRWdistr{}(\Mtddalt )) | > \varepsilon M\big\} \cap \Omega^{\sss\rm (SC)}(\Mtddalt) \cap \typicalevent{\rm\sss (ER)}{\Mtddalt},\\
\mathcal{M}_2(\varepsilon) &= \big\{\exists\, t_{L} \in (\Mtddalt, \Mtddalt + a_n)\colon\,  \mathfrak{X}_1^{(n)}(t_L) > 1-\varepsilon^2 \big\}.
\end{aligned}
\end{equation}
\hfill$\spadesuit$
\end{definition}

\begin{lemma}[Mixing induced by a single large cycle]
\label{obs:mix}
Recall Definition~\ref{def:obsmix}. Let $(a_n)_{n\in\N}$ be such that $\lim_{n\to\infty} a_n = \infty$ slowly enough. Then
\begin{align}
\label{eq:lmix:inclusion}
\big\{ \tlocalmix{\varepsilon} \in (\Mtddalt, \Mtddalt + a_n) \big\} \supseteq \mathcal{M}_1(\varepsilon, \delta) \cap \mathcal{M}_2(\varepsilon).
\end{align}
Furthermore, on the event $\mathcal{M}_1(\varepsilon, \delta) \cap \mathcal{M}_2(\varepsilon)\cap \Omega^{\sss\rm (SC)}(\Mtddalt)$, there exists a $t_L\in (\Mtddalt, \Mtddalt + a_n)$ such that
\begin{equation}
\label{TV-bound-close-to-dropdown}
1-\frac{1}{n} |\maxgraph{\AG{\Pi_n}(t_L)}{}|-\vep\leq \dnTVD{t_L} \leq 1-\frac{1}{n} |\maxgraph{\AG{\Pi_n}(t_L)}{}|+\vep.
\end{equation}
\end{lemma}

\begin{proof}
Recall that the event $\Omega^{\sss\rm (SC)}(\Mtddalt) \subset \mathcal{M}_1(\varepsilon, \delta)$ implies that all the mass of the ISRW-distribution enters the giant component on a single cycle. Therefore the event $\mathcal{M}_1(\varepsilon, \delta)$ implies that
\begin{align}
\forall\,u &\in \supp(\ISRWdistr{}(\Mtddalt ))\colon  \ISRWdistr{u}(\Mtddalt) \leq \frac{1}{\varepsilon M},
%\forall\,w &\in \maxgraph{\AG{\Pi_n}(\Mtddalt)}{} \setminus \supp(\ISRWdistr{}(\Mtddalt ))\colon &&  \ISRWdistr{w}(\Mtddalt) = 0,
\label{bounds-mu-uniform}
\end{align}
where we recall that $M=|\maxgraph{\AG{\Pi_n}(\Mtddalt)}{}|$. The event $\mathcal{M}_1(\varepsilon, \delta)\cap\mathcal{M}_2(\varepsilon)$ indicates that a cycle of size at least $(1-\varepsilon^2)|\maxgraph{\AG{\Pi_n}(\Mtddalt)}{}|$ has appeared by time $\Mtddalt+a_n$. We denote this large permutation cycle by $\mathfrak{X}_{1}^{(n)}(t_L)$. This cycle necessarily contains some mass of the ISRW-distribution because, due to the event $\mathcal{M}_1(\varepsilon, \delta)$, the mass was initially spread out over a cycle that is larger than the region not covered by $\mathfrak{X}_{1}^{(n)}(t_L)$. We compute the effect of the event $\mathcal{M}_1(\varepsilon, \delta)\cap\mathcal{M}_2(\varepsilon)$ on the decay of the total variation distance. The worst possible scenario is when the $\varepsilon^2M$ elements not covered by $\mathfrak{X}_{1}^{(n)}(t_L)$ each carry mass $1/(\varepsilon M)$. Note that the definition of ISRW requires that the remaining mass is spread out uniformly over $\mathfrak{X}_{1}^{(n)}(t_L)$. A simple calculation (see Appendix~\ref{app:sc}) shows that, at time $t_L$ (introduced in the definition of the event $\mathcal{M}_2(\varepsilon)$) and for $n$ large enough,
\begin{align}
\label{TV-bound-onM2}
\TVD\Big(\ISRWdistr{}(t_L), \Unif{\maxgraph{\AG{\Pi_n}(t_L)}{}}\Big) < \varepsilon,
\end{align}
and \eqref{eq:lmix:inclusion} follows. 

To prove \eqref{TV-bound-close-to-dropdown} we use that, for probability mass functions $p=(p_x)_{x\in \mathcal{X}}$ and $q=(q_x)_{x\in \mathcal{X}}$,
\begin{equation}
\label{TV-representation}
\TVD(p,q)=\sum_{x\in \mathcal{X}} [\,p_x-(p_x\wedge q_x)\,].
\end{equation}

\medskip\noindent
$\bullet$ For the upper bound in \eqref{TV-bound-close-to-dropdown}, we use the triangle inequality to estimate
\begin{equation}
\label{est1}
\begin{aligned}
&\TVD\left(\ISRWdistr{}(t_L), \Unif{[n]}\right)\\
&\leq \TVD\left(\ISRWdistr{}(t_L), \Unif{\maxgraph{\AG{\Pi_n}(t_L)}{}}\right)
+\TVD\big(\Unif{\maxgraph{\AG{\Pi_n}(t_L)}{}}, \Unif{[n]}\big).
\end{aligned}
\end{equation}
Note that, by \eqref{TV-representation},
\begin{equation}
\label{est2}
\TVD\big(\Unif{\maxgraph{\AG{\Pi_n}(t_L)}{}}, \Unif{[n]}\big) = 1 - \frac{1}{n}|\maxgraph{\AG{\Pi_n}(t_L)}{}|,
\end{equation}
while, by \eqref{TV-bound-onM2},
\begin{equation}
\label{est3}
\TVD\left(\ISRWdistr{}(t_L), \Unif{\maxgraph{\AG{\Pi_n}(t_L)}{}}\right)\leq \vep.
\end{equation}
Combing \eqref{est2}--\eqref{est3} with \eqref{est1}, we get the upper bound in \eqref{TV-bound-close-to-dropdown}.

\medskip\noindent
$\bullet$
For the lower bound in \eqref{TV-bound-close-to-dropdown}, we note that
\begin{equation}
\label{est4}
\begin{aligned}
&\TVD\left(\ISRWdistr{}(t_L), \Unif{[n]}\right)\\
&\geq \TVD\left(\ISRWdistr{}(t_L), \Unif{\maxgraph{\AG{\Pi_n}(t_L)}{}}\right)
-\TVD\big(\Unif{\maxgraph{\AG{\Pi_n}(t_L)}{}}, \Unif{[n]}\big).
\end{aligned}
\end{equation}
Combining \eqref{est2}--\eqref{est3} with \eqref{est4}, we get the lower bound in \eqref{TV-bound-close-to-dropdown}.
\end{proof}

Before proceeding we make the following observation. Since $\frac{1}{n}\Mtddalt \convdist u^{\Downarrow}$ by Lemma \ref{lem:cfdp:tdddistr}, and $\prob(u^{\Downarrow}\leq \tfrac{1}{2} +\delta)=2\delta(1+o(1))$ for $\delta>0$ small enough, we note that $\prob(\tfrac{1}{n}\Mtddalt\leq \tfrac{1}{2}+\delta)$ can be made arbitrarily small by picking $\delta>0$ small. Furthermore, for $\delta>0$ small enough, $\maxgraph{\AG{\Pi_n}((\tfrac{1}{2}+\delta)n)}{}\geq \delta n$ with high probability, which follows from the properties of the \ER giant, specifically from the fact that $\zeta^\prime(\tfrac12)=2$. Thus, \eqref{bounds-mu-uniform} implies that, on the event
\begin{equation}
\Big\{\Mtddalt>(\tfrac{1}{2}+\delta)n,\ |\maxgraph{\AG{\Pi_n}((\tfrac{1}{2}+\delta)n)}{}|\geq \delta n\Big\},
\end{equation}
we have
\begin{equation}
\begin{aligned}
\forall\,u &\in \supp(\ISRWdistr{}(\Mtddalt ))\colon && \ISRWdistr{u}(\Mtddalt) \leq \frac{1}{\varepsilon \delta n},\\
%\forall\,w &\in \maxgraph{\AG{\Pi_n}(\Mtddalt)}{} \setminus \supp(\ISRWdistr{}(\Mtddalt ))\colon &&  \ISRWdistr{w}(\Mtddalt) = 0,\\
\forall\,w &\not\in \maxgraph{\AG{\Pi_n}(\Mtddalt)}{}\colon &&  \ISRWdistr{w}(\Mtddalt) = 0,
\label{bounds-mu-uniform-ext}
\end{aligned}
\end{equation}
where the last equality holds on the event $\Omega^{\sss\rm tree}(\Mtddalt)$. 

Note that for times $t\geq \Mtddalt$, just as for CDP,
\begin{equation}
\label{bounds-mu-uniform-ext-b}
\ISRWdistr{w}(t) = 0 \quad \forall\,w \not\in \maxgraph{\AG{\Pi_n}(t)}{} \qquad \forall\,t\geq \Mtddalt,
\end{equation}
since, by the construction of the associated graph process, the support of $\ISRWdistr{u}(\Mtddalt)$ always lies on a single connected component of the associated graph process. The uniform bounds above will prove to be essential below.

\medskip
Lemma~\ref{obs:mix} allows us to quantify the probability of $\varepsilon$-mixing after a single appearance of a cycle of size $(1-\varepsilon^2)M$:

\begin{proposition}
\label{prop:epsbound}
Fix $\delta>0$. Let $(a_n)_{n\in\N}$ be such that $\lim_{n\to\infty} a_n = \infty$ slowly enough. Then there exists a function $\varepsilon \mapsto f(\varepsilon)$ satisfying $\lim_{\varepsilon\downarrow0} f(\varepsilon) = 0$ such that
\begin{align}
\label{eq:lmix:bound}
\P\Big(\tlocalmix{\varepsilon} \not\in (\Mtddalt, \Mtddalt + a_n), \Mtddalt\geq (\tfrac{1}{2}+\delta)n\Big) 
\leq f(\varepsilon) + o(1), \qquad n\to\infty.
\end{align}
%Furthermore, there exists a threshold value $\varepsilon^\star > 0$ such that 
%\begin{align}
%\forall\, \varepsilon \in (0, \varepsilon^\star)\colon\, f(\varepsilon^\star) < \varepsilon^{1/3}.
%\end{align}
Consequently, on the event that $\Mtddalt\geq (\tfrac{1}{2}+\delta)n$, the conclusion of \eqref{TV-bound-close-to-dropdown} fails with probability at most $f(\varepsilon)$.
\end{proposition}

\begin{proof}
We will derive an upper bound for the probability of the event $(\mathcal{M}_1(\varepsilon, \delta)^\cmpl \cup \mathcal{M}_2(\varepsilon)^\cmpl)\cap \{\Mtddalt\geq (\tfrac{1}{2}+\delta)n\}$, which by \eqref{eq:lmix:inclusion} includes the event in the left-hand side of \eqref{eq:lmix:bound}. To do so, we will work with a further sub-event.

Denote  the number of vertices in $\maxgraph{\AG{\Pi_n}(t)}{}$ that are in cycles of size smaller than $\varepsilon n$ by $S(\varepsilon n, t)$. We use \cite[Lemma~2.4]{BKLM2019}, which states that for any $t> cn$ with $c>\tfrac12$ there exists a $C>0$ such that, for any $\varepsilon\in(0,1)$ and $n$ large enough,
\begin{align}\label{eq:S}
\E \left[ S(\varepsilon n, t) \right] < C\varepsilon \log(\tfrac{1}{\varepsilon})\,n.
\end{align}
Define the event 
\begin{equation}
\mathcal{M}_3(\varepsilon, t) = \{ S(\varepsilon n, t)  < \sqrt{\varepsilon}\,n\}.
\end{equation} 
Recall that the mass at $\Mtddalt$ enters the largest component of the associated graph process on a cycle that belongs to a uniform element. Trivially,
\begin{align}
\label{incl1}
&(\mathcal{M}_1(\varepsilon, \delta)^\cmpl  \cup \mathcal{M}_2(\varepsilon)^\cmpl)
\cap  \{\Mtddalt\geq (\tfrac{1}{2}+\delta)n\} \nonumber\\
&\qquad \subseteq 
\Big(\mathcal{M}_1(\varepsilon, \delta)^\cmpl  
\cup \mathcal{M}_2(\varepsilon)^\cmpl) \cap \mathcal{M}_3(\varepsilon, \Mtddalt) 
\cap  \{\Mtddalt\geq (\tfrac{1}{2}+\delta)n\}\Big)
\cup \mathcal{M}_3(\varepsilon, \Mtddalt)^\cmpl.
\end{align}
We estimate the probability of these events one by one. First, use the Markov inequality and \eqref{eq:S} to estimate, for $n$ large enough,
\begin{align}
\label{incl2}
\P\left(\mathcal{M}_3(\varepsilon, \Mtddalt)^\cmpl \right) \leq C\sqrt{\varepsilon} \log(\tfrac{1}{\varepsilon}).
\end{align}
Second, estimate
\begin{align}
\label{incl3}
\P \left(  \mathcal{M}_1(\varepsilon, \delta)^\cmpl \cap \mathcal{M}_3(\varepsilon, \Mtddalt) \right) 
\leq  \P \left( \mathcal{M}_1(\varepsilon, \delta)^\cmpl \mid  \mathcal{M}_3(\varepsilon, \Mtddalt) \right)
\leq \sqrt{\varepsilon},
\end{align}
with the last inequality following from the definition of $\mathcal{M}_1(\varepsilon, \delta)$. Third, the key estimate stated in Proposition~\ref{MaxFlux}, whose proof turns out to be rather delicate, yields that, for $\delta>0$ fixed,
\begin{align}
\label{incl4}
\P\left(\mathcal{M}_2(\varepsilon)^\cmpl, \Mtddalt\geq (\tfrac{1}{2}+\delta)n\right) = o(1).
\end{align}
Indeed, the key event that is estimated in Proposition~\ref{MaxFlux} is 	
\begin{equation}
\label{En-Prop-D10-def}
{\cal E}_n(c,\varepsilon, \kappa) = \big\{\exists\,(t_{k})_{k=1}^{\kappa}\in (cn, cn+a_n)\colon\, 
\mathfrak{X}_1^{(n)}(t_{k}-1)< 1-\varepsilon,\,\mathfrak{X}_1^{(n)}(t_{k})
\geq 1-\varepsilon\,\big\},
\end{equation}
which states that there are at least $\kappa\in \mathbb{N}$ times in the interval $(cn, cn+a_n)$ such that the size of the maximal cycle crosses  $(1-\varepsilon)$ upwards, i.e., $\mathfrak{X}_1^{(n)}(t_{k}) \geq 1-\vep$. Proposition~\ref{MaxFlux} states that ${\cal E}_n(c,\varepsilon, \kappa)$ occurs with high probability for all $c\in(1/2,\infty)$, $\kappa\in \mathbb{N}$ and $\vep>0$. We apply Corollary~\ref{cor:MaxFlux}, which is a consequence of Proposition~\ref{MaxFlux}, to obtain \eqref{incl4}.
% $c=\frac{1}{n}\Mtddalt$, which, on the event $\Mtddalt\geq (\tfrac{1}{2}+\delta)n$, is bounded away from $\tfrac{1}{2}$. \tc{\bf [Give the correct reference.]}

Combining \eqref{incl1}--\eqref{incl4}, we find that there exists a $C>0$ such that
\begin{equation}
\begin{aligned}
&\P\big((\mathcal{M}_1(\varepsilon, \delta)^\cmpl  \cup \mathcal{M}_2(\varepsilon)^\cmpl)\cap  \{\Mtddalt\geq (\tfrac{1}{2}+\delta)n\}\big)\\
&\qquad \leq \P\left(\mathcal{M}_3(\varepsilon, \Mtddalt)^\cmpl \right)
+\P \left(  \mathcal{M}_1(\varepsilon, \delta)^\cmpl \cap \mathcal{M}_3(\varepsilon, \Mtddalt)  \right)
+\P \left( \mathcal{M}_2(\varepsilon)^\cmpl, \Mtddalt\geq (\tfrac{1}{2}+\delta)n\right)\\
&\qquad \leq C\sqrt{\varepsilon} \log(\tfrac{1}{\varepsilon})+\sqrt{\varepsilon} + o(1) < \varepsilon^{1/3}
\end{aligned}
\end{equation}
for $\varepsilon$ small enough, which in turn decays to $0$ as $\varepsilon\to 0$.
\end{proof}

\paragraph{$\bullet$ Adaptation of Lemma~\ref{obs:mix} and Proposition~\ref{prop:epsbound}.}  Finally, we adapt Lemma~\ref{obs:mix} and Proposition~\ref{prop:epsbound}. Note that Lemma~\ref{obs:mix} is true at time $t=cn$ when we replace the events $\mathcal{M}_1(\vep, \delta), \mathcal{M}_2(\vep)$ by (compare with \eqref{eq:Mevents})
\begin{equation}
\begin{aligned}
\mathcal{M}^\prime_1(cn, \varepsilon, \delta) 
&= \big\{|\supp(\ISRWdistr{}(cn))| > \varepsilon  n\big\}\cap \Omega^{\sss\rm (SC)}(\Mtddalt),\\
\mathcal{M}^\prime_2(cn, \varepsilon,\delta) 
&= \big\{\exists\,t_{L} \in (cn, cn + a_n)\colon\,\mathfrak{X}_1^{(n)}(t_L) > 1-\tfrac{\varepsilon^{2}}{\delta} \big\}.
\label{Mdefs}
\end{aligned}
\end{equation}
Here, we recall the event $\Omega^{\sss\rm (SC)}(t)$ from Lemma~\ref{lem:1cycle} (which occurs with high probability for any $t =cn$ with $c \in (\tfrac12,\infty)$), and the extra factor $1/\delta$ is added to accommodate the extra factor $1/\delta$ in the first line of \eqref{bounds-mu-uniform-ext}. It remains to redo the calculations in the proofs of Lemma~\ref{obs:mix} and Proposition~\ref{prop:epsbound} with these modified events. Take $t=cn$  with $c \in (\tfrac12,\infty)$, and define
\begin{align}
\tlocalmix{\varepsilon}(t) = \min \left\{s>t\colon\, \TVD\big(\ISRWdistr{}(s), \Unif{\maxgraph{t}{}}\big)< \varepsilon\right\}.
\end{align}

We start by adapting Lemma~\ref{obs:mix}:

\begin{lemma}[Mixing induced by a single large cycle]
\label{obs:mix-cn}
Let $(a_n)_{n\in\N}$ be such that $\lim_{n\to\infty} a_n = \infty$ slowly enough, and let $c\in(1/2, \infty)$. Then, for any $\delta\in (0, c-\tfrac12)$,
\begin{align}
\label{eq:lmix:inclusionalt}
&\big\{\tlocalmix{\varepsilon}(cn) \in (cn, cn+ a_n) \big\}  \cap \{(\tfrac{1}{2}+\delta)n\leq\Mtddalt\leq cn\}\\
&\qquad\supseteq \mathcal{M}^\prime_1(cn, \varepsilon, \delta)
\cap \mathcal{M}^\prime_2(cn, \varepsilon,\delta) \cap \{(\tfrac{1}{2}+\delta)n\leq \Mtddalt\leq cn-a_n\}.\nonumber
\end{align}
Furthermore, on the event $\mathcal{M}_1(\varepsilon, \delta) \cap \mathcal{M}_2(\varepsilon)\cap \{(\tfrac{1}{2}+\delta)n\leq \Mtddalt\leq cn-a_n\} \cap \Omega^{\sss\rm (SC)}(\Mtddalt)$ there exists a~$t_L\in (cn, cn+a_n)$ such that
\begin{equation}
\label{TV-bound-after-dropdown}
1-\frac{1}{n} |\maxgraph{\AG{\Pi_n}(t_L)}{}|-\vep \leq \dnTVD{t_L}\leq 1-\frac{1}{n} |\maxgraph{\AG{\Pi_n}(t_L)}{}|+\vep.
\end{equation}
\end{lemma}

\begin{proof}
The main ingredient in the proof of Lemma~\ref{obs:mix} was \eqref{bounds-mu-uniform}. Recall the extension of \eqref{bounds-mu-uniform} in \eqref{bounds-mu-uniform-ext}. With \eqref{bounds-mu-uniform-ext} in hand, we can simply follow the proof of Lemma~\ref{obs:mix}.
\end{proof}

We continue by adapting Proposition~\ref{prop:epsbound}:

\begin{proposition}
\label{prop:DTepsbound}
 Let $(a_n)_{n\in\N}$ be such that $\lim_{n\to\infty} a_n = \infty$ slowly enough, and let $c>\tfrac{1}{2}$. Then,  for any $\delta\in (0, c-\tfrac12)$, with $\varepsilon \mapsto f(\varepsilon)$ as in Proposition \ref{prop:epsbound}, 
\begin{align}
\label{eq:lmix:boundalt}
\P\Big(\tlocalmix{\varepsilon}(cn) \not\in (cn, cn + a_n),(\tfrac{1}{2}+\delta)n\leq \Mtddalt  \leq cn-a_n\Big) 
\leq f(\varepsilon) + o(1), \qquad n\to\infty.
\end{align}
Consequently, on the event that $(\tfrac{1}{2}+\delta)n\leq \Mtddalt\leq cn-a_n$, the conclusion of \eqref{TV-bound-after-dropdown} fails with probability at most $f(\varepsilon)$.
\end{proposition}

\begin{proof}
We follow the proof of Proposition \ref{prop:epsbound}, which relies on the inclusion in Lemma \ref{obs:mix}. Instead, we now rely on the inclusion in Lemma \ref{obs:mix-cn}. Recall from the proof of Lemma \ref{obs:mix} that $S(\varepsilon n, t)$ denotes the number of vertices in $\maxgraph{\AG{\Pi_n}(t)}{}$ that are in cycles of size smaller than $\varepsilon n$, and that, by \eqref{eq:S}, $\E \left[ S(\varepsilon n, t) \right] < C\varepsilon \log(\tfrac{1}{\varepsilon})\,n$.

Recall $\tlocalmix{\varepsilon}$ from \eqref{eq:LMatDD}. Define the event 
\begin{equation}
\mathcal{M}_3^\prime(\varepsilon) = \{ \tlocalmix{\varepsilon}\in (\Mtddalt, \Mtddalt + a_n)\}.
\end{equation}
Trivially,
\begin{equation}
\begin{aligned}
&\Big(\mathcal{M}^\prime_1(cn, \varepsilon, \delta)^\cmpl  \cup \mathcal{M}^\prime_2(cn, \varepsilon,\delta)^\cmpl\Big)
\cap \Big\{(\tfrac{1}{2}+\delta)n\leq \Mtddalt  \leq cn-a_n\Big\}\\
&\qquad \subseteq 
\Big(\mathcal{M}_3^\prime(\varepsilon) \cap \Big(\mathcal{M}^\prime_1(cn, \varepsilon, \delta)^\cmpl  \cup \mathcal{M}^\prime_2(cn, \varepsilon,\delta)^\cmpl\Big) 
\cap \Big\{(\tfrac{1}{2}+\delta)n\leq \Mtddalt  \leq cn-a_n\Big\} \Big)\\
&\qquad\qquad \cup \Big(\mathcal{M}_3^\prime(\varepsilon)^\cmpl\cap \{(\tfrac{1}{2}+\delta)n\leq \Mtddalt  \leq cn-a_n\}\Big).
\label{split-Prop-38}
\end{aligned}
\end{equation}
We estimate the probability of these events one by one. First, for $n$ large enough,
\begin{align}
\label{eq:3:bound1}
\P\Big(\mathcal{M}_3^\prime(\varepsilon)^\cmpl, (\tfrac{1}{2}+\delta)n\leq \Mtddalt  \leq cn-a_n \Big) \leq \P\Big(\mathcal{M}_3^\prime(\varepsilon)^\cmpl, (\tfrac{1}{2}+\delta)n\leq \Mtddalt \Big) \leq f(\varepsilon) + o(1),
\end{align}
where the last inequality follows from Proposition \ref{prop:epsbound}. Second, if $\tlocalmix{\varepsilon}\in (\Mtddalt, \Mtddalt + a_n)$ and $\Mtddalt\geq (\tfrac{1}{2}+\delta)n$, then
\begin{align}
\label{keybd}
\P\Big(\mathcal{M}^\prime_1(cn, \varepsilon, \delta)^\cmpl ~\Big|~  \mathcal{M}^{\prime}_3(\varepsilon),(\tfrac{1}{2}+\delta)n 
\leq \Mtddalt  \leq cn-a_n\Big) = 0.
\end{align}
Indeed, $\mathcal{M}^\prime_1(cn, \varepsilon, \delta)^\cmpl$ and $\Mtddalt  \leq cn-a_n$ imply that $|\supp(\ISRWdistr{}(\Mtddalt+a_n))| \leq \varepsilon  n$. By an application of \eqref{TV-representation} with $\mathcal{X}=[n]$, $p=\Unif{[n]}$ (for which $p_v=\frac{1}{n}$ for all $v\in [n]$) and $q_v=\ISRWdistr{}(\Mtddalt+a_n)$, this implies that
\begin{equation}
\dnTVD{\Mtddalt+a_n}\geq 1-\vep.
\end{equation}
However, the latter is incompatible with \eqref{TV-bound-close-to-dropdown} when $\Mtddalt\geq (\tfrac{1}{2}+\delta)n$, since
\begin{equation}
\begin{aligned}
\dnTVD{\Mtddalt+a_n}&\leq \dnTVD{t_L}\leq 1-\frac{1}{n} |\maxgraph{\AG{\Pi_n}(t_L)}{}|+\vep \\
&\leq 1-\frac{1}{n} |\maxgraph{\AG{\Pi_n}((\tfrac{1}{2}+\delta)n)}{}|+\vep \\
&\leq 1 - 2\delta +o(1) + \vep<1-\vep,
\end{aligned}
\end{equation}
where the second inequality uses the definition of the event $\mathcal{M}^{\prime}_3(\vep)$, and the last inequality is valid for $\vep$ small enough depending on $\delta$. Third, apply the key estimate stated in Proposition~\ref{MaxFlux} (see the explanation below \eqref{En-Prop-D10-def}), to get
\begin{align}
\label{eq:3:bound3}
\P\Big(\mathcal{M}^\prime_2(cn, \varepsilon,\delta),(\tfrac{1}{2}+\delta)n\leq \Mtddalt  \leq cn-a_n\Big) = o(1).
\end{align}

Combining \eqref{eq:lmix:inclusionalt}, \eqref{split-Prop-38}--\eqref{keybd} and \eqref{eq:3:bound3}, we obtain
\begin{equation}
\begin{aligned}
&\P\Big(\tlocalmix{\varepsilon}(cn) \not\in (cn, cn+ a_n), (\tfrac{1}{2}+\delta)n\leq \Mtddalt  \leq cn-a_n\Big)\\
&\quad\leq \P\Big((\mathcal{M}^\prime_1(cn, \varepsilon, \delta)^\cmpl \cup \mathcal{M}^\prime_2(cn,\varepsilon,\delta)^\cmpl)
\cap  \{(\tfrac{1}{2}+\delta)n\leq \Mtddalt  \leq cn-a_n \} \Big) \\
&\quad \leq \P\Big(\mathcal{M}_3^\prime(\varepsilon)^\cmpl, (\tfrac{1}{2}+\delta)n\leq \Mtddalt  \leq cn-a_n \Big)
+ \P\Big(\mathcal{M}^\prime_1(cn, \varepsilon, \delta)^\cmpl\cap \mathcal{M}_3^\prime(\varepsilon),(\tfrac{1}{2}+\delta)n
\leq \Mtddalt  \leq cn-a_n\Big)\\
&\quad\qquad + \P\Big(\mathcal{M}^\prime_2(cn,\varepsilon,\delta)^\cmpl,(\tfrac{1}{2}+\delta)n\leq \Mtddalt  \leq cn-a_n\Big)\\
&\quad \leq f(\vep) + o(1) + 0 + o(1) = f(\vep)+o(1),
\end{aligned}
\end{equation}
where the first inequality uses the inclusion in Lemma \ref{obs:mix-cn}.
\end{proof}

%%%

\subsection{Mixing profile}
\label{sec:CFDP:mixprofile}

Like in the case of CDP, the results on the mixing profile are established in two steps. First we establish pointwise convergence, afterwards we extend to process convergence. The following lemma settles Theorem~\ref{thm:frag:main}(2):

\begin{lemma}[Pointwise convergence of the mixing profile for ISRW on CFDP]
\label{lem:CFDP:profile}
Uniformly in $\startvertex \in [n]$,
\begin{equation}
\label{eq:mixing-stat}
\dnTVD{sn} \stackrel{d}{\to} 1-\zeta(s)W(s), \qquad s\in[0,\infty),
\end{equation}   
where $W(s)\sim\mathsf{Bernoulli}(\zeta(s))$.
\end{lemma}

\begin{proof}
Fix $s\in[0,\infty)$ and split the random variable $\dnTVD{sn} -1$ as
\begin{align}
\label{split-dTV-before-after}
\dnTVD{sn} - 1 &= [\dnTVD{sn} - 1] \left( \1_{\{ \Mtddalt > sn\}} + \1_{\{ \Mtddalt \leq sn\}} \right)\\
&\qquad\times \left( \1_{ \typicalevent{\sss(ER)}{sn}}  
+ \1_{ [\typicalevent{\sss(ER)}{sn}]^{\rm c}} \right) \left( \1_{\Omega^{\rm tree}(sn)}  
+ \1_{ [\Omega^{\rm tree}(sn)]^{\rm c}} \right),\nonumber
\end{align}
where the event $\Omega^{\rm tree}(sn)$ is defined in the proof of Lemma~\ref{lem:1cycle}, and we recall the \ER typicality event (see \eqref{ER-event})
\begin{align}
\typicalevent{\sss(ER)}{t} &= \left\{ |\maxgraph{n,t}{ER}| = n(\zeta(\tfrac{t}{n}) - \varepsilon_n, \zeta(\tfrac{t}{n}) + \varepsilon_n, ) \right\} 
\cap \left\{ |\secgraph{n,t}{ER}| \leq n\varepsilon_n \right\}.
\end{align} 
Because $\typicalevent{\sss(ER)}{sn}$ and $\Omega^{\rm tree}(sn)$ both occur with high probability, the terms containing the indicators $\1_{ [\typicalevent{ER}{sn}]^{\rm c}}$ and $\1_{ [\Omega^{\rm tree}]^{\rm c}}$  converge to 0 in probability, and hence
\begin{align}
\label{eq:simplified}
\dnTVD{sn} - 1 &= [\dnTVD{sn} - 1] \left( \1_{\{ \Mtddalt > sn\}} + \1_{\{ \Mtddalt \leq sn\}} \right) \1_{ \typicalevent{\sss(ER)}{sn}} \1_{\Omega^{\rm tree}(sn)} + o_{\P}(1).
\end{align}

To deal with the first term in \eqref{eq:simplified}, we note that 
\begin{align}
\label{eq:mix:limit1}
&[\dnTVD{sn} -1] \1_{\{\Mtddalt > sn\}} \1_{ \typicalevent{\sss(ER)}{sn}} \1_{\Omega^{\rm tree}(sn)} \\ 
&\qquad\deq \left[\left(1 - \frac{O_{\P}(n^{2/3})}{n}\right) -1\right] \1_{\{\Mtddalt > sn\}} 
\1_{\typicalevent{\sss(ER)}{sn}} \1_{\Omega^{\rm tree}(sn)} \convprob 0,\nonumber
\end{align}
since, on the above events, the distribution of ISRW is uniform over a single permutation cycle \emph{outside} of the largest component of the associated graph process, whose size is $O_{\P}(n^{2/3})$. 

To deal with the second term in \eqref{eq:simplified}, which only contributes when $s>\tfrac{1}{2}$, we use Lemma \ref{lem:cfdp:tdddistr}. For $\delta>0$ sufficiently small and $a_n$ as in Proposition \ref{prop:DTepsbound}, we split
\begin{equation}
\label{eq:mix:3waysplit}
\1_{\{ \Mtddalt \leq sn\}}  = \1_{\big\{(\tfrac{1}{2}+\delta) n\leq \Mtddalt \leq sn-a_n\big\}} + \1_{\big\{sn-a_n< \Mtddalt \leq sn\big\}}
+ \1_{\big\{\Mtddalt<(\tfrac{1}{2}+\delta) n\big\}}.
\end{equation}
We rely on \eqref{TV-bound-after-dropdown} in Lemma \ref{obs:mix-cn}, which holds with high probability due to Proposition \ref{prop:DTepsbound}. (It is here that we need $ \Mtddalt\geq (\tfrac{1}{2}+\delta) n$, since this appears as an assumption in Proposition \ref{prop:DTepsbound}.) We claim that
\begin{equation}
\label{convergence-supremum-Cmax}
\sup_{t \in \N} \Big|\frac{1}{n} |\maxgraph{\AG{\Pi_n}(t)}{}|-\zeta(\tfrac{t}{n})\Big|=\op(1).
\end{equation}
Indeed,  \eqref{convergence-supremum-Cmax} holds because $\frac{1}{n} |\maxgraph{\AG{\Pi_n}(sn)}{}|\convp \zeta(s)$ for all $s>0$ fixed, $s\mapsto \zeta(s)$ is non-decreasing and continuous, and $s\mapsto \frac{1}{n} |\maxgraph{\AG{\Pi_n}(sn)}{}|$ is non-decreasing. By \eqref{TV-bound-after-dropdown} and \eqref{convergence-supremum-Cmax}, we obtain, for all $s>\tfrac{1}{2}+\delta$ and on the event $\{\tlocalmix{\varepsilon}(cn)\in (sn, sn + a_n),(\tfrac{1}{2}+\delta)n\leq \Mtddalt  \leq sn-a_n\}$, that there exists a $t_L\in (sn, sn+a_n)$ such that
\begin{equation}
\label{TV-bound-after-dropdown-ext}
1-\zeta(\tfrac{t_L}{n})-\vep-\op(1) \leq \dnTVD{t_L}\leq 1-\zeta(\tfrac{t_L}{n})+\vep+\op(1).
\end{equation}
Since $\vep>0$ is arbitrary, we conclude that, on the event $\{\tlocalmix{\varepsilon}(sn)\in (sn, sn + a_n),(\tfrac{1}{2}+\delta)n\leq \Mtddalt  \leq sn-a_n\}$, there exists a~$t_L\in (sn, sn+a_n)$ such that
\begin{equation}
\label{TV-bound-after-dropdown-ext2}
\dnTVD{t_L}=1-\zeta(\tfrac{t_L}{n})+\op(1).
\end{equation}
Since the above is true for all $s>\tfrac{1}{2}+\delta$, and $t\mapsto \dnTVD{t}$ is non-increasing, while $s\mapsto 1-\zeta(s)$ is non-increasing and continuous, \eqref{TV-bound-after-dropdown-ext2} implies that, for all $c>\tfrac{1}{2}+\delta$ and on the event $\{(\tfrac{1}{2}+\delta)n\leq \Mtddalt  \leq sn-a_n\}$,
\begin{equation}
\label{TV-bound-after-dropdown-ext3}
\dnTVD{sn}=1-\zeta(s)+\op(1).
\end{equation}
Since $\1_{ \typicalevent{\sss(ER)}{sn}} \1_{\Omega^{\rm tree}(sn)}\convp 1$, it follows that
\begin{equation}
[1-\dnTVD{sn}]\1_{\{(\tfrac{1}{2}+\delta) n\leq \Mtddalt \leq sn-a_n\}}
 \1_{ \typicalevent{\sss(ER)}{sn}} \1_{\Omega^{\rm tree}(sn)}  -\zeta(s)\1_{\{\tfrac{1}{2}+\delta \leq  \tfrac{1}{n}\Mtddalt \leq s-\tfrac{a_n}{n}\}}=\op(1).
\end{equation}
By Lemma \ref{lem:cfdp:tdddistr} and Slutsky's theorem, we thus conclude that
\begin{equation}
\label{eq:mix:limit2}
[1-\dnTVD{sn}]\1_{\{ \Mtddalt \leq sn-a_n\}} \1_{ \typicalevent{\sss(ER)}{sn}} \1_{\Omega^{\rm tree}(sn)} \convd \zeta(s)
\1_{\{\tfrac{1}{2}+\delta\leq u^{\Downarrow} \leq s\}}.
\end{equation}

Finally, by Lemma \ref{lem:cfdp:tdddistr},
\begin{equation}
\label{eq:mix:limit3}
\prob(sn-a_n< \Mtddalt \leq sn)+\prob(\tfrac{1}{n}\Mtddalt<\tfrac{1}{2}+\delta) \rightarrow \zeta(\tfrac{1}{2}+\delta),
\end{equation}
which tends to 0 as $\delta\downarrow 0$. The claim in \eqref{eq:mixing-stat} follows by combining \eqref{eq:simplified}--\eqref{eq:mix:3waysplit} and \eqref{eq:mix:limit2}--\eqref{eq:mix:limit3}.
\end{proof}

Finally, an argument based on monotonicity and a growing sequence of compact intervals settles Theorem~\ref{thm:frag:main} and concludes this section.

\begin{proof}[Proof of Theorem~\ref{thm:frag:main}(2)]
Observe that for every $n\in\N$, any realisation of $\dnTVD{\cdot}$ is a monotone c\`adl\`ag path on the set $[0,\infty)$. The pointwise convergence proven in Lemma~\ref{lem:CFDP:profile} implies, by \cite[Corollary 12.5.1.]{W2002}, pathwise convergence in the Skorokhod $M_1$-topology on any compact set $[0, t]$ such that $t>0$ is with probability 1 a continuity point of the limiting process. But the latter is true for any $t>0$ because the limiting process has almost surely \emph{one} point of discontinuity, whose position is distributed randomly according to the non-atomic distribution identified in Lemma~\ref{lem:cfdp:tdddistr}. Taking a sequence $(t_k)_{k\in\N}$ of such continuity points with $t_k \to\infty$ as $k\to\infty$, we also obtain pathwise convergence in the Skorokhod $M_1$-topology on the non-compact set $[0,\infty)$. For details, see \cite[p.~414]{W2002}.
\end{proof}

%%%%%%%%%%%% APPENDICES %%%%%%%%%%%%%%%

\appendix

%%%%%%% APPENDIX A %%%%%%%%%%%%%%%%%%%%%%%%%%

\section{Infinite-speed random walk as limit of finite-speed random walk}
\label{apx:limit}

It is natural to ask what the relation is between an infinite-speed random walk (ISRW) and a finite-speed standard random walk.

\begin{definition}[Finite-speed random walk on $\Pi_n$]
Let $\Pi_n$ be a dynamic permutation. Fix an element $\startvertex$. Recall that $\cycle{0}{\startvertex}$ is the cycle of $\Pi_n(0)$ that contains $\startvertex$. Pick $\rho\in\N$ as the speed ratio between the evolution of the random walk and the random graph. More formally, let $\Pi_n$ satisfy the condition
\begin{equation}
\forall\, n \in \N_0\ \forall\, i\, \in \intv[0]{\rho-1}\colon\, \Pi_n(\rho n) \equiv \Pi_n(\rho n+i).
\end{equation}
Denote by $\Pi_n(t)[i]$ the image of the element $i$ under the permutation $\Pi_n(t)$. The finite-speed random walk is the random process $(Y_n^{\startvertex}(t))_{t\in\N_0}$ on $[n]$ that has $Y_n^{\startvertex}(0) = \startvertex$ and, for any time $t>0$, 
\begin{equation}
\P( {Y_n^{\startvertex}}(t) = i\,|\, {Y_n^{\startvertex}}(t-1) = j) =
\begin{cases}
1, &\text{if } \Pi_n(t)[i] = j\text{ and }\Pi_n(t)[j] = i \text{ (which is the case for 1- and 2-cycles)}, \\
\frac{1}{2}, &\text{if } \Pi_n(t)[i] \neq j \text{ and either } \Pi_n(t)[j] = i \text{ or } \Pi^{-1}_n(t)[j] = i,\\
0, &\text{otherwise}.
\end{cases}
\end{equation}
i.e., a simple symmetric random walk on the elements of any given cycle, where the underlying permutation changes once every $\rho$ steps. \hfill$\spadesuit$
\end{definition}

With this definition in hand, we can explain what we mean by saying that ISRW has infinite speed:

\begin{proposition}[ISRW arises as limit]
Let $\Pi_n$ be a dynamic permutation. Recall $(\ISRWdistr{}(t))_{t\in\N_0}$ defined in Definition~\ref{def:ISRW}, and let $Y_n^{\startvertex} = \{ {Y_n^{\startvertex}}(t)\}_{t\in\N_0}$ be a finite-speed random walk with speed ratio $\rho$ on the dynamic permutation $\Pi^Y_n$, with distribution $\distr{Y_n^{\startvertex}}{}(t)$ at time $t$. Note that $ \startvertex \in [n]$ is the starting element for both the ISRW and the finite-speed random walk. Finally, suppose that the following relation holds between $\Pi_n$ and $\Pi^Y_n$:
\begin{equation}
\forall\, n\in\N_0\,\, \forall\, i\, \in \intv[0]{\rho-1}\colon\, \Pi_n(n+i) \equiv \Pi^Y_n(n).
\end{equation}
Then, with high probability, for any $i\in\N$,
\begin{equation}
d_{\rm \scriptstyle TV}\left(\ISRWdistr{}(i),\distr{Y_n^{\startvertex}}{}(\rho i)\right) \stackrel{\rho \to \infty}{\to} 0.
\end{equation}
\end{proposition}

\begin{proof}
It is well-known that a simple random walk on a circle of size $m$ mixes in $O(m^2)$ steps. Pick $\varepsilon_{RW} > 0$, and pick $\rho(\varepsilon_{RW})$ larger than the $\varepsilon_{RW}$-mixing time of a simple random walk on a circle of size $n$. Denote by $Y_n^{\startvertex}({\rho(\varepsilon_{RW}) i})\restrict{\xcycle{\rho i}{v}{Y}}$ the restriction of $Y_n^{\startvertex}({\rho(\varepsilon_{RW}) i})$ to the cycle $\xcycle{\rho i}{v}{Y} \in \Pi_n^Y(\rho i)$, and let $c\,\Unif{S}$ be the uniform distribution multiplied component-wise by the constant $c$. Then
\begin{equation}
d_{\rm \scriptstyle TV}\left(Y_n^{\startvertex}({\rho(\varepsilon_{RW}) i})\restrict{\xcycle{\rho i}{v}{Y}}, 
\frac{1}{\sum_{v \in \xcycle{i}{v}{}} \ISRWdistr{v}(i)} \Unif{[|\xcycle{i}{v}{}|}]\right) < \varepsilon_{RW}.
\end{equation}
Since, by definition, 
\begin{equation}
\ISRWdistr{}(i)\restrict{\xcycle{i}{v}{}} = \frac{1}{\sum_{v \in \xcycle{i}{v}{}} \ISRWdistr{v}(i)}\,\Unif{[|\xcycle{i}{v}{}|]},
\end{equation} 
the claim follows. 
\end{proof}

%%%%%%%% APPENDIX B %%%%%%%%%%%%%%%%%%%
 
\section{Normalisation of the jump-time distribution}
\label{apx:norm}

In this appendix we return to Definition~\ref{def:ERstruct} and show that $\lim_{s\to\infty} \phi(s)=1$, i.e., the laws of the jump-down times in Theorems~\ref{thm:main}--\ref{thm:frag:main} are normalised. The proof amounts to showing that 
\begin{equation}
\label{eq:integral-zeta}
\int_{\tfrac12}^\infty \dd s\,[1-\zeta^2(s)] = \tfrac12,
\end{equation}
where $\zeta(u)$ is the unique solution of the equation
\begin{equation}
\label{eq:def-zeta}
\ee^{-2s\zeta(s)} = 1 - \zeta(s), \qquad s\in\left[\tfrac12,\infty\right).
\end{equation}
To evaluate \eqref{eq:integral-zeta}, introduce a new variable $u = \zeta(s)$ and write
\begin{align}
\int_{\tfrac12}^\infty \dd s\,[1-\zeta^2(s)] = \int_{0}^{1} \dd u\,\frac{\dd s}{\dd u}\,(1-u^2).
\end{align}
To compute the Jacobian $\frac{\dd s}{\dd u}$, we take the logarithm of \eqref{eq:def-zeta} and differentiate:
\begin{align}
s = - \frac{1}{2u} \log\left(1-u\right), \qquad \frac{\dd s}{\dd u} = \frac{1}{2u(1-u)} + \frac{\log(1-u)}{2u^2}.
\end{align}
This gives
\begin{equation}
 \label{eq:int-pows} 
\begin{aligned}
\int_{0}^{1} \dd u\,\frac{\dd s}{\dd s}\,(1-u^2)
&=  \int_{0}^{1} \dd u\,(1+u) \frac{1}{2u} \left[ 1 - \frac{1-u}{u} \sum \limits_{k\in\N} \frac{u^k}{k}\right]\\
&= \frac{1}{2} \int_{0}^{1} \dd u\, \left[ (1+u) - (1-u^2) \sum_{l\in\N_0} \frac{u^l}{l+2}\right]
= \frac{1}{2} \int_0^1 \dd u \sum _{l\in\N_0} c_l u^l
\end{aligned}
\end{equation}
with the coefficients $c_l$ given by
\begin{equation}
c_l = \delta_{0l} + \delta_{1l} -\frac{1}{l+2} \1_{\{l\geq0\}} + \frac{1}{l}\1_{\{l\geq2\}}, \qquad l \in \N_0.
\end{equation}
The sum has radius of convergence $1$, and so term by term integration gives
\begin{equation}
\begin{aligned}
\frac{1}{2} \int_0^1 \dd u \sum \limits_{l\in\N_0} c_l u^l
&= \frac{1}{2} \left[ \frac{1}{2} \times 1 + \frac{2}{3} \times \frac{1}{2} + \sum_{l\geq2} \frac{1}{l+1} \left( \frac{1}{l} - \frac{1}{l+2}\right) \right] \\
&= \frac{1}{4} +\frac{1}{6} + \frac{1}{2} \sum_{l\geq2} \left[ \frac{1}{l(l+1)} - \frac{1}{(l+1)(l+2)} \right]
= \frac{1}{4} + \frac{1}{6} + \frac{1}{12} = \frac{1}{2},
\end{aligned} 
\end{equation}
as required.

 %%%%%%%% APPENDIX C %%%%%%%%%%%%%%%%%%%

\section{Cycle structure of coagulative-fragmentative dynamic permutations and Schramm's coupling}
\label{apx:schramm}

In  \cite{S2005}, Schramm introduced a remarkable coupling to study the cycle structure of CFDP. The coupling is realised between the part of the dynamic permutation that is supported on elements that lie in the largest component of the associated graph process and an independent $\mathsf{PoiDir}(1)$-sample, both seen as partitions of the unit interval $[0,1]$ evolving via coagulative-fragmentative dynamics. The main idea of the coupling is to evolve both partitions according to the same dynamics, but to rearrange and \emph{match} components of the partition that become close in size.  A detailed account of Schramm's coupling is given in \cite[Section 5]{BKLM2019}, to which we refer the reader. Below we provide a short overview of the main features.

Appendix~\ref{sec:summary} provides a short summary of Schramm's coupling and collects a number of estimates that are needed later on. Appendix~\ref{sec:recur} proves a key proposition (Proposition~\ref{MaxFlux} below) stating that large cycles  recur arbitrarily often in any time interval of diverging length. 

%%%

\subsection{Short summary of Schramm's coupling}
\label{sec:summary}

Schramm's coupling is defined for a pair $\vec{Y}(t),\vec{Z}(t)$ of evolving partitions of the unit interval $[0,1]$ into countable many subintervals. The dynamics of the evolution is chosen such that the Poisson-Dirichlet distribution $\mathsf{PoiDir}(1)$ is invariant under the dynamics (see e.g.\ \cite{T2000}). As we will see below, such a dynamics corresponds to the coagulative-fragmentative dynamics considered in this paper. 

%%%%%%%%%%%%%%%%%%%%%%%%%%%%%%%%%%%%%%%%%%%%%%%%%%%%%%%%%%%%
\begin{figure}[htbp]
\centering
\includegraphics[width=0.5\linewidth]{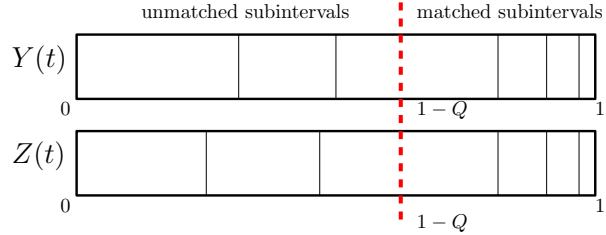}
\caption{\small Illustration of the setup for Schramm's coupling defined in Definition~\ref{def:asc}.}
\label{fig:asc}
\end{figure}
%%%%%%%%%%%%%%%%%%%%%%%%%%%%%%%%%%%%%%%%%%%%%%%%%%%%%%%%%%%%

\begin{definition}[Abstract version of Schramm's coupling]
\label{def:asc}
Take two partitions $\vec{Y}(0), \vec{Z}(0)$ of the unit interval $[0,1]$ into countably many subintervals. A single step of the coupling proceeds as follows:
\begin{enumerate}
\item 
If there is a subinterval $a \in \vec{Y}(0)$ that has the same length as some other subinterval $b\in \vec{Z}(0)$, then declare $a$ and $b$ to be \emph{matched}. (Note that the relation of being matched is symmetric.) 
\item 
Reorder the subintervals within $\vec{Y}(0)$ and $\vec{Z}(0)$ as follows:
\begin{enumerate}
\item 
Let $Q$ be the total length of all the matched intervals. In both $\vec{Y}(0)$ and $\vec{Z}(0)$, place within $(1-Q,1]$ all the \emph{matched} subintervals, ordered by their size such that the longest matched subinterval is on the left.
\item 
In both $\vec{Y}(0)$ and $\vec{Z}(0)$, place within $[0,1-Q)$ the respective \emph{unmatched} subintervals from the respective partitions, once again ordered by their size such that the longest matched subinterval is on the left. 
\end{enumerate}
\item 
Sample $U,U^\prime \sim \Unif{[0,1)}$ and use these random variables to evolve the partitions $\vec{Y}(0), \vec{Z}(0)$ as follows:
\begin{enumerate}
\item 
Call the subintervals $h_1\in \vec{Y}(0)$ or $h_2 \in \vec{Z}(0)$ \emph{highlighted} if $U$ falls into these subintervals after the reordering described above.
\item 
If $U^\prime$ falls into a subinterval $g_1 \in \vec{Y}(0)$ such that $g_1 \neq h_1$, then merge $h_1$ and $g_1$. Do the same for $h_2$ and $\vec{Z}(0)$.
\item  
If $U^\prime$ falls into a highlighted subinterval $h_1 \in \vec{Y}(0)$, then split $h_1$ at $U^\prime$. Do the same for $h_2$ and $\vec{Z}(0)$.
\end{enumerate}
\item 
Call $\vec{Y}(1)$ and $\vec{Z}(1)$ the new partitions that are obtained by the reordering and the application of the dynamics, and repeat.
\end{enumerate}
See Fig.~\ref{fig:asc} for an illustration. \hfill$\spadesuit$
\end{definition}

The construction specified in Definition~\ref{def:asc} has to be modified slightly to fit the setting of CFDP. First, we need to convert a permutation into a partition of the unit interval.

\begin{definition}[Cycle structure on a dynamic permutation]
\label{def:CS}
Denote by $|\gamma^{(i)}(t)|$ the size of the $i^{\mathrm{th}}$-largest permutation cycle at time $t\geq 0$, with the following conventions: if there is only a finite number $k\in\N$ of permutation cycles, then $|\gamma^{(j)}(t)| = 0$ for all $j>k$, while if two or more elements have the same size, then their ordering (among each other) is arbitrary. Define
\begin{align}
\label{eq:frakX}
\vec{\mathfrak{X}}^{(n)}(t) = \left(\mathfrak{X}_i^{(n)}(t)\right)_{i\in\N}
= \left(\frac{|\gamma^{(i)}(t)|}{|\maxgraph{\AG{\Pi_n}(t)}{}|}\right)_{i\in\N},
\end{align}
i.e., the ordered sequence of sizes of permutation cycles normalised by the size of $\maxgraph{\AG{\Pi_n}(t)}{}$. Abusing notation, we denote by $\vec{\mathfrak{X}}^{(n)}(t)$ also the partition of the interval $[0, n/|\maxgraph{\AG{\Pi_n}(t)|}{}]$ into subintervals induced by $\vec{\mathfrak{X}}^{(n)}(t)$. \hfill$\spadesuit$
\end{definition}

\begin{remark}[Normalisation of $\vec{\mathfrak{X}}^{(n)}(t)$]\label{rem:norm}
Note that, with this normalisation, the sum of all the elements in $\vec{\mathfrak{X}}^{(n)}(t)$ is $n/|\maxgraph{\AG{\Pi_n}(t)}{}|$, which is roughly $\tfrac{1}{\zeta(t/n)}$, and the permutation cycles with elements on the largest component of the associated graph process (which is unique with high probability) correspond approximately to the subinterval $[0,1]$ of the partition induced by $\vec{\mathfrak{X}}^{(n)}(t)$.
\end{remark}

All the details required to modify Schramm's coupling to the setting of CFDP are explained in \cite[Sections 5.2 and 5.3]{BKLM2019}, and we will not reproduce them here in full. The main modifications of the coupling, this time between $\vec{\mathfrak{X}}^{(n)}(t)$ seen as a partition of $[0, n/|\maxgraph{\AG{\Pi_n}(t)}{}|]$ and a random partition of the interval $[0,1]$ distributed according to $\mathsf{PoiDir}(1)$, are as follows:
\begin{enumerate}
\item 
Allow for \emph{approximate} matching of components, i.e., allows for a margin of order $O(n^{-1/2})$.
\item 
Markers $U, U^\prime$ that generate the dynamics are sampled uniformly from $[0, n/|\maxgraph{\AG{\Pi_n}(t)}{}|]$, but if one of them falls outside the interval $[0,1]$, then the move is \emph{not carried out} in the initially $\mathsf{PoiDir}(1)$-distributed partition.
\item 
A \emph{forbidden} set $F(t)$ is the set of points that $U, U^\prime$ must avoid for the coupling to be successful. This set takes care of the possible errors that may arise due to, for example, the discrete nature of the permutations or the growing size of the largest component in the associated graph process.
\end{enumerate}

\begin{remark}[Limiting distribution of the cycle structure $\vec{\mathfrak{X}}^{(n)}$]
In \cite[Theorem 1.1]{S2005}, Schramm's coupling in the form adapted to CFDP was used to show that, at any time $t\geq cn$ with $c>\tfrac12$, the restriction of $\vec{\mathfrak{X}}^{(n)}(t)$ to the interval $[0,1]$ has the distributional limit
\begin{align}
\vec{\mathfrak{X}}_{\mid [0,1]}^{(n)}(t) \convdist \mathsf{PoiDir}(1), 
\end{align}
where $\mathsf{PoiDir}(1)$ is the Poisson-Dirichlet distribution with parameter $1$, which is the \emph{unique} invariant distribution of $\vec{\mathfrak{X}}^{(n)}$ w.r.t.\ the permutation dynamics \cite{S2005}.
\end{remark}

Note that the largest cycle $\mathfrak{X}_1^{(n)}(cn)$, for $c>1/2$, is typically large (see e.g. \cite{H2001}). In the lemma stated below, which is adapted from \cite[Section 5.3]{BKLM2019}, we collect results that give us control over $\vec{\mathfrak{X}}^{(n)}(t)$. To understand the limitations in the statement of this lemma, we need to introduce some extra notation: 

\begin{definition}[Events relating to the success of Schramm's coupling]
\label{def:schrammfail}
Let $c>\tfrac12$, fix $\beta>0$ and set $T_\beta=\lceil \beta^{-1/2}\rceil - 1$ and $\hat{I}_{c,\beta}=[cn, cn+T_\beta]$. Consider Schramm's coupling 
$$
\big(\vec{\mathfrak{X}}^{(n)}(t), \vec{Z}(t)\big)_{t\in \hat{I}_{c,\beta}}
$$ 
with $\vec{Z}(0)\sim\mathsf{PoiDir}(1)$. For $t \in \hat{I}_{c,\beta}$, define
\begin{itemize}
\item 
$N_{\beta}(\vec{\mathfrak{X}}^{(n)}(t))$ to be the number of unmatched blocks in $\vec{\mathfrak{X}}^{(n)}(t)$ whose sizes are larger than $\beta$. Analogously for $N_{\beta}(\vec{Z}(t))$.
\item
$\overline{N}_{\beta}(t) = N_{\beta}(\vec{\mathfrak{X}}^{(n)}(t)) + N_{\beta}(\vec{Z}(t))$.
\item 
$\sigma(\beta,\vec{\mathfrak{X}}^{(n)}(t))$ to be the total length of the blocks of size smaller than $\beta$ in the partition $\vec{\mathfrak{X}}^{(n)}(t)$. Analogously for~$\vec{Z}(t)$.
\item
$\overline{\beta}(t) = \beta + \sigma(\beta, \vec{\mathfrak{X}}^{(n)}(t)) + \sigma(\beta,\vec{Z}(t))$. 
\item
$\mathfrak{X}^{(n, {\rm\sss UM})}_{1}(t)$ to be the largest unmatched segment in $\vec{\mathfrak{X}}^{(n)}(t)$.
\end{itemize}
Using the above notation, define the following events:
\begin{align}
\label{eq:Aevents}
\mathcal{A}_1(t) = \{ \overline{\beta}(t) \leq \beta^{3/4}\}, \qquad\mathcal{A}_2(t)
= \{ \overline{N}_{\beta}(t) \leq \beta^{-1/4}\}.
\end{align}
\hfill$\spadesuit$
\end{definition}

The next definition captures the key regularity event used in the upcoming arguments.

\begin{definition}[Dynamics regularity event for Schramm's coupling]
\label{def:SCreg}
Define the event\footnote{ See \cite[p.~44]{BKLM2019}, where this event is denoted by $\mathcal{G}$ and the three defining properties correspond to the three items before Eq.~(5.7) therein.} $\mathcal{A}_3(\hat{I}_{c,\beta})$ by requiring that for any time $t\in \hat{I}_{c,\beta}$ the following three facts hold:
\begin{enumerate}
\item 
$U, U^\prime$ sampled at time $t$ do not fall in the forbidden set~$F(t)$.
\item 
The dynamics does not split a component of $\vec{\mathfrak{X}}^{(n)}(t)$ of size $\leq \sqrt{n}/|\maxgraph{\AG{\Pi_n}(t)}{}|$.
\item 
If the dynamics induces a fragmentation in one of the partitions, then it does so also in the other partition. 
\end{enumerate}
\hfill$\spadesuit$
\end{definition}

Since the event $\mathcal{A}_3(\hat{I}_{c,\beta})$ is crucial for the success of Schramm's coupling, we will need the following quantitative estimate (see also \cite[eq.~(5.7)]{BKLM2019}) to get uniform control over the time scale on which Schramm's coupling remains successful with high probability.
 
\begin{lemma}[Probability bound for $\mathcal{A}_3(\hat{I}_{c,\beta})$]
\label{lem:A3bound}
Let $c>\tfrac12$ and $\beta>0$. Then
\begin{align}
\label{eq:1/26}
\P\left( \mathcal{A}_3^\cmpl(\hat{I}_{c,\beta}) \right) \leq 16 \beta^{-1}n^{-1/13}.
\end{align}
\end{lemma}

With this notation and information in hand, we can now state the following key lemma, adapted from \cite[Section~5.3]{BKLM2019}, which gives us control over the cycle structure $\vec{\mathfrak{X}}^{(n)}(t)$:

\begin{lemma}[$\mathsf{PoiDir}(1)$ approximation of  $\vec{\mathfrak{X}}^{(n)}(t)$]
\label{lem:schramm}
Fix $c>\tfrac12$ and $\beta>0$. Let $\vec{\mathfrak{X}}^{(n)}(t)$ be as in Definition~\ref{def:CS}. Consider $(\vec{Z}(t))_{t\in \hat{I}_{c,\beta}}$ such that $\vec{Z}(cn) \sim \mathsf{PoiDir}(1)$ is sampled independently of anything else, and at later times the evolution of $\vec{Z}(t)$ is governed by the dynamics of the underlying permutation (see \cite[Section 5]{BKLM2019}). Consider Schramm's coupling $(\vec{\mathfrak{X}}^{(n)}(t), \vec{Z}(t))_{t\in \hat{I}_{c,\beta}}$, and let $q \sim \Unif{[cn, \ldots, cn + T_\beta]}$ be a uniform random variable in $\hat{I}_{c,\beta}$, independent of anything else. For $\delta\in(0,1)$, define the event
\begin{equation}\label{deltaNorm}
\mathscr{D}_{\delta, \beta} = \left\{ \norm{\vec{\mathfrak{X}}^{(n)}(q) - \vec{Z}(q)}_{\infty} < \delta\right\}.
\end{equation}
Recall the events in \eqref{eq:ERtypical} and in Definition~\ref{def:schrammfail}. There exist constants $C, C' >0$ such that, for $n$ sufficiently large, the following estimate is valid uniformly in $\beta$ and $\delta$: 
\begin{align}
\label{eq:supnorm}
\P( \mathscr{D}_{\delta, \beta}^\cmpl) 
& \leq \P\left(\typicalevent{\mathcal{C}}{cn}\right) + \P\left(\typicalevent{}{cn} \cap \mathcal{A}_1^\cmpl(cn) \right) 
+\P\left(\typicalevent{}{cn} \cap \mathcal{A}_2^\cmpl(cn)\right)
+\P\left(\typicalevent{}{cn} \cap \mathcal{A}_3^\cmpl( \hat{I}_{c,\beta})\right)\\ 
& \notag \qquad + \P\left( \left\{ \mathfrak{X}^{ (n,{\rm\sss UM})}_{1}(q) 
\geq \delta/2\right\} \cap  \typicalevent{}{cn} \cap \mathcal{A}_1(cn) \cap \mathcal{A}_2(cn) 
\cap \mathcal{A}_3( \hat{I}_{c,\beta}) \right)\\ 
\label{eq:schramm:RHS}
&\leq o(1) + 3C \beta^{1/4} \log\left( \beta^{-1} \right) + 2C' \beta^{1/4} \log^2\left( \beta^{-1} \right) 
+ 16 \beta^{-1}n^{-1/13} + \frac{2 C_1}{\delta \log(\beta^{-3/4})}.
\end{align}
\end{lemma}

The proof of this lemma can be extracted from the various arguments and statements presented in \cite[Section~3]{S2005}. Alternatively, the precise bounds in \eqref{eq:supnorm} and \eqref{eq:schramm:RHS} are explicitly derived in \cite[Section 5]{BKLM2019}. In particular, the first inequality in \eqref{eq:supnorm} is obtained by applying a union bound in the second term of the first bound in the last display of the proof of \cite[Theorem 1.1, p.~45]{BKLM2019}. The quantitative bounds on the non-trivial terms in \eqref{eq:schramm:RHS} can be found in, respectively, the second to last display in the proof of \cite[Theorem 1.1, p.~45]{BKLM2019} for the terms involving $\mathcal{A}_1^\cmpl$ and $\mathcal{A}_2^\cmpl$, Lemma \ref{lem:A3bound} for the term involving $\mathcal{A}_3^\cmpl$, and \cite[Corollary 5.7]{BKLM2019} for the last term related to the maximal unmatched block $\mathfrak{X}^{(n, {\rm\sss UM})}_{1}$.

%%%

\subsection{Recurrence of large cycles}
\label{sec:recur}

We will make use of Schramm's coupling in the proof of Proposition~\ref{MaxFlux} below, which is at the very core of our argument in the proofs of Propositions~\ref{prop:epsbound} and \ref{prop:DTepsbound}. To this aim, we first use the quantitative estimate in Lemma \ref{lem:schramm} to derive the following technical lemma, stating that the Poisson-Dirichlet approximation of the underlying permutation dynamics is good with high probability over a properly chosen time interval of diverging length. 

\begin{lemma}[Pathwise approximation property of Schramm's coupling]
\label{lem:SCpath}
Fix $\varepsilon \in (0,\tfrac12)$ and $c>\tfrac12$. Take any time interval of the form $I_{c,n}= [cn, cn +a_n]$ with $\lim_{n\to\infty} a_n=\infty$ and $a_n=o(n^{1/26})$, and denote by $\smash{\mathcal{F}^{\sss \rm SC} = \left( \mathcal{F}^{\sss \rm SC}_t \right)_{t\in I_c}}$ the natural filtration of Schramm's coupling. For $a,b \in I_{c,n}$, define the event that for all $t\in[a,b]$ the cycle structure of the underlying permutation $ \vec{\mathfrak{X}}^{(n)}(t)$ is well-approximated by the process $\vec{Z}(t)$ introduced above. More precisely, define
\begin{align}
\Omega^{\sss\rm WA}(a,b,\varepsilon) = \left\{\forall\, t \in [a,b]\colon\,  \norm{\vec{\mathfrak{X}}^{(n)}(t) - \vec{Z}(t)}_{\infty} < \varepsilon \right\} \in \mathcal{F}^{\sss \rm SC}_b ,
\end{align}
\noindent pick a sequence $(d_n)_{n\in\N}$ such that $d_n = o(\log a_n)$ and $d_n\to\infty$ as $n\to\infty$. Then
\begin{align}
\P \left(\Omega^{\sss\rm WA}(q_n, q_n+ \sqrt{d_n},\varepsilon) \right) = 1 - o(1)
\end{align}
for  $q_n\sim\mathsf{Unif}(I_{c,n})$. 
\end{lemma}

\begin{proof}
We will use the estimates in Appendix~\ref{sec:summary}. In particular, with the choice of $a_n$ and $d_n$ as in the statement of the lemma, we set $\delta=\delta_n=\varepsilon/d_n,$ and $\beta=\beta_n$ such that $a_n=\lceil \beta_n^{-1/2}\rceil - 1$, and note that Lemma \ref{lem:A3bound} with $\hat{I}_{c,\beta}=I_{c,n}$ guarantees that
\begin{equation}
\label{ciao}
\P\left(\mathcal{A}_3( I_{c,n}) \right) \geq 1-o(1)
\end{equation}
as soon as $a_n=o(n^{1/26})$.
    
The key event $\mathcal{A}_3( I_{c,n})$, occurring with high probability, guarantees regularity of the dynamics uniformly over the time interval $I_{c,n}$.  On this event, we can find a (much smaller) random subinterval $[q_n, q_n+ \sqrt{d_n}\,]$, still of diverging length, for which the Poisson-Dirichlet approximation up to the threshold $\varepsilon$ is valid with high probability. Indeed, with the notation as in \eqref{deltaNorm}, we can bound
\begin{align} 
\notag
\P\big(\Omega^{\sss\rm WA}(q_n, q_n+ \sqrt{d_n},\varepsilon)\big)
&= \P\Big(\forall\, t \in [q_n,q_n+\sqrt{d_n}]\colon\,  \norm{\vec{\mathfrak{X}}^{(n)}(t) - \vec{Z}(t)}_{\infty} < \varepsilon\Big) \\ 
& \notag \geq \P\big(\mathscr{D}_{\varepsilon/d_n, \beta_n } \cap \mathcal{A}_3(I_{c,n})\big)
= \P\big( \mathscr{D}_{\varepsilon/d_n, \beta_n}\big) - \P\big( \mathscr{D}_{\varepsilon/d_n, \beta_n}  
\cap \mathcal{A}^\cmpl_3(I_{c,n}) \big) \\
& \notag \geq 1- \P\big( \mathscr{D}^\cmpl_{\varepsilon/d_n, \beta_n }\big) - o(1) \geq 1 - o(1),
\end{align}
where the first inequality is true because, on the event $\mathcal{A}_3(I_{c,n})$, if at time $q_n\in I_{c,n}$ we have $\lVert\vec{\mathfrak{X}}^{(n)}(q_n) - \vec{Z}(q_n)\rVert_{\infty} < \varepsilon/d_n$, then the next $\sqrt{d_n}$ (coagulation or fragmentation) moves of the dynamics can only increase the sup-norm in a bounded way. To see why, observe that on the event $\mathcal{A}_3(I_{c,n})$ only the following can occur:
\begin{itemize}
\item 
Due to Definition~\ref{def:SCreg}[1.], we can leave out the effect of fragmentations of tiny components. Furthermore, Definition~\ref{def:SCreg}[3.] guarantees consistency of moves between the two coupled cycle structures.
\item 
Coagulation or fragmentation within the \emph{matched} components can increase the sup-norm by at most  $\tfrac{C}{\sqrt{n}}$, which is the margin allowed by the approximate matching rule. This discrepancy can be made arbitrarily small by taking $n$ large enough. This follows from Definition~\ref{def:SCreg}[1.] and the definition of a forbidden set $F(t)$.
\item 
On the event $\mathscr{D}_{\varepsilon/d_n, \beta_n }$, the largest unmatched component satisfies $\mathfrak{X}^{\sss (n, {\rm\sss UM})}_{1}(q_n) \leq  \tfrac{\varepsilon}{2d_n}$. This guarantees that, during the interval $ [q_n,q_n+\sqrt{d_n}]$, the coagulation or fragmentation moves involving at least one \emph{unmatched} component can, in the worst case, increase the sup-norm by at most $\sum_{k=1}^{\sqrt{d_n}}\tfrac{k \varepsilon}{2d_n} \leq \varepsilon$. 
\end{itemize}
The other inequalities follow, respectively, from \eqref{ciao} and \eqref{eq:schramm:RHS}, with the choice of parameters as here.
\end{proof}

Using the pathwise approximation result in Lemma~\ref{lem:SCpath}, we can finally state the core proposition of this appendix (which was used to obtain \eqref{keybd}):

\begin{proposition}[Infinite crossings of level $1-\varepsilon$ for $\mathfrak{X}^{(n)}_1$]
\label{MaxFlux}
Fix $\varepsilon \in (0,\tfrac12)$ and $c>\tfrac12$. Take any time interval of the form $I_{c,n}= [cn, cn +a_n]$ with $\lim_{n\to\infty} a_n=\infty$ and $a_n=o(n^{1/26})$. Consider the event ${\cal E}_n(c,\varepsilon, \kappa)$ that the process $(\mathfrak{X}_1^{(n)}(t)-(1-\varepsilon))_{t = cn}^{cn+a_n}$ changes sign $\kappa \in \N$ times along an increasing sequence of times $(t_k)_{k=1}^{\kappa}$, i.e.,
\begin{equation}
\label{KeyEvent}
{\cal E}_n(c,\varepsilon, \kappa) = \big\{\exists\,(t_{k})_{k=1}^{\kappa}\colon\, \mathfrak{X}_1^{(n)}(t_{k}-1)< 1-\varepsilon,\,\mathfrak{X}_1^{(n)}(t_{k})
\geq 1-\varepsilon\,\big\}.
\end{equation}
For any $\kappa\in\N$,
\begin{equation}
\label{Ewhp}
\lim_{n\to\infty} \P({\cal E}_n(c,\varepsilon, \kappa)) = 1.
\end{equation}
\end{proposition}

\begin{proof}
The idea behind the proof is that if the coupling described above holds, then the behaviour of the $\mathsf{PoiDir}(1)$-sample will induce the desired large cycles in the dynamic permutation.

Let $I_{c,n} = [cn, cn+a_n]$, and denote by $\smash{\mathcal{F}^{\sss \rm SC} = (\mathcal{F}^{\sss \rm SC}_t)_{t\in I_c}}$ the natural filtration of Schramm's coupling. For any $a,b \in I_{c,n}$, $\varepsilon \in (0,\tfrac12)$ and $\kappa\in\N$, define the event that large cycles in the $Z$ process \emph{recur} more that $\kappa$ times:
\begin{align}
\Omega^{\sss\rm ZR}(a,b,\varepsilon, \kappa) = \big\{\#\{ t \in[a,b]\colon\, Z_1(t) > 1-\varepsilon\} > \kappa \big\}  \in \mathcal{F}^{\sss \rm SC}_b.
\end{align}
The key observation is that, for any $0 < \lambda < \varepsilon$,
\begin{align}
\label{eq:indevent}
\big\{ \#\{ t \in[a,b]\colon\, \mathfrak{X}_1^{(n)}(t) > 1-\varepsilon\} > \kappa \big\}  \supseteq 
\Omega^{\sss \rm WA}(a,b,\lambda) \cap \Omega^{\sss \rm ZR}(a,b,\varepsilon+\lambda, \kappa).
\end{align}
This is true because if the supremum norm between $\mathfrak{X}^{(n)}(t)$ and $Z(t)$ is at all times bounded by $\lambda$, then every occurrence of a cycle larger than $1-\varepsilon+\lambda$ in $Z(t)$ induces a cycle of size at least $1-\varepsilon$ in $\mathfrak{X}^{(n)}(t)$. 

From Lemma~\ref{lem:SCpath}, it follows that for $q_n\sim\mathsf{Unif}(I_{c,n})$ with $d_n\to\infty$ and $d_n = o(\log a_n)$, the well-approximation event $\Omega^{\sss \rm WA}(q_n,q_n+\sqrt{d_n},\varepsilon/2)$ occurs with high probability. Recall that, at the beginning of the coupling, $Z(cn) \sim \mathsf{PoiDir}(1)$ independently of everything else, and the dynamics of the coupling is such that $\mathsf{PoiDir}(1)$ is invariant. Denote by $\mathcal{P}_1$ the space of size-ordered countable partitions of the unit interval $[0,1]$ into subintervals, which is the space over which the measure $\nu = \mathsf{PoiDir}(1)$ is defined. Introduce the notation
\begin{align}
L_{\varepsilon} = \left\{ P \in \mathcal{P}_1\colon\, P_1 > 1-\varepsilon \right\},
\end{align}
where $P_1$ denotes the first, and therefore the largest, element of $P$. By \cite[Theorem 1.2]{MP2020}, for any $\varepsilon \in (0,\tfrac12)$,
\begin{align}
\nu\left( L_{\varepsilon} \right) = -\log(1-\varepsilon) > 0.
\end{align}
Recall from \cite{S2005} that the evolution of $Z$ is a time-homogeneous Markov process and that $\nu = \mathsf{PoiDir}(1)$ is the unique invariant measure of this process. Since, for any $\varepsilon\in(0,\tfrac12)$, the set $L_{\varepsilon}$ has a strictly positive $\nu$-measure and the starting point of the process is sampled according to $\nu$, it follows that the \emph{hitting} time of $L_{\varepsilon}$ is finite a.s.\ and has finite expectation. Furthermore, the \emph{return} times to $L_{\varepsilon}$ are all finite a.s.\ and, by the \emph{Kac recurrence time lemma} \cite[Theorem~4.6]{P1989}, have finite expectation, equal to $1/\nu(L_{\varepsilon})$. Consequently, for any fixed number of occurrences $\kappa\in \N$, any fixed threshold for cycle sizes $\varepsilon \in (0,\tfrac12)$, random times $q_n\sim\mathsf{Unif}(I_{c,n})$ independent of anything else, and any $d_n\to\infty$, 
\begin{align}
\lim_{n\to\infty} \P\big(\Omega^{\sss \rm ZR}(q_n, q_n+d_n, \varepsilon, \kappa)\big) = 1.
\end{align}
Therefore, for $q_n\sim\mathsf{Unif}(I_{c,n})$, any $\kappa\in\N$ and any sequence $(d_n)_{n\in\N}$ such that $d_n = o(\log a_n)$ and $d_n\to\infty$ as $n\to\infty$,  
\begin{align}
\P \Big( \Omega^{\sss \rm WA}(q_n, q_n+\sqrt{d_n},\varepsilon) \cap \Omega^{\sss \rm ZR}(q_n, q_n+\sqrt{d_n},\tfrac\varepsilon2, \kappa) \Big) = 1-o(1),
\end{align}
which by \eqref{eq:indevent} yields the desired result.
\end{proof}

\begin{remark}[Extension of Proposition~\ref{MaxFlux} to the degree-two graphs with rewiring]
\label{CMMaxFlux}
In  \cite{GUW2011} and \cite[Section 5.2]{BKLM2019} it is noted that Schramm's coupling can be adapted to the setting of coagulation-fragmentation dynamics that keep the measure $\mathsf{PoiDir(\theta)}$, $\theta\in(0,1]$, invariant. An example is a dynamic graph model with all degrees equal to two and endowed with a rewiring dynamics, which corresponds to a coagulation-fragmentation dynamics with invariant measure $\mathsf{PoiDir}(1/2)$ (see \cite{O2021}). Since $\nu_{\theta}(L_{\varepsilon}) > 0$ for any $\theta,\varepsilon\in(0,1)$ and $\nu_{\theta} = \mathsf{PoiDir}(\theta)$, the proof of Proposition~\ref{MaxFlux} can be adapted to the aforementioned situation.
\end{remark}

In conclusion of this section, we state a corollary of the previous proposition, which is useful in Section~\ref{sec:CFDP}:

\begin{corollary}\label{cor:MaxFlux}
In the setting of Proposition~\ref{MaxFlux},
\begin{align}
 \P \left( {\cal E}^\cmpl_n(\tfrac{\Mtddalt}{n},\varepsilon, \kappa) \cap \{ \Mtddalt \geq (\tfrac12 + \delta)n \}\right) = o(1).
\end{align}
\end{corollary}

\begin{proof}
First, let us rewrite the desired expression as
\begin{align}
&\P \left( {\cal E}^\cmpl_n(\tfrac{\Mtddalt}{n},\varepsilon, \kappa) \cap \{ \Mtddalt \geq (\tfrac12 + \delta)n \}\right)  \\ &\qquad\qquad = \E\left[ \P \left( {\cal E}^\cmpl_n(\tfrac{\Mtddalt}{n},\varepsilon, \kappa) \mid \Mtddalt \right) \1_{\{ \Mtddalt \geq (\tfrac12 + \delta)n \}} \right].\nonumber
\end{align}
Since the inner part of the right-hand side is bounded between 0 and 1, by dominated convergence theorem it is enough to show that 
\begin{align}
\label{eq:insideE}
\P \left( {\cal E}^\cmpl_n(\tfrac{\Mtddalt}{n},\varepsilon, \kappa) \mid \Mtddalt \right) \1_{\{ \Mtddalt \geq (\tfrac12 + \delta)n \}} = o_{\prob}(1).
\end{align}

The core observation is that arguments based on Schramm coupling from Proposition~\ref{MaxFlux} work even under the conditioning in \eqref{eq:insideE}. From the perspective of the cycle structure, on the event $\{ \Mtddalt \geq (\tfrac12 + \delta)n \}$, at time $\Mtddalt$ the giant component of the associated graph process is enlarged by a connected component whose size can be, with high probability, uniformly bounded by $C\log^2 n$, for some $C>0$. But the effect of these moves is already captured by the forbidden set construction (recall Definition~\ref{def:SCreg} and the discussion below Remark~\ref{rem:norm}), and a quantitative bound on the effect of these moves is expressed in Lemma~\ref{lem:A3bound}. Therefore, by repeating the arguments from Lemma~\ref{lem:SCpath} and Proposition~\ref{MaxFlux}, we obtain \eqref{eq:insideE}, which yields the desired expression.

\end{proof}

%%%%%%% APPENDIX D %%%%%%%%%%%%%%%%%%%%%%%%%%%%%%%

\section{Mixing upon dropdown on the largest cycle}
\label{app:sc}

Recall that for two laws $\mu,\nu$ defined on the same countable probability space 
\begin{align}
\label{eq:L1}
d_{\rm TV}(\mu, \nu) = \frac12 \sum \limits_{v} |\mu_v - \nu_v|.
\end{align}
Fix $n\in\N$. Recall the short-hand notation $M=|\maxgraph{\AG{\Pi_n}(\Mtddalt)}{}|$. We wish to show that, on the events (recall \eqref{eq:Mevents})
\begin{equation}
\label{eq:Ms}
\begin{aligned}
\mathcal{M}_1(\varepsilon, \delta) &= \big\{|\supp(\ISRWdistr{}(\Mtddalt )) | > \varepsilon M\big\} \cap \Omega^{\sss\rm (SC)}(\Mtddalt) \cap \typicalevent{\rm\sss (ER)}{\Mtddalt},\\
\mathcal{M}_2(\varepsilon) &= \big\{\exists\, t_{L} \in (\Mtddalt, \Mtddalt + a_n)\colon\,  \mathfrak{X}_1^{(n)}(t_L) > 1-\varepsilon^2 \big\},
\end{aligned}
\end{equation}
the following estimate (recall  \eqref{TV-bound-onM2}) is valid:
\begin{align}
\label{eq:apxE:target}
\TVD\Big(\ISRWdistr{}(t_L), \Unif{\maxgraph{\AG{\Pi_n}(t_L)}{}}\Big) < \varepsilon.
\end{align}

On the event $\mathcal{M}_1(\varepsilon, \delta)$, the ISRW-distributiuon at time $\Mtddalt$ is uniform over a single cycle supported on the largest component of the associated graph process (due to $\Omega^{\sss\rm (SC)}(\Mtddalt)$), which by assumption is larger than $\varepsilon M$ (due to $\smash{\{|\supp(\ISRWdistr{}(\Mtddalt))| > \varepsilon M\}}$). Therefore,
\begin{align}
\label{eq:distbound}
\ISRWdistr{v}(\Mtddalt)  
\begin{cases}
\leq \frac{1}{\varepsilon M}, &v \in \supp(\ISRWdistr{}(\Mtddalt)),\\
= 0, & \text{otherwise}.
\end{cases} 
\end{align}
We will use \eqref{eq:distbound} to provide us with an upper bound on the mass carried by the elements on the giant component of the associated graph process.

\begin{figure}[h]
\centering
\includegraphics[width=0.5\textwidth]{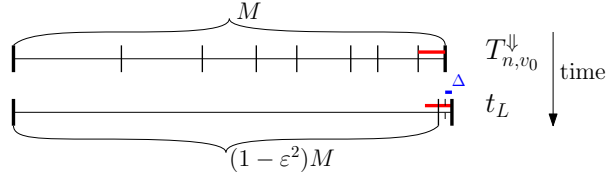}
\caption{The first line segment represents the cycle structure of the permutation restricted to the giant component of the associated graph process at time $\Mtddalt$. The red line represents the size of the single cycle that contains the full mass of the ISRW-distribution. The second line segment represents the same at time $t_L$: a cycle of size $(1-\varepsilon^2)M$ has appeared, which necessarily carries some of the mass of the ISRW-distribution.}
\label{fig:apxE:mix}
\end{figure}

On the event $\mathcal{M}_1(\varepsilon, \delta) \cap \mathcal{M}_2(\varepsilon)$, the largest cycle at time $t_L$ necessarily carries some of the mass of the ISRW-distribution, since the initial size of the support is larger than the size of the subset not covered by the largest cycle (see Fig.~\ref{fig:apxE:mix} for a visual explanation).

%%%%%%%%%%%%%%%%%%%%%%%%%%%%%%%%%%%%%%%%%%%%%%%%%%
%%%%%%%%%%%%%%%%%%%%%%%%%%%%%%%%%%%%%%%%%%%%%%%%%%%

We proceed by computing the worst-case $l^1$-distance between the two distributions at time $t_L$, which will yield the desired inequality in \eqref{eq:apxE:target}. Put $\Delta = |\maxgraph{\AG{\Pi_n}(t_L)}{}| - M$, which represents the growth of the largest component of the associated graph process compared to its size at time $\Mtddalt$. Note that the uniform distribution in \eqref{eq:apxE:target} gives the individual elements the following mass:
\begin{align}
\label{eq:unifdistr}
\Unif{\maxgraph{\AG{\Pi_n}(t_L)}{}}(v) 
= \begin{cases}
\frac{1}{M + \Delta}, &v\in \maxgraph{\AG{\Pi_n}(t_L)}{},\\
0, &\text{otherwise}.
\end{cases}
\end{align}
To compute the worst-case $l^1$-distance at time $t_L$, we maximise the individual summands in \eqref{eq:L1}. To do so, we work with a distribution $\distr{\dagger}{}$ supported on the giant component of the associated graph process that puts mass $1/(\varepsilon M)$ (see \eqref{eq:Ms}) on the elements outside the largest cycle of size $(1-\varepsilon^2)M$, and spreads the remaining mass $\Sigma = 1 - \frac{\varepsilon^2 M+\Delta}{\varepsilon M}$ over the largest cycle. In symbols,
\begin{align}
\label{eq:dagger}
\distr{\dagger}{v} 
= \begin{cases}
\frac{\Sigma}{(1-\varepsilon^2)M}, &v \in  \mathfrak{X}_1^{(n)}(t_L),\\
0, &v \notin \maxgraph{\AG{\Pi_n}(t_L)}{},\\
\frac{1}{\varepsilon M}, &\text{otherwise}.
\end{cases}
\end{align}
The measure $\distr{\dagger}{}$ constructed in this way gives the bound
\begin{align}
\TVD\Big(\ISRWdistr{}(t_L),\Unif{\maxgraph{\AG{\Pi_n}(t_L)}{}}\Big) 
\leq \TVD\Big(\distr{\dagger}{},\Unif{\maxgraph{\AG{\Pi_n}(t_L)}{}}\Big).
\end{align}

Finally, we compute, using \eqref{eq:L1} and \eqref{eq:unifdistr}--\eqref{eq:dagger},
\begin{equation}
\begin{aligned}
&\TVD\Big(\distr{\dagger}{},\Unif{\maxgraph{\AG{\Pi_n}(t_L)}{}}\Big) = \frac12 \norm{ \Unif{\maxgraph{\AG{\Pi_n}(t_L)}{}} - \distr{\dagger}{}}_{l^1}\\ 
&\qquad = \frac12 (\varepsilon^2 M + \Delta) \left[ \frac{1}{\varepsilon M} - \frac{1}{M+\Delta} \right]  
+ \frac12 (1-\varepsilon^2)M \left[ \frac{1}{M+\Delta} - \frac{\Sigma}{(1-\varepsilon^2)M}\right]\\ 
&\qquad = \varepsilon + \frac{\Delta}{\varepsilon M} - \frac{\varepsilon^2 M}{M + \Delta} - \frac{\Delta}{M+\Delta}.
\end{aligned}
\end{equation}
On the event $\typicalevent{\rm\sss (ER)}{\Mtddalt}$, the sizes of all but the largest component are uniformly bounded by $Cn^{2/3}$ for some $C>0$. Therefore, for $a_n$ growing slowly enough,  more concretely $a_n = o(n^{1/3})$, it follows that $\Delta = o(M) = o(n)$. Therefore
\begin{equation}
\TVD\Big(\ISRWdistr{}(t_L), \Unif{\maxgraph{\AG{\Pi_n}(t_L)}{}}\Big) \leq \varepsilon - \varepsilon^2 + o(1) < \varepsilon,
\end{equation}
where the last inequality is true for $n$ large enough.

%%%%%%%% APPENDIX E %%%%%%%%%%%%%%%%%%%
 
\section{Mixing in dynamic degree-two graphs}
\label{apx:deg2}

\subsection{Permutations and degree-two graphs}

Let $\pi \in S_n$ be a permutation of $[n]=\{1, \ldots, n\}$. Such a permutation admits a decomposition into distinct permutation \emph{cycles}, and this decomposition can be used to create a mapping between permutations and graphs whose vertex degrees are all equal to two. Note that this mapping is \emph{not} a bijection between the set of all degree-two graphs on $n$ vertices $\mathfrak{G}_n^{\rm deg2}$ and the set of all permutations of $n$ elements $S_n$. Such a bijection is not possible in general: while the connected components of degree-two graphs indeed are  cycles, permutation cycles also carry information about the direction of traversal within them.  

Let us nonetheless construct two mappings between these sets. There is a natural bijection between $[n]$ seen as the set of permutation \emph{elements} and the same set seen as the set of graph \emph{vertices}. Thanks to this bijection, we will (abusing notation) make no distinction between the two sets. Moreover, given a permutation $\pi \in S_n$, we can construct a degree-two graph with $n$ vertices as follows:
\begin{enumerate}
\item Set $e$ to be the least element of $\pi$, i.e., the element \enquote{1}.
\item Add an edge between the vertices $e$ and $\pi(e)$.
\item Set $e$ to be $\pi(e)$,  and repeat step 2 until the entire permutation cycle containing $e$ is traversed. 
\item Once the entire permutation cycle is traversed, set $e$ to be the least non-traversed element of $\pi$ and repeat from step $2$ onwards. If there is no such element, then the entire permutation is traversed and the algorithm terminates.
\end{enumerate}

Note that the mapping induced by this algorithm is a surjection, and hence the relation \enquote{being represented by the same degree-two graph} is an equivalence relation on $S_n$. Equivalence classes are formed by permutations whose cycle decomposition differs only in the reversal of the cyclic order between some of the cycles. For example, take $\pi_1, \pi_2 \in S_3$ such that
\begin{equation}
\pi_1 = (1,2,3), \qquad \pi_2 = (3,2,1).
\end{equation}
By following the algorithm described above, we see that these two permutations are represented by the same degree-two graph.

When constructing the mapping from $\mathfrak{G}_n^{\rm deg2}$ to $S_n$, the problem of cycle orientation manifests itself again. Since all the connected components of the degree-two graph are cycles, it is trivial to represent them as permutation cycles, but we are free to choose the direction of traversal of these connected components. These discrepancies are not relevant, since none of the mathematical objects used in the present paper depends on the direction of traversal of the permutation cycles. The most important property for our results is the size of the permutation cycles, or equivalently the size of the connected components. 

\subsection{ISRW mixing on degree-two graphs}

To connect the result of the present paper to our previous work describing mixing of random walks on top of configuration models endowed with \emph{rewiring} dynamics (see \cite{AGHH2018,AGHH20182,AGHHN2020}), we note that rewiring of degree-two graphs and CFDP can easily be related to each other (recall Fig.~\ref{FigPerRew}). Our techniques and results can be adapted to the setting where the underlying geometry is modelled by a graph process starting from the configuration where \emph{all the vertices have a self-loop}, equipped with rewiring dynamics. This process has been previously studied in \cite{O2021}. 

To state these results, let us first define the process and the underlying geometry. The following definition are an adaptation of Definitions~\ref{def:cfdp} and \ref{def:ISRW}:

\begin{definition}[ISRW on a graph]
\label{def:graphISRW}
Take a sequence of graphs $(G_n(t))_{t\in\N_0}$ and a vertex $\startvertex \in [n]$. Denote by  $\conncomp{G_n(t), v}$ the connected component of the graph $G_n(t)$ that contains the vertex~$v$. Formally, the infinite-speed random walk (ISRW) starting from $\startvertex$ is a sequence of probability distributions $(\ISRWdistr{}(t))_{t\in\N_0}$ supported on $[n]$, with initial distribution at time $t=0$ given by
\begin{equation}
\ISRWdistr{}(0) = \left(\ISRWdistr{w}(0)\right)_{w\in[n]},
\end{equation}
where $\ISRWdistr{w}(0)$, the mass at $w \in [n]$ at time $t=0$, is given by
\begin{equation}
\ISRWdistr{w}(0) = 
\begin{cases}
\frac{1}{|\conncomp{G_n(0), v} |},
&w \in \conncomp{G_n(0), v},\\ 
0, 
&w \notin  \conncomp{G_n(0), v},
\end{cases}
\end{equation}
and with distribution at later time $t \in \N$ given by
\begin{equation}
\ISRWdistr{}(t) = \left(\ISRWdistr{w}(t)\right)_{w\in[n]},
\end{equation}
where
\begin{equation} 
\ISRWdistr{w}(t) = \frac{1}{|\conncomp{G_n(t), w}|} \sum_{u \in \conncomp{G_n(t), w}}\ISRWdistr{u}(t-1).
\end{equation}
Informally, ISRW \emph{spreads infinitely fast over the connected component it resides on}. \hfill$\spadesuit$ 
\end{definition}

\begin{definition}[Dynamic degree-two graph with rewiring dynamics]
Fix the vertex set $\vertex = [n]$ and let $G_n(0)$ be the graph with all vertices having a single self-loop. At any later time $t\in\N$, $G_n(t)$ is obtained from $G_n(t-1)$ as follows:
\begin{enumerate}
\item Pick two edges $e_1, e_2$ uniformly at random without replacement from the set of edges within $G_n(t-1)$.
\item Break  edge $e_1$ into half-edges $h_1^A, h_1^B$ and edge $e_2$ into half-edges $h_2^A, h_2^B$ (see \cite[Section~7.2]{vdH2016} for a definition of a half-edge).
\item Graph $G_n(t)$ has all the unbroken edges of $G_n(t-1)$, while the broken edges $e_1, e_2$ are replaced by a uniform choice from 2 possible sets of new edges: $\{(h_1^A, h_2^A), (h_1^B, h_2^B) \}, \{(h_1^A, h_2^B), (h_1^B, h_2^A)\}$.
\end{enumerate}
We call the sequence $(G_n(t))_{t\in\N_0}$ the \emph{dynamic degree-two graph with rewiring dynamics.}

\end{definition}

In this setting, we state the following theorem, analogous to Theorem~\ref{thm:frag:main}: 

\begin{theorem}[Mixing profile for ISRW on dynamic degree-two graphs with rewiring]
\label{thm:cm2:main} 
$\mbox{}$ 
\begin{enumerate}
\item[{\rm (1)}] 
Uniformly in $\startvertex\in[n]$,
\begin{align}
\frac{\tddalt}{n} \convdist u^\Downarrow,
\end{align}
where $u^\Downarrow$ is the non-negative random variable with distribution (recall Definition~\ref{def:ERstruct}(1))
\begin{equation}
\P(u^\Downarrow \leq u) =  \zeta(u), \qquad u\in[0,\infty).
\end{equation}
\item[{\rm (2)}] 
Uniformly in $\startvertex\in[n]$, 
\begin{align}
(\dnTVD{un})_{u \in [0,\infty)} \convdist \big(1-\zeta(u)\1_{\{u>u^\Downarrow\}}\big)_{u \in [0,\infty)}
\quad \mbox{in the Skorokhod $M_1$-topology}.
\end{align}
\end{enumerate}
\end{theorem}

\begin{proof}
Since Schramm's coupling and related arguments are fully applicable in this setting (see Remark~\ref{CMMaxFlux}), the arguments from Section~\ref{sec:CFDP} can be used again with only slight modifications, as outlined below:
\begin{enumerate}
\item 
\textbf{Associated graph process and drop-down time:} The two constructions carry over without modifications. Even if a cycle-creating edge in the associated graph process produces a fragmentation only with probability $1/2$ (since such an edge can with equal probability tear a segment of a cycle apart and create a new cycle or put the segment back into its original cycle in reversed order), this does not influence the arguments that underlie Lemma~\ref{lem:cfdp:tdddistr}.
\item 
\textbf{Fast mixing upon drop-down:} Since the proof of Proposition \ref{prop:DTepsbound} is based on Proposition~\ref{MaxFlux}, which can be adapted to the alternative setting (see Remark~\ref{CMMaxFlux}), we conclude that Proposition \ref{prop:DTepsbound} also carries over.
\item 
\textbf{Drop-down in a single cycle:} Lemma~\ref{lem:1cycle} requires no adaptations.
\item 
\textbf{Mixing profile:} Since all the ingredients used in the proof of Lemma~\ref{lem:CFDP:profile} and Theorem~\ref{thm:frag:main} carry over, the proofs themselves do likewise.
\end{enumerate}
\end{proof}

\begin{remark}[ISRW on a dynamic degree-two configuration model]
The setting of Theorem~\ref{thm:cm2:main} is \emph{different} from a degree-two configuration model with rewiring, where the edges in the initial graph would be created by a random matching of half-edges. We conjecture that if the underlying geometry were modelled by a degree-two configuration model with  rewiring,  then ISRW would mix in $o_{\P}(n)$ steps and the mixing profile at scale $n$ would be trivial.
\end{remark}

%%%%%%%%%%%%%% REFERENCES %%%%%%%%%%%%%%%%%%%%%%%

\bibliographystyle{abbrv}
\bibliography{references}{}

%%%%%%%%%%%%%%%%%%%%%%%%%%%%%%%%%%%%%%%%%%%%%

\end{document}